\documentclass{article}
\usepackage{makeidx}
\usepackage{amssymb}
\usepackage{amsfonts}
\usepackage{graphicx}
\usepackage{amsmath}
\usepackage{algpseudocode}
\usepackage{algorithmicx}
\usepackage{algorithm}
\usepackage{authblk}
\usepackage{blindtext}
\usepackage[colorlinks=true]{hyperref}
\usepackage{hyperref}
\usepackage{gensymb}
\usepackage{mathptmx}      % use Times fonts if available on your TeX system

\def\IRT{\textsc{IR Tools}}
\def\pen{\mathrm{\Omega}}

\begin{document}

\title{\IRT:\ A MATLAB Package of Iterative Regularization Methods and
Large-Scale Test Problems\thanks{We acknowledge funding from
Advanced Grant No.\ 291405 from the European Research Council and
US National Science Foundation under grant no.\ DMS-1522760.}}

%\author{%
%  Silvia Gazzola\footremember{UoB}{Department of Mathematical Sciences,
%University of Bath, Bath BA2 7AY, UK. }%, \email{s.gazzola@bath.ac.uk}}%
%\and Per Christian Hansen
%\footremember{DTU}{Department of Applied Mathematics and Computer Science,
%Technical University of Denmark, 2800 Kgs.\ Lyngby, Denmark. }
%%\email{pcha@dtu.dk}}%
%\and James G. Nagy\footrecall{Emory}
%\footnote{Department of Mathematics and Computer Science,
%Emory University, Atlanta, USA. }%, \email{jnagy@emory.edu}}%
%}

\author{
  {\Large Silvia Gazzola} \\ \vspace{-8pt} Department of Mathematical Sciences \\
  University of Bath, Bath BA2 7AY, UK \\ 
  Email:  \texttt{s.gazzola@bath.ac.uk}
\and  {\Large Per Christian Hansen}
\\ \vspace{-8pt} Department of Applied Mathematics and Computer Science\\
Technical University of Denmark, 2800 Kgs.\ Lyngby, Denmark \\ 
  Email:  \texttt{pcha@dtu.dk}
\and {\Large James G. Nagy}
\\ \vspace{-8pt} Department of Mathematics and Computer Science\\ Emory University, Atlanta, USA \\ 
  Email:  \texttt{jnagy@emory.edu}

}


\date{}

\maketitle

\begin{abstract}
This paper describes a new MATLAB software package of iterative regularization
methods and test problems for large-scale linear inverse problems.
The software package, called \IRT, serves two related purposes:\
we provide implementations of a range of iterative solvers,
including several recently proposed methods that are not
available elsewhere, and we provide a set of large-scale test
problems in the form of discretizations of 2D linear inverse problems.
The solvers include iterative regularization methods where the
regularization is due to the semi-convergence of the iterations,
Tikhonov-type formulations where the regularization is explicitly
formulated in the form of a regularization term, and methods that can impose
bound constraints on the computed solutions.
All the iterative methods are implemented in a very flexible fashion
that allows the problem's coefficient matrix to be available as
a (sparse) matrix, a function handle, or an object.
The most basic call to all of the various iterative methods requires
only this matrix and the right hand side vector; if the method
uses any special stopping criteria, regularization parameters, etc.,
then default values are set automatically by the code. Moreover,
through the use of an optional input structure,
the user can also have full control of any of the algorithm parameters.
The test problems represent realistic large-scale problems found
in image reconstruction and several other applications.  Numerical
examples illustrate the various algorithms and test problems
available in this package.
\end{abstract}

\section{Introduction}
\label{sec:intro}

In this paper we are concerned with discretizations of linear inverse problems of
the form
  \begin{equation}
  \label{eq:Axb}
    A\, x \approx b, \qquad A \in \mathbb{R}^{M \times N} ,
  \end{equation}
where the vector $b$ represents measured data (typically with noise)
and the matrix $A$ represents the forward mapping.
There are no restrictions on $M$ and $N$.
Given $A$ and $b$,
the aim is to compute an approximation of the unknown vector $x$.
We are concerned with large-scale problems, where $A$ is either
represented by a sparse matrix, or is given in some other form
(i.e., a user-defined object or a function handle)
in which matrix-vector products with $A$, and also possibly
$A^T$, can be performed efficiently.
Such problems arise, e.g., in computed tomography \cite{Buzug},
image deblurring \cite{ChKN}, and geoscience \cite{Zhdanov}.

Although the iterative methods described in this paper can be used
for any large-scale linear system, we are mainly interested in
problems that are ill-posed in the sense that the singular values
of $A$ gradually decay and cluster at zero.
The decay rate depends on the problem, and many large-scale problems
tend to have a rather slow decay -- however, due to the large problem
dimensions the matrix is very ill conditioned and hence the
computed $x$ is very sensitive to errors in~$b$.
Regularization is therefore needed in order to produce stable
solutions to (\ref{eq:Axb}).

Regularization is often achieved by solving a
penalized least-squares problem of the form
  \begin{equation}
  \label{eq:penalized}
    \min_x \left\{ \| A\, x - b \|_2^2 + \lambda^2\, \pen(x) \right\} ,
  \end{equation}
where the penalty term $\pen(x)$ is chosen to reflect the
specific type of regularization that is suited for the problem.
In the case where $\pen(x) = \| x \|_2^2$ and $\pen(x) = \| L \, x \|_2^2$
we obtain the classical Tikhonov regularization problem.
A different way to achieve regularization is to apply an iterative
method directly on the fit-to-data term (e.g., $\min \| A\, x - b \|_2^2$),
and terminate the iterations when semi-convergence is achieved;
that is, terminate when a desired approximation is
obtained, but before noise starts to show up in the solution.
Using an iterative method in this way is often referred to as iterative regularization.
For more details on these issues see, e.g., \cite{Hansen} and \cite{Vogel}.

As the computational problems associated with (\ref{eq:Axb}) become large,
it is crucial to formulate the forward computation -- represented by
$A$ -- in a convenient and storage-efficient way.
For example, problems in various types of computed tomography applications
typically lead to sparse matrices.
For other problems, such as image deblurring and inverse diffusion, it is most
convenient to formulate the forward problem -- and possibly its adjoint --
as computations performed by a function (in MATLAB via a
function handle or an object).
Our package allows all these representations of~$A$, thus making
it suitable for many large-scale problems.

The software is distributed as a compressed archive; uncompressing the file
will create a directory that contains the code.
More information can be found in the \texttt{README.txt} file contained
in the package.
The software is available from Netlib
\url{http://www.netlib.org/numeralgo/}
as the \texttt{na49} package.
Maintenance of the code is available from GitHub:\
\url{https://github.com/jnagy1/IRtools}.
To obtain full functionality it is recommended to also install the
MATLAB package \textsc{AIR Tools II} \cite{AIRTools}
available from Netlib as the \texttt{na47} package.

This package has two significant aims:
The first {one} is to provide model implementations
of a range of iterative algorithms that can be used for large-scale
ill-posed linear inverse problems, including several recently
proposed methods that are not available elsewhere.
The second aim is to provide a set of new
test problems for large-scale linear inverse problems that can be used
to experiment with the iterative methods in this package, or
as benchmark test problems for newly developed algorithms.
Our software satisfies the following design objectives:
\begin{itemize}
\item
The software is easy to use:\ the installation is very simple and
there are no files to be compiled.
There is no need for commercial MATLAB toolboxes.
\item
Additional iterative methods and test problems are provided via interface
to the package \textsc{AIR Tools II} \cite{AIRTools} which implements a
number of algebraic iterative reconstruction methods.
\item
Calls to all iterative solvers and all
test-problem generators are simple, and essentially identical.
\item
Strict naming conventions are used for all functions,
such as \texttt{IR\_\_\_} for the iterative solvers
and and \texttt{PR\_\_\_}  for the test-problem generators.
\item
We include realistic 2D test problems, presented in such a way that they
require no special background knowledge of the applications from which they
arise.
\item
The functions are easy to use; default values are provided for any
parameters needed by the iterative solvers and problem generators.
\item
At the same time, the user can take full control of the functionality
by changing these parameters through an optional
 \texttt{options} input structure.
\item
Stopping rules and paradigms for choosing regularization parameters
are integrated within the iterative methods.
\item
Information about the performance of the iterative methods
is returned in an optional \texttt{Info} output structure.
\item
Visualization of the right-hand side $b$ (the data) and the approximate solution
$x$ for all test problems is done by two functions
\texttt{PRshowb} and \texttt{PRshowx}.
\item
Users can easily expand the package to include new solvers and/or
new test problems.
\end{itemize}
Other MATLAB packages are available for inverse problems, but
they can either be used only on small-scale problems,
or they focus on one specific application or type of regularization scheme
(e.g., image denoising, or tomographic reconstruction, or
$\ell_1$-regularization, or total variation).
We are not aware of other packages that fully contain the broad range of
iterative solvers in this new \IRT\ package, including several
recently proposed methods that are not available elsewhere.
The solvers include iterative regularization methods where the regularization
is due to the semi-convergence of the iterations, Tikhonov-type formulations
where the regularization is explicitly formulated in the form of a
regularization term (e.g., a 1-, or 2-norm, or total variation penalization),
and methods that can impose bound constraints on computed solutions.
Compared to our earlier software packages for regularization, we make
the following remarks:
\begin{itemize}
\item
\textsc{Regularization Tools} \cite{RegTools} does not allow $A$ to be a
function handle or an object, and was designed for small-scale
problems. In addition, the small-scale test
problems included in \textsc{Regularization Tools} are outdated and
do not represent current important applications.
\item
\textsc{Restore Tools} \cite{RestoreTools} focuses solely on
image deblurring problems, and $A$ must be a MATLAB object.
\item
\textsc{AIR Tools II} \cite{AIRTools} (a drastically expanded version of the
original \textsc{AIR Tools} package)
is primarily aimed at tomographic image reconstruction.
\end{itemize}

This paper is organized as follows.
Section~\ref{sec:Solvers} gives an overview of the iterative solvers provided in \IRT,
while Section~\ref{sec:Problems} describes the various test problems.
Examples using the solvers and test problems {available in \IRT} are given
in Section~\ref{sec:Examples},
and Section~\ref{sec:Conclusions} contains concluding remarks.

\section{Overview of the Iterative Solvers}
\label{sec:Solvers}

The overall goal for our package is to provide robust and flexible implementations
of regularization algorithms based on iterative solvers for linear problems,
in a common framework.
We do not intend to survey the details and performance of all the iterative solvers
in this paper; for full details of the algorithms we refer to the papers
listed in Table~\ref{table:solvers} below.
In our framework all calls {are of the form}
\begin{verbatim}
   [X, Info] = IR___(A, b, K, options);
\end{verbatim}
Here, \texttt{A} is the discrete forward operator, \texttt{b} is the measured data,
the vector \texttt{K} determines which iterations are stored
as columns in \texttt{X}, \texttt{options} is a structure that defines
the algorithm parameters, and \texttt{Info} is a structure containing
information about the iterations, such as residual norms,
and what stopping criterion led to the iterations being terminated.

Throughout the package we follow the convention that all
error norms and residual norms are \textit{relative}.
This means that, if the true solution $x$ is provided to the
iterative method through the {\tt options} structure (see below for
an explanation on how to do this), and $x^{(k)}$ is the $k$th
iteration vector, then in \texttt{Info.Enrm} we return
  \[
    \bigl\| x - x^{(k)} \bigr\|_2 / \| x \|_2 , \qquad k=1,2,3,\ldots
  \]
Similarly, if $b$ is the right-hand side of a least squares problem
then in \texttt{Info.Rnmr} and \texttt{Info.NE\_Rnrm} (when relevant)
we return
  \[
    \bigl\| b - A\, x^{(k)} \bigr\|_2 / \| b \|_2 \quad \hbox{and}
    \quad \bigl\| A^T ( b - A\, x^{(k)} ) \bigr\|_2 /
    \bigl\| A^T b \bigr\|_2 , \qquad k=1,2,3,\ldots
  \]

Inputs \texttt{K} and \texttt{options}, and output \texttt{Info} are optional,
so that all solvers can be used with the simple call:
\begin{verbatim}
   X = IR___(A, b);
\end{verbatim}
In this case (depending on the method), default values are used
for regularization parameters and stopping criteria, and \texttt{X}
contains the approximate solution at the final iteration.
The inclusion of the input parameter \texttt{options} has the effect of
overriding various default options, depending on the considered solver
and on the fields specified in \texttt{options}.
Moreover, if the user stores in \texttt{options} additional information
about the test problem, additional information about the behavior of
the solver can be stored in the output structure \texttt{Info};
for instance, if the true solution is stored in \texttt{options},
then the relative errors are computed at each iteration and
returned in \texttt{Info}. To determine what the possible default options for the various
test problems are, use:
\begin{verbatim}
   options = IR___('defaults')
\end{verbatim}
One can then change the default options either by directly changing a specific field,
for example,
\begin{verbatim}
   options.field_name = field_value;
\end{verbatim}
or by using the function \texttt{IRset},
\begin{verbatim}
   options = IRset(options, 'field_name', field_value);
\end{verbatim}
Note that, in the above example using {\tt IRset}, it is assumed that the structure
\texttt{options} is already defined, and only one of its field values is changed.
It is possible to change multiple field values using {\tt IRset}, for example,
\begin{verbatim}
   options = IRset(options, 'field_name1', field_value1, ...
         'field_name2', field_value2, 'field_name3', field_value3);
\end{verbatim}
It is also possible to use \texttt{IRset} without a pre-defined \texttt{options}
structure, such as
\begin{verbatim}
   options = IRset('field_name', field_value);
\end{verbatim}
In this case, all default options are used, except \texttt{field\_name}.

Our package includes some standard Krylov subspace algorithms,
as well as their hybrid versions where regularization
is applied to the problem projected in a Krylov subspace. Other algorithms are based on
flexible Krylov subspace methods, where an iteration-dependent preconditioner
is used to penalize or impose constraints on the solution; sometimes these methods are
combined with restarts. For both approaches, the regularization comes from projecting onto
the Krylov subspace (possibly combined with regularization of
the projected problem) or from applying the method to a
penalized problem of the form (\ref{eq:penalized}).
Tables \ref{table:solvers} and \ref{table:problems} give an overview
of each of the iterative solvers in the package, and some additional
discussion is provided in the following subsections.

\begin{table}
\renewcommand{\arraystretch}{1.5}
\caption{\label{table:solvers} List of iterative methods in \IRT; the
two functions \texttt{IRart} and \texttt{IRsirt} require \textsc{AIR Tools II}.
The naming convention in the \textbf{Type} column is as follows.
``Semi-convergence'':\ methods that rely on semi-convergence, cf.\ \S\ref{sec:semi}.
``Penalized'':\ methods that solve the full penalized problem, cf.\ \S\ref{sec:penalized}.
``Hybrid'':\ methods that penalize the projected problem, cf.\ \S\ref{sec:hybrid}.
``PRI'':\ methods based on penalized and/or projected restarted iterations,
cf.\ \S\ref{sec:PRI}.}
\begin{tabular}{|p{26mm}|p{81mm}|p{27mm}|p{5mm}|}\hline
\textbf{Method} & \textbf{Description} & \textbf{Type} & \textbf{Ref.} \\ \hline\hline
\texttt{IRart} & The algebraic reconstruction technique, also known as
  Kaczmarz's method. & Semi-convergence & \cite{ElHN} \\
\texttt{IRcgls} & The conjugate gradient algorithm applied implicitly to
  the normal equations. Priorconditioning allowed.
  & Penalized ($\lambda \neq 0$) \newline Semi-conv.\ ($\lambda=0$) & \cite{Hansen} \\
\texttt{IRconstr\_ls} & Projected-restarted iteration method that incorporates
  box and/or energy constraints. Priorconditioning allowed.
  & PRI & \cite{CLRS04} \\
\texttt{IRell1} & Simplified driver for \texttt{IRhybrid\_fgmres} for
  computing a 1-norm penalized solution. &
  Hybrid & \cite{GaNa14} \\
\texttt{IRenrich} & Similar to \texttt{IRcgls} but enriches the CGLS Krylov
  subspace with a low-dim.\ subspace that represents desired features of the solution.
  & Semi-convergence & \cite{enrichedCGNR} \\
\texttt{IRfista} & First-order optim.\ method FISTA that solves the Tikhonov problem
  with box and/or energy constraints; $L = I$ only.
  & Penalized ($\lambda \neq 0$) \newline Semi-conv.\ ($\lambda=0$) & \cite{FISTA} \\
\texttt{IRhtv} & Penalized restarted iteration method that
  incorporates a heuristic TV penalization term. & PRI & \cite{GaNa14} \\
\texttt{IRhybrid\_fgmres} & Hybrid version of flexible GMRES
  that applies a 1-norm penalty term to the original problem.
  & Hybrid & \cite{GaNa14} \\
\texttt{IRhybrid\_gmres} & Hybrid version of GMRES that applies a
  2-norm penalty term to the projected problem.
  Priorconditioning allowed. & Hybrid & \cite{CMRS00}, \newline \cite{GNR15} \\
\texttt{IRhybrid\_lsqr} & Hybrid version of LSQR that applies a
  2-norm penalty term to the projected problem.
  Priorconditioning allowed. & Hybrid & \cite{HyBR} \\
\texttt{IRirn} & Iteratively reweighted norm approach (penalized restarted
  iterations) for computing a 1-norm penalized solution.
  & PRI & \cite{RoWo08}  \\
\texttt{IRmrnsd} & Modified residual norm steepest descent method
  to solve nonnegatively constrained least squares problems. & Semi-convergence
  & \cite{mrnsd} \\
\texttt{IRnnfcgls} & Flexible CGLS method to solve nonnegatively constrained least squares problems.
& Semi-convergence & \cite{GaWi17} \\
\texttt{IRrestart} & A general framework for penalized and/or projected
  restarted iteration methods.
  & PRI & \cite{CLRS04}, \newline \cite{GaNa14} \\
\texttt{IRrrgmres} & Range restricted GMRES method. & Semi-convergence
  & \cite{rrgmres} \\
\texttt{IRsirt} & Simultaneous iterative reconstruction techniques (CAV,
  Cimmino, DROP, Landweber, SART). & Semi-convergence & \cite{AIRTools} \\
\hline
\end{tabular}
\end{table}

\begin{table}
\renewcommand{\arraystretch}{1.5}
\caption{\label{table:problems} Overview of the types of problems that
can be solved with this software.
The set $\mathcal{C}$ is either the box $[\mathtt{xMin},\mathtt{xMax}]^N$
or the set defined by $\| x \|_1 = \mathtt{xEnergy}$.
The matrix $L$ must have full rank.
A star $^*$ means that the function computes an \textit{approximation} to
the solution.}
\begin{tabular}{|p{53mm}|p{62mm}|} \hline
\textbf{Problem type} & \textbf{Functions} \\ \hline\hline
$\min_x \| A \, x - b \|_2^2$ \newline + semi-convergence
  & \texttt{IRart}, \texttt{IRcgls}, \texttt{IRenrich}, \texttt{IRsirt},
  \newline \texttt{IRrrgmres} ($M=N$ only)\\ \hline
$\min_x \| A \, x - b \|_2^2$ \ s.t.\ $x \geq 0$ \newline + semi-convergence
  & \texttt{IRmrnsd}, \texttt{IRnnfcgls} \\ \hline
$\min_x \| A \, x - b \|_2^2$ \ s.t.\ $x\in\mathcal{C}$ \newline
  + semi-convergence
  & \texttt{IRconstr\_ls}$^*$, \texttt{IRfista} \\ \hline
$\min_x \| A \, x - b \|_2^2 + \lambda^2 \| L\, x \|_2^2 $
  & \texttt{IRcgls}, \texttt{IRhybrid\_lsqr}, \newline
  \texttt{IRhybrid\_gmres} ($M=N$ only) \\ \hline
$\min_x \| A \, x - b \|_2^2 + \lambda^2 \| L \, x \|_2^2$ \ s.t.\ $x\in\mathcal{C}$
  & \texttt{IRconstr\_ls}$^*$, \texttt{IRfista} ($L=I$ only)\\ \hline
$\min_x \| A \, x - b \|_2^2 + \lambda \| x \|_1$
  & \texttt{IRell1}$^*$ ($M=N$ only),\newline \texttt{IRhybrid\_fgmres}$^*$ ($M=N$ only), \texttt{IRirn}$^*$ \\ \hline
$\min_x \| A \, x - b \|_2^2 + \lambda \| x \|_1$ s.t.\ $x\geq 0$,
  & \texttt{IRirn}$^*$ \\ \hline
$\min_x \| A \, x - b \|_2^2 + \lambda \mathsf{TV}(x)$ \newline with or without
  constraint $x \geq 0$ & \texttt{IRhtv}$^*$ \\ \hline
\end{tabular}
\end{table}

\subsection{Methods Relying on Semi-Convergence}
\label{sec:semi}

For many iterative methods regularization can be enforced by
terminating the process before asymptotic convergence
to the un-regularized and undesired (least squares) solution.
The underlying mechanism, which is typically referred to as \textit{semi-convergence},
is well understood, cf.\ \cite[Chapter~6]{Hansen} and the references therein.
Three of the methods in this package compute the solution $x^{(k)}$ of the problem
  \begin{equation}
  \label{eq:Krylovmethod}
    \min_x \| A\, x - b \|_2^2 \qquad \hbox{{subject to} (s.t.)} \qquad
    x \in \mathcal{S}_k \ ,
  \end{equation}
where $\mathcal{S}_k$ is a linear subspace of dimension $k$ that
takes one of the following forms:
\begin{equation}\label{eq:KrylovSp}
  \begin{array}{ll}
    \mathtt{IRcgls} & :
      \mathcal{S}_k = \mathcal{K}_k =
      \mathrm{span} \{ A^Tb, A^TA\, A^Tb, (A^TA)^2 A^Tb, \ldots (A^TA)^{k-1} A^Tb\}, \\
    \mathtt{IRenriched} & :
      \mathcal{S}_k = \mathcal{K}_k + \mathcal{W}_p, \\
    \mathtt{IRrrgmres} & :
      \mathcal{S}_k = \widehat{\mathcal{K}}_k =
      %\mathrm{span} \{ b, A \, b, A^2 \, b, \ldots A^{k-1} \, b\},
      \mathrm{span} \{A \, b, A^2 \, b, \ldots A^{k} \, b\}.
  \end{array}
\end{equation}
Here $\mathcal{K}_k$ and $\widehat{\mathcal{K}}_k$ are $k$-dimensional
Krylov subspaces, and $\mathcal{W}_p$ is a low-dimensional subspace
whose $p$ basis vectors are chosen  by the user to represent
desired features in the solution.

For \texttt{IRcgls} it is possible to apply \textit{priorconditioning} --
a type of preconditioning that modifies the underlying Krylov subspace.
Consider a Tikhonov penalization/regularization term of the form
$\pen(x) = \| L\, x \|_2^2$ with an invertible matrix~$L$\@.
In order to produce conforming iterates we introduce a new
variable $\xi$ such that $x = L^{-1} \xi$ and, implicitly, apply CGLS to the
modified problem $\min_\xi \| A\, L^{-1} \xi - b \|_2^2$, and then
compute $x^{(k)} = L^{-1} \xi^{(k)}$.
This is equivalent to solving (\ref{eq:Krylovmethod}) with $\mathcal{K}_k$
in (\ref{eq:KrylovSp}) replaced by the Krylov subspace
  \begin{equation}
  \label{eq:modKrylov}
    \mathcal{K}_{L,k} = \mathrm{span} \{ P\, A^Tb,
    (P\, A^T A)\, P\,A^Tb, (P\, A^T A)^2 P\,A^Tb, \ldots
    (P\, A^T A)^{k-1} P\,A^Tb\} ,
  \end{equation}
where $P = (L^TL)^{-1}$; see \cite[Chapter 8]{Hansen} for motivations
and details.
In this package $L$ can represent the 1D and 2D Laplacian with
zero boundary conditions, or $L$ can be a user-specified matrix
with $\mathrm{rank}(L)=N$.

Four other methods relying on semi-convergence
are based on first-order optimization methods (with step length $\omega$),
and they can incorporate constraints that can be
formulated as a projection $P_{\mathcal{C}}$ onto a convex set $\mathcal{C}$:
\begin{itemize}
\item
\texttt{IRart}, the algebraic reconstruction technique, is a row-action
method that involves each row $a_i^T$ of $A$ in a cyclic fashion:
  \begin{align*}
    y^{(k,0)} &= x^{(k)} \\
    y^{(k,i)} &= P_{\mathcal{C}}\! \left( y^{(k,i-1)} +
      \omega \, \frac{b_i - a_i^T y^{(k,i-1)}}{\| a_i \|_2^2}\, a_i \right)
      \quad \hbox{for} \quad i=1,2,\ldots,m \\
    x^{(k+1)} &= y^{(k,m)} .
  \end{align*}
The convention in this package is that one iteration
involves one sweep through all the rows.
This method {can be} understood as a projected incremental gradient descent method~\cite{AnHa}.
\item
\texttt{IRsirt} is a class of projected gradient methods of the form
  \[
    {x^{(k+1)} =  P_{\mathcal{C}}\! \left( x^{(k)} + \omega \,
    D_1\, A^T D_2\, \bigl( b - A \, x^{(k)} \bigr) \right) } ,
  \]
and the five different realizations in this package arise from different
choices of the positive diagonal matrices $D_1$ and $D_2$.
The default is the SART algorithm for which the elements of $D_1$
and $D_2$ are the 1-norms of the columns and rows of $A$, respectively.
\item
\texttt{IRfista} with regularization parameter $\lambda=0$ implements a
particular instance of the FISTA algorithm of the form
  \begin{eqnarray*}
    t_{k+1} & = & \frac{1}{2} \left( 1 + \sqrt{1 + 4\, t_k^2} \right) \\
    y^{(k+1)} & = & x^{(k)} + \frac{t_k-1}{t_{k+1}}\, \Bigl( x^{(k)} - x^{(k-1)} \Bigr) \\[2mm]
    x^{(k+1)} & = & P_{\mathcal{C}} \Bigl( y^{(k+1)} + \omega_k \,
      A^T ( b - A\, y^{(k+1)} ) \Bigr) \ .
  \end{eqnarray*}
where $\omega_k$ depends on the iteration number.
This scheme accelerates the convergence of first-order optimization methods.
\item
\texttt{IRmrnsd} is an unconstrained and modified steepest-descent algorithm
of the form
  \[
    {x^{(k+1)} = x^{(k)} + \omega_k\,
    \mathrm{diag}\bigl(x^{(k)}\bigr)\, A^T\bigl(b - A\, x^{(k)}\bigr) } ,
  \]
where the nonnegativity is imposed by the ``weight matrix''
$\mathrm{diag}(x^{(k)})$ and by bounding the step length $\omega_k$.
All elements of the initial vector must be nonnegative.
\end{itemize}
Yet another method depends on semi-convergence:
\texttt{IRnnfcgls} is a particular implementation of the flexible CGLS
algorithm that uses a judiciously constructed preconditioner,
which changes in every iteration,
to ensure convergence to a non-negative solution \cite{GaWi17}.

Nonnegativity constraints are hardwired into \texttt{IRmrnsd} and \texttt{IRnnfcgls},
while the other three methods can incorporate general box constraints
(with nonnegativity as a special case), as well as a so-called
\textit{energy constraint}, which has the form
  \[
    \| x \|_1 = \mathrm{constant} \ ,
  \]
where the constant is specified by the user.

\subsection{Methods for Solving the Penalized Least Squares Problem}
\label{sec:penalized}

Three of the methods in the above category, \texttt{IRcgls}, \texttt{IRenrich}
and \texttt{IRfista},
can also be used to solve the penalized least-squares problem (\ref{eq:penalized})
with $\pen(x) = \| L \, x \|_2^2$ (i.e., Tikhonov regularization),
which corresponds to the least squares problem
  \begin{equation}
  \label{eq:Tikhonov}
    \min_x \left\| \begin{pmatrix} A \\ \lambda L \end{pmatrix} x -
    \begin{pmatrix} b \\ 0 \end{pmatrix} \right\|_2^2 \qquad
    \Leftrightarrow \qquad x = (A^TA+{\lambda^2}\,L^TL)^{-1}A^Tb .
  \end{equation}
In this case we ignore semi-convergence and instead rely on
asymptotic convergence to the penalized solution in (\ref{eq:Tikhonov}).
In \texttt{IRcgls} the matrix $L$ can be either the identity matrix or any of the
matrices described as priorconditioners in Section \ref{sec:semi}.
In \texttt{IRenrich} and \texttt{IRfista} only $L=I$ is allowed,
and \texttt{IRfista} has the option to
also incorporate box constraints and/or the energy constraint.

Two other penalization functions can be handled:\ the 1-norm,
$\pen(x) = \| x \|_1$, which enforces sparsity on $x$, and
$\pen(x) = \mathsf{TV}(x)$, where the total variation (TV) function
is defined in a discrete setting by
  \[
    \mathsf{TV}(x) = \sum \sqrt{ [ D_{\mathrm{h}} x]_i^2 + [D_{\mathrm{v}} x]_i^2 } \ .
  \]
Here the two matrices $D_{\mathrm{h}}$ and $D_{\mathrm{v}}$ compute finite
difference approximations to the horizontal and vertical partial derivatives,
respectively, and the sum is over all elements of $x$ for which these can
be computed.
These penalized problems are solved approximately by means of our
implementations of particular hybrid methods
\texttt{IRell1}, \texttt{IRhtv} and \texttt{IRirn};
hybrid methods are described in more detail in the following subsection.

\subsection{Hybrid Krylov Subspace Methods that Regularize the Projected Problem}
\label{sec:hybrid}

In hybrid Krylov subspace methods the penalization is moved from the
``original problem'' (\ref{eq:penalized}) to the ``projected problem'',
i.e., the least squares problem restricted to the Krylov subspace~\cite{ChKN}.
The main advantage is that the search for a good
regularization parameter is done on the projected problem,
which has relatively small dimensions and is therefore
less computationally demanding than working with the original large-scale problem.
This means that the regularization parameter is iteration dependent,
and is adjusted as the Krylov subspace grows.
Therefore, we use the notation $\lambda_k$ to denote the regularization
parameter corresponding to the $k$th iteration.
We provide three hybrid methods:
\begin{itemize}
\item
\texttt{IRhybrid\_lsqr} is, similarly to \texttt{IRcgls},
based on the Krylov subspace~$\mathcal{K}_k$ defined in (\ref{eq:KrylovSp});
the underlying LSQR
method explicitly builds an orthonormal basis for this space
allowing us to easily formulate and solve the penalized
projected problem.
The default approach for choosing the regularization parameter $\lambda_k$ for the projected problem is
weighted generalized cross validation (GCV).
\item
\texttt{IRhybrid\_gmres} follows the same idea, except that it is based on
the Krylov subspace~$\mathrm{span}\{b,\,\widehat{\mathcal{K}}_{k-1}\}$,
where $\widehat{\mathcal{K}}_{k-1}$ is analogous to the subspace defined
in (\ref{eq:KrylovSp}). By default it uses GCV to determine the
regularization parameter $\lambda_k$ for the projected problem.
\item
\texttt{IRhybrid\_fgmres} is based on a flexible
version of the approximation subspace used for \texttt{IRhybrid\_gmres},
which incorporates an iteration dependent preconditioner
whose role is to emulate a 1-norm (sparsity) penalty term on the solution.
By default it uses GCV to determine the regularization parameter $\lambda_k$
for the projected problem.
\end{itemize}
We note that, with $\lambda=0$, the hybrid LSQR algorithm in \texttt{IRhybrd\_lsqr}
is mathematically equivalent to LSQR -- as well as CGLS, available in
\texttt{IRcgls} with $\lambda=0$.
Similarly,
{with $\lambda=0$} the hybrid GMRES algorithm in \texttt{IRhybrid\_gmres} is
mathematically equivalent to the GMRES algorithm.

We also note that when $\lambda \neq 0$ and $L\neq I$, the Krylov subspace
in (\ref{eq:modKrylov}) that underlies the hybrid LSQR algorithm is different
from the Krylov subspace underlying CGLS when applied to the Tikhonov problem
(\ref{eq:Tikhonov}) -- although they are identical when $L=I$;
see \cite{KiHE07} for details.

\subsection{Penalized and/or Projected Restarted Iterations (PRI)}
\label{sec:PRI}

These functions are based on restarted inner-outer iterations.
Semi-convergent or penalized Krylov methods, or hybrid iterative solvers,
are used in the inner iterations,
and every outer iteration produces a new {approximate solution} that incorporates the
desired properties or constraints.
This general framework is implemented in the function \texttt{IRrestart},
which is called by other functions (\texttt{IRconstr\_ls}, \texttt{IRhtv} and
\texttt{IRirn}) with more specific goals.
The experienced user can run \texttt{IRrestart} in such a way that a
variety of combinations of inner solvers and outer constraints are
heuristically incorporated into the approximate solution,
and may wish to add further application-specific constraints.
\texttt{IRrestart} can handle penalized and/or projected schemes as
detailed below.

\begin{tabbing}
xxx \= xxx \= xxx \kill
\> \underline{Penalized Restarted Iterations} \\[1mm]
\> Initialize $x^{(0)}$ and $L_0$ \\
\> for $\ell=0,1,2,\ldots$ \\
\> \> $r^{(\ell)} = b - A\, x^{(\ell)}$ \\
\> \> $w^{(\ell)} = \arg\min_w \| A\, w-r^{(\ell)}\|_2^2 +
      \lambda_\ell^2 \, \| L_\ell\, w \|_2^2$ \\
\> \> $x^{(\ell+1)} = \left\{ \begin{array}{l} x^{(\ell)} + w^{(\ell)} \\
       P_{\mathcal{C}} \bigl( x^{(\ell)} + w^{(\ell)} \bigr) \end{array} \right.$
      depending on the user's choice \\
\>\> {update $L_{\ell+1}$} \\
\> end
\end{tabbing}
\begin{tabbing}
xxx \= xxx \= xxx \kill
\> \underline{Projected Restarted Iterations} \\[1mm]
\> Initialize $x^{(0)}$  \\
\> for $\ell=0,1,2,\ldots$ \\
\> \> $r^{(\ell)} = b - A\, x^{(\ell)}$ \\
\> \> $w^{(\ell)}$ = $\arg\min_w \| A\, w - r^{(\ell)}\|_2$ \\
\> \> $x^{(\ell+1)} = P_{\mathcal{C}} \bigl( x^{(\ell)} + w^{(\ell)} \bigr)$ \\
\> end
\end{tabbing}
{
Note that, in addition to updating the regularization matrix $L_\ell$,
the user can also choose to incorporate a projection
at each outer iteration of the Penalized Restarted Iterations.
For methods based on restarts, the concept of total number of iterations,
i.e., the number of iterations performed jointly in the inner and outer iterations,
should be considered.
}

Computation of the update $w^{(\ell)}$ at the $\ell$th outer iteration is performed
by means of some of the iterative solvers in this package.
The number of inner iterations in these solvers acts as a regularization
parameter and is always chosen by one of the stopping-rule methods discussed in
Section~\ref{sec:StopR} below.
We emphasize that this has the consequence that even without a stopping rule for the
outer iterations (except for the maximum number of iterations),
the specific mandatory choice of stopping rule for the \textit{inner iterations}
influences the behavior and convergence of the outer iterations.

We note that most of these restarted iterations
can be regarded as a heuristic approach to resemble first-order
optimization methods and, in particular, they are reminiscent of an
alternating projection scheme onto convex sets.
We will not pursue this aspect further here.

There are three functions that act as easy-to-use drivers to \texttt{IRrestart}
for specific purposes.

The function \texttt{IRconstr\_ls} uses the restarted iterations to
enforce box and energy constraints, by projection onto the relevant convex sets
at each outer iteration.
The two functions \texttt{IRhtv} and \texttt{IRirn} use the restarted
iterations to approximate a penalized solution with penalty term
$\pen(x) = \mathsf{TV}(x)$ and $\pen(x)=\| x \|_1$, respectively.
The penalty is enforced through a 2-norm $\| L_\ell \cdot \|_2$,
where the matrix $L_\ell$ is chosen to enforce the desired penalty;
it depends on the current iterate~$x^{(\ell)}$ as follows:
\begin{itemize}
\item
\texttt{IRhtv}: $L_\ell = \begin{pmatrix}
  \widehat{L}_\ell D_{\mathrm{h}} \\ \widehat{L}_{\ell} D_{\mathrm{v}} \end{pmatrix}$
  \ with \ $\widehat{L}_\ell = \mathrm{diag} \! \left( \left( (D_{\mathrm{h}}x^{(\ell)})_i^2 +
  (D_{\mathrm{v}}x^{(\ell)})_i^2 \right)^{-1/4} \right)$.
\item
\texttt{IRirn}: $L_\ell = \mathrm{diag}\!\left( \bigl| x^{(\ell)}_i \bigr|^{-1/2}
 \right)$.
\end{itemize}

\subsection{Stopping Rules {and Parameter Choice Strategies}}\label{sec:StopR}

Since the iterative solvers in this package are designed for
{regularization of} inverse problems,
we provide well-known stopping rules for such problems.
Also parameter choice strategies for setting the regularization parameter
$\lambda_k$ for hybrid methods are surveyed:
these are related to the discrepancy principle and generalized cross validation.

The basic idea behind the \textit{discrepancy principle} is to stop as soon as the
norm of the residual $b - A\, x^{(k)}$ is sufficiently small, typically of
the same size as the norm of the perturbation $e$ of the right-hand
side, cf.\ \cite[\S 5.2]{Hansen}.
In this package, where all norms are relative, this takes the form
  \begin{center}
    stop as soon as $\| b - A\, x^{(k)} \|_2/\| b \|_2 \leq \eta
    \cdot \mathtt{NoiseLevel}$ ,
  \end{center}
where $\eta$ is a ``safety factor'' slightly larger than 1,
and \texttt{NoiseLevel} is the relative noise level $\| e \|_2/\| b \|_2$.
If the noise level is not specified, then the default value used in all codes is 0.
To solve a noise-free problem with a given threshold $\tau$, the user may
set $\eta = 1$ and $\mathtt{NoiseLevel} = \tau$.
The specific implementation of this stopping criterion takes different forms,
depending on the circumstances:
\begin{itemize}
\item
For the functions that leverage semi-convergence, \texttt{IRart}, \texttt{IRcgls},
\texttt{IRenrich}, \texttt{IRfista}, \texttt{IRmrnsd}, \texttt{IRnnfcgls},
\texttt{IRrrgmres} and \texttt{IRsirt}, the implementation is done in
a straight-forward way.
\item
For the functions that use hybrid methods,
\texttt{IRell1}, \texttt{IRhbyrid\_fgmres}, \linebreak[4]\texttt{IRhybrid\_gmres}
and \texttt{IRhybrid\_lsqr},
we implemented the ``secant method'' from \cite{GaNo14}, which updates the
regularization parameter for the projected problem in such a way
that stopping by the discrepancy principle is ensured.
\item
For the functions that use inner-outer iterations,
\texttt{IRconstr\_ls}, \texttt{IRhtv}, \texttt{IRirn} and
\texttt{IRrestart}, the discrepancy principle can be applied to
the solver in the inner iterations, and the outer iterations are terminated
when either $\| x^{(\ell)} \|_2$, $\| L\, x^{(\ell)} \|_2$,  or the value of
the regularization parameter, at each restart, has stabilized.
The choice is controlled by \texttt{options.stopOut} which accepts
the values \texttt{'xstab'}, \texttt{'Lxstab'} and \texttt{'regPstab'}.
\end{itemize}

The basic idea behind \textit{generalized cross validation} (GCV) is to choose
the solution that gives the best prediction of the unperturbed data, cf.\
\cite[\S 5.4]{Hansen}.
This method is practical only for the hybrid methods, where it can
be applied to the projected problem.
Let $W_k$ be a matrix with orthonormal columns that span the relevant
Krylov subspace for the approximation of the solution, and let $A\, W_k$ have the factorization
$A\, W_k = Z_{k+1}\, R_k$, where $Z_{k+1}$ has orthonormal columns and $R_k$ is
either lower bidiagonal or upper Hessenberg, depending on the chosen Krylov method.
Then we apply Tikhonov regularization to the projected problem
$\min_{y\in\mathbb{R}^k} \| R_k\, y - Z_{k+1}^T b \|_2$ to obtain
$y_{\lambda}^{(k)} = R_k^{\#}(\lambda) Z_{k+1}^T b$, where $R_k^{\#}(\lambda)$
is a ``fictive'' matrix that defines the regularized solution.
The regularization parameter $\lambda_k$ minimizes the GCV function
  \[
    G_k(\lambda) = \frac{\| R_k\, y_{\lambda}^{(k)} - Z_{k+1}^T b \|_2}
      {Q - w\, \mathrm{trace}\bigr( R_k\, R_k^{\#}(\lambda) \bigl)} \ ,
  \]
and we provide three different variants of this function:
\begin{tabbing}
xxx \= xxxxxxxxxxxxxxxx \= xxxxxxxxxxxxx \= xxx \kill
\> standard GCV: \> $Q=k+1$ \> $w=1$ \quad (cf.~\cite{GHW}), \\
\> modified GCV: \> $Q=M-k$ \> $w=1$ \quad (cf.~\cite{NoRu14}), \\
\> weighted GCV: \> $Q=k+1$ \> $w<1$ \quad (cf.~\cite{HyBR}).
\end{tabbing}
Once $\lambda_k$ is determined we put $x^{(k)} = W_k\, y_{\lambda_k}^{(k)}$.
The iterations are terminated as soon as one of these conditions is satisfied:
\begin{itemize}
\item
The minimum of $G_k(\lambda)$, as a function of $k$, stabilizes or starts
to increase within a given iteration window.
\item
The residual norm $\| b - A\, x^{(k)} \|_2$ stabilizes.
\end{itemize}
When GCV is applied to methods that use inner-outer iterations,
similarly to the discrepancy principle case, the GCV is applied to
the inner iterations, and the outer iterations are terminated when some
stabilization occurs in $\|x^{(\ell)}\|_2$, $\|Lx^{(\ell)}\|_2$, or
the regularization parameter.

In addition to these stopping rules, there are cases where semi-convergence
is not relevant -- either because the data is
noise-free or because we iteratively solve the Tikhonov problem.
In these cases it is preferable to terminate the iterations when the
residual for the (penalized) normal equations is small, i.e.,
  \begin{center}
    stop as soon as  $\| A^Tb - (A^TA + \lambda^2 L^T L) \, x^{(k)} \|_2/{\|A^T b \|_2}
    \leq \mathtt{NE\_Rtol}$ ,
  \end{center}
including the case $\lambda=0$.
This stopping rule can be used in the functions
\texttt{IRcgls}, \texttt{IRenrich}, \texttt{IRfista} and \texttt{IRmrnsd}.

\section{Overview of the Test Problems}
\label{sec:Problems}

\begin{table}
\renewcommand{\arraystretch}{1.5}
\caption{\label{table:testproblems} Overview of the types of test
problems provided in this package, plus some related functions.
The problems \texttt{PRseismic}, \texttt{PRspherical} and \texttt{PRtomo}
require \textsc{AIR Tools II} \cite{AIRTools}.}
\begin{tabular}{|p{47mm}|p{40mm}|p{50mm}|} \hline
\textbf{Test problem type} & \textbf{Function} & \textbf{Type of $A$} \\ \hline\hline
Image deblurring & \texttt{PRblur} (generic function) & \\
-- spatially invariant blur &
  \texttt{PRblurdefocus}, \texttt{PRdeblurgauss}, \newline
  \texttt{PRdeblurmotion}, \texttt{PRdelburshake}, \newline
  \texttt{PRdeblurspeckle} & Object \\
-- spatially variant blur & \texttt{PRblurrotation} & Sparse matrix \\
Inverse diffusion & \texttt{PRdiffusion} & Function handle \\
Inverse interpolation & \texttt{PRinvinterp2} & Function handle \\
NMR relaxometry & \texttt{PRnmr} & Function handle \\
Tomography & & \\
-- seismic travel-time tomography & \texttt{PRseismic} &
  Sparse matrix or function handle \\
-- spherical means tomography & \texttt{PRspherical} & Sparse matrix or function handle \\
-- X-ray computed tomography & \texttt{PRtomo}
  & Sparse matrix or function handle \\
\hline
Add noise to the data: \newline Gauss, Laplace, %\newline
  multiplicative & \texttt{PRnoise} & \\
Visualize the data $b$ and the solution $x$ in appropriate formats
  & \texttt{PRshowb}, \texttt{PRshowx} & \\
\hline
Auxiliary functions for some test problems
  & \texttt{OPdiffusion}, \texttt{OPinvinterp2}, \newline \texttt{OPnmr} & \\
\hline
\end{tabular}
\end{table}

While realistic test problems are crucial for testing, debugging and
demonstrating algorithms to solve inverse problems, there are very few
collections available.
One exception is the set of simple 1D test problems in
\textsc{Regularization Tools}, but they are outdated and do not
represent current large-scale applications.
For this reason, we find it necessary to provide a new set of
more realistic 2D test problems that are better suited for
testing algorithms that
are designed especially for large-scale applications,
such as the iterative methods implemented in this package.
When choosing these test problems we had the following criteria in mind:
\begin{itemize}
\item
The functions for generating the test problems must be easy to use,
with good choices of default parameters.
\item
At the same time, the user {should} have full control of the underlying
model parameters via an \texttt{options} input.
\item
The test problems can be used as ``black boxes'' without any specific
knowledge about the application domain.
\item
It must be easy to add noise to the data.
\item
The right-hand side $b$ (the data) and the solution $x$ can be easily
visualized.
\end{itemize}

The functions for generating the test problems, together with a few
auxiliary functions, are listed in Table~\ref{table:testproblems}.
Although the test problems represent a variety of applications, they all use
the same calling sequence,
\begin{verbatim}
   [A, b, x, ProbInfo] = PR___(n, options);
\end{verbatim}
with two inputs:\ \texttt{n}, which defines the problem size,
and the structure \texttt{options} for setting the model parameters.
Either or both can be omitted, and default options produce a suitable
test problem of medium difficulty.
{Note that throughout the paper, and in all of the implemented test problems,
the input {\tt n} (lower case) defines the problem size, but does not necessarily
give explicit information about the actual sizes of the matrix {\tt A} and vectors {\tt x} and {\tt b}.
We use the convention that $M \times N$ (i.e., with the use of upper case letters $M$ and $N$)
denotes the dimensions of the matrix {\tt A}; the precise relationship
between {\tt n} and $M$ and $N$ depends
on the application. For example, in an image deblurring problem, the input {\tt n}
creates a test problem with images having $n\times n$ pixels,
and $M = N=n^2$. The help documentation for each of the {\tt PR\_\_\_} test problems provides
more details, and can be viewed with MATLAB's {\tt help} or {\tt doc} commands.}

In the output parameters, \texttt{A} represents the forward operation,
\texttt{b} is a vector with the noise-free data,
\texttt{x} is a vector with the true solution, and
\texttt{ProbInfo} is a structure that contains useful information about the problem
(such as image dimensions, problem type, and important model parameters).
The type of \texttt{A} depends on the test problem:
\begin{itemize}
\item
For image deblurring, \texttt{A} is either an object that follows the
conventions from \textsc{Restore Tools} \cite{RestoreTools}, or
a sparse matrix (depending on the type of blurring).
\item
For inverse diffusion, inverse interpolation,  and NMR relaxometry, \texttt{A} is
a function handle that gives easy access to functions, written by us
and stored as \texttt{OP\_\_\_} files, that perform matrix-vector multiplications.
\item
For the tomography problems, the user can choose \texttt{A} to be a sparse matrix
or a function handle; the former gives faster execution but requires
more memory, while the latter executes slowly but has very limited
memory requirements.
\end{itemize}
When a function handle is used for {\tt A}, then our iterative methods
expect {\tt A} to conform to the following definitions:
\begin{itemize}
\item[]
\begin{tabular}{ll}
   {\tt u = A(x, 'notransp');} & computes the matrix-vector multiplication $u = A\,x$, \\
   {\tt v = A(y, 'transp');} & computes the matrix-vector multiplication $v = A^Ty$, \\
   {\tt dims = A([], 'size');} & returns the dimensions of the matrix $A$,
\end{tabular}
\end{itemize}
that is, \texttt{dims(1)} = $M$ and \texttt{dims(2)} = $N$, the dimensions
of~\texttt{A}. In some cases (e.g.,  inverse diffusion) it may be
difficult to implement the multiplication with $A^T$. In these cases,
only transpose-free iterative methods should be used.
Note that our test problems illustrate the three possibilities (sparse matrix,
user-defined object, and function handle)
for representing the problem{s} that can be handled by our software,
and they provide templates for users who want to write code
for their own problems.

The input parameter \texttt{options} is a structure that can be used to override
various default options.
To determine what the possible default options for the various
test problems are, use:
\begin{verbatim}
   options = PR___('defaults');
\end{verbatim}
One can then change the default options either by directly changing a specific field,
for example,
\begin{verbatim}
   options.field_name = field_value;
\end{verbatim}
or by using the function \texttt{PRset},
\begin{verbatim}
   options = PRset(options, 'field_name', field_value);
\end{verbatim}
Note that in the above example using {\tt PRset}, it is assumed that the structure
\texttt{options} is already defined, and only one of its field values is changed.
It is possible to change multiple field values using {\tt PRset}, for example,
\begin{verbatim}
   options = PRset(options, 'field_name1', field_value1, ...
         'field_name2', field_value2, 'field_name3', field_value3);
\end{verbatim}
It is also possible to use \texttt{PRset} without a pre-defined \texttt{options}
structure, such as
\begin{verbatim}
   options = PRset('field_name', field_value);
\end{verbatim}
In this case, all default options are used, except \texttt{field\_name}.
In the following subsections we provide some additional specific examples.

\subsection{Image Deblurring}
\label{sec:blur}

Image deblurring (which is sometimes referred to as image restoration or deconvolution)
is an inverse problem that reconstructs an image from a blurred and noisy observation.
Image deblurring problems arise in many {important} applications, such as astronomy, microscopy,
crowd surveillance,
{just to name a few} \cite{AnHu,BeBo98,HaNaOL06,LagBie}.
A mathematical model of this problem can be expressed in
the continuous setting as an integral equation
  \begin{equation}
  \label{eq:IFK}
    g(s) = \int k(s,t) \, f(t) \, ds + e(s)\,,
  \end{equation}
where $s,t \in \mathbb{R}^2$.
The kernel $k(s,t)$ is a function that specifies how the points in the image
are distorted, and is therefore called the \textit{point spread function} (PSF).
If the kernel has the property that $k(s,t) = k(s-t)$, then the PSF is said
to be spatially invariant; otherwise, it is said to be spatially variant.

In a realistic {setting}, images are collected only at discrete points (pixels),
and are also only available in a finite region.
Therefore one must usually work directly with the discrete model (\ref{eq:Axb})
where $b$
{and $x$ are vectors that represent the blurred and
sharp images,}
and $A$ is a large,
usually ill-conditioned matrix that models the blurring operation.

From equation (\ref{eq:IFK}) it can be observed that each pixel in the blurred
image is formed by integrating the PSF with pixel values of the true
image scene. Generally the integration operation is local, and so pixels
in the center of the viewable region are well defined by the linear system
(\ref{eq:Axb}).
However, pixels of the blurred image near the boundary of the viewable
region are affected by information outside the viewable region.
Therefore, in constructing
the matrix $A$, one needs to incorporate boundary conditions to model how
the image scene extends beyond the boundaries of the viewable region.
Typical boundary conditions include zero, periodic, and reflective \cite{HaNaOL06}.
Note that it is generally not possible to know precisely what values should be
assigned to $x$ outside the borders of the viewable region, and so
even in the noise-free case (i.e., $e=0$), the product $Ax$ is unlikely to be
exactly equal to $b$.

\IRT\ includes several test problems with various blurring operations:
\begin{itemize}
\item
{\tt PRblurdefocus} simulates a spatially invariant, out-of-focus blur.
\item
{\tt PRblurgauss} simulates a spatially invariant Gaussian blur.
\item
{\tt PRblurmotion} is a spatially invariant blur that simulates relative
{linear} motion, at a 45 degree angle, between an imaging device and the scene.
\item
{\tt PRblurrotation} simulates a spatially variant rotational motion blur
{around the center of the image}.
\item
{\tt PRblurshake} simulates spatially invariant motion blur caused by shaking of a
camera.  The path of motion is generated randomly, so repeated
calls to {\tt PRblurshake} will create different blurring operators, unless
the random number generator is manually set to a specific value using
MATLAB's built-in {\tt rng} function.
\item
{\tt PRblurspeckle} simulates spatially invariant blurring caused by atmospheric turbulence.
\end{itemize}
As stated earlier in this section, these test problems can be called as follows:
\begin{verbatim}
   [A, b, x, ProbInfo] = PRblur___(n, options);
\end{verbatim}
The two inputs, {\tt n} and {\tt options} are optional; if they are not specified, default
values are used (e.g., the default value for {\tt n} is 256).
In the case of spatially invariant blur examples, {\tt A} is a {\tt psfMatrix} object that
overloads various standard MATLAB operations\footnote{Another possibly useful
overloaded operation is the {\tt full} command that, as is done for sparse matrices, transforms 
a {\tt psfMatrix} to full storage organization.  This operation, however, should only be
used for relatively small problems.}, such as * to efficiently
implement matrix vector multiplications with $A$ and $A^T$; for further details, see
\cite{RestoreTools}.
In the case of spatially variant rotational motion blur, {\tt A} is a sparse matrix
\cite{HaNT14}. 

As will be illustrated in Section~\ref{sec:Examples}, it is very easy to use
the iterative methods we provide in the package with {\tt A} for either
the {\tt psfMatrix} object or sparse matrix format.
It is also easy for users to test
their own iterative methods with these problems because
matrix-vector multiplies can be computed using standard MATLAB operators, such as
\begin{verbatim}
   r = A'*(b - A*x);
\end{verbatim}
In addition, the effective matrix size of {\tt A} (if it could be constructed explicitly as a full matrix)
can be found using MATLAB's built-in {\tt size} function. For example,
with the default \texttt{n} = 256, then
\begin{verbatim}
   dims = size(A);
\end{verbatim}
returns the vector {\tt dims = [65536, 65536]}.

The \texttt{options} structure can be used to set the boundary conditions
to zero, periodic, or reflective; if nothing is specified, the default choice
is reflective.
It is possible to construct a problem where $A\, x$ is exactly equal to $b$;
that is, the specified boundary conditions used to construct $A$
exactly model how $x$ behaves outside the viewable region.
Because this situation is unrealistic, we consider it to be a classic example of
committing an ``inverse crime".  To construct such an example, use the options
structure:
\begin{verbatim}
   options = PRset('CommitCrime', 'on');
   [A, b, x, ProbInfo] = PRblur___(options);
\end{verbatim}
The structure {\tt options} can also be used to modify a variety of other default
parameters, including:
\begin{itemize}
\item
\texttt{options.trueImage} can be used to choose one of several (true scene) test images
provided in the package, or it can be a user-defined test image; it is returned
in the vector {\tt x}.  Default is an image of
the Hubble Space Telescope.
\item
\texttt{options.PSF} can be used to choose one of several point spread functions implemented in the package, or it can be used to set a user-defined PSF, stored as a two-dimensional array. Default is a Gaussian PSF.
\item
\texttt{options.BlurLevel} sets the severity of blur; choices are
{\tt 'mild'}, {\tt 'medium'} (default), or {\tt 'severe'}.
\item
\texttt{options.BC} sets the boundary conditions; choices are
{\tt 'zero'}, {\tt 'periodic'}, or {\tt 'reflective'} (default).
\end{itemize}

\begin{figure}
\begin{center}
\begin{tabular}{ccc}
\includegraphics[width=0.3\textwidth]{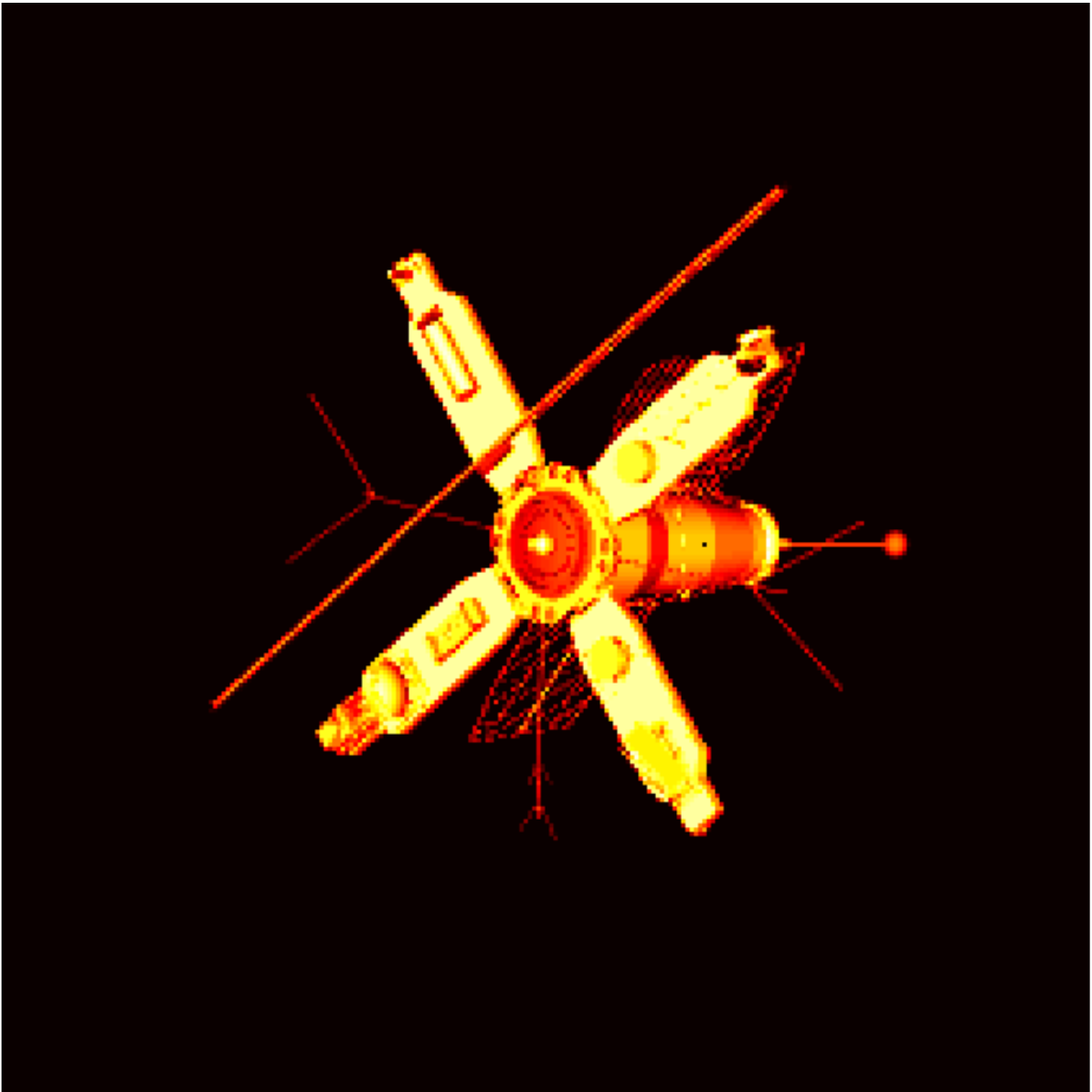} &
\includegraphics[width=0.3\textwidth]{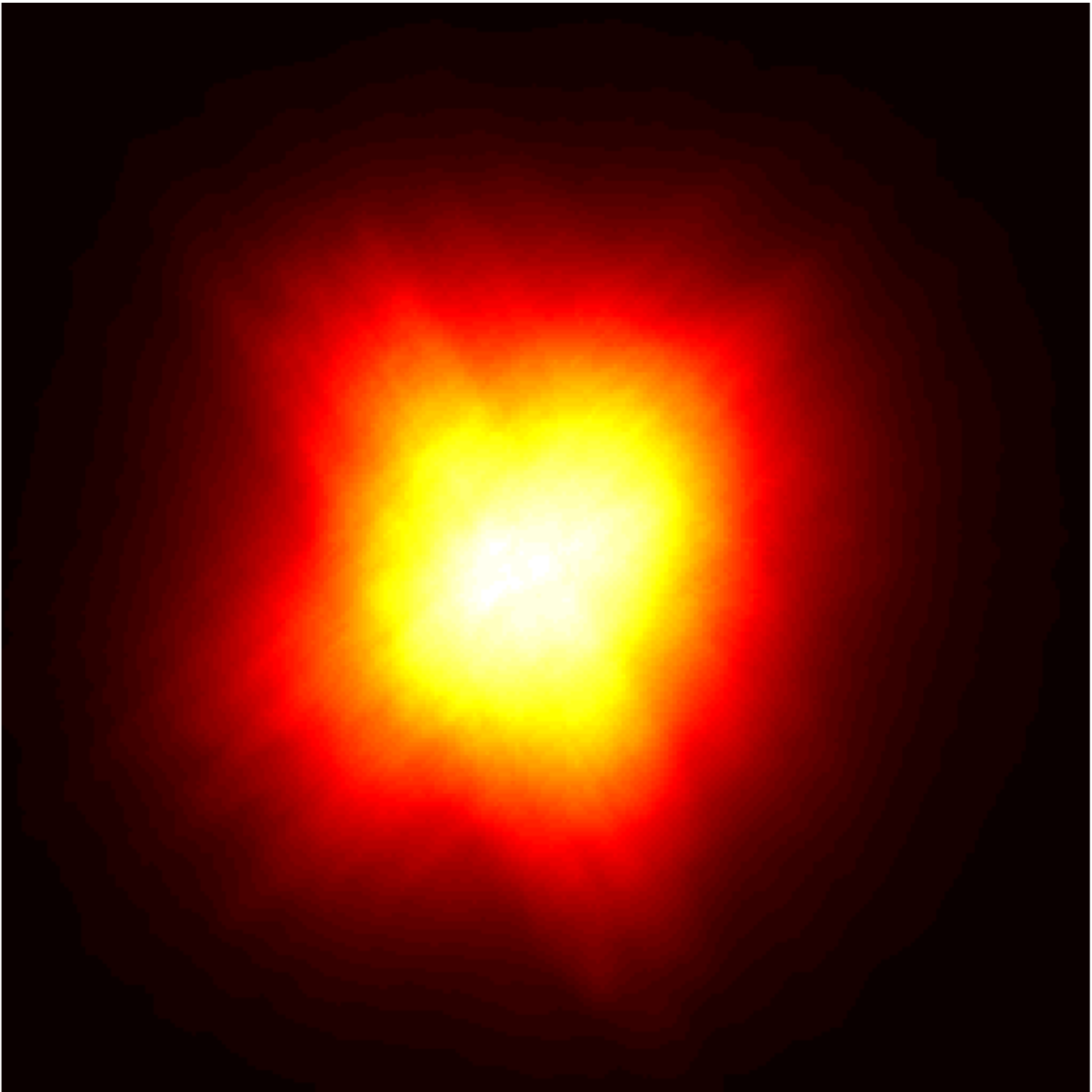} &
\includegraphics[width=0.3\textwidth]{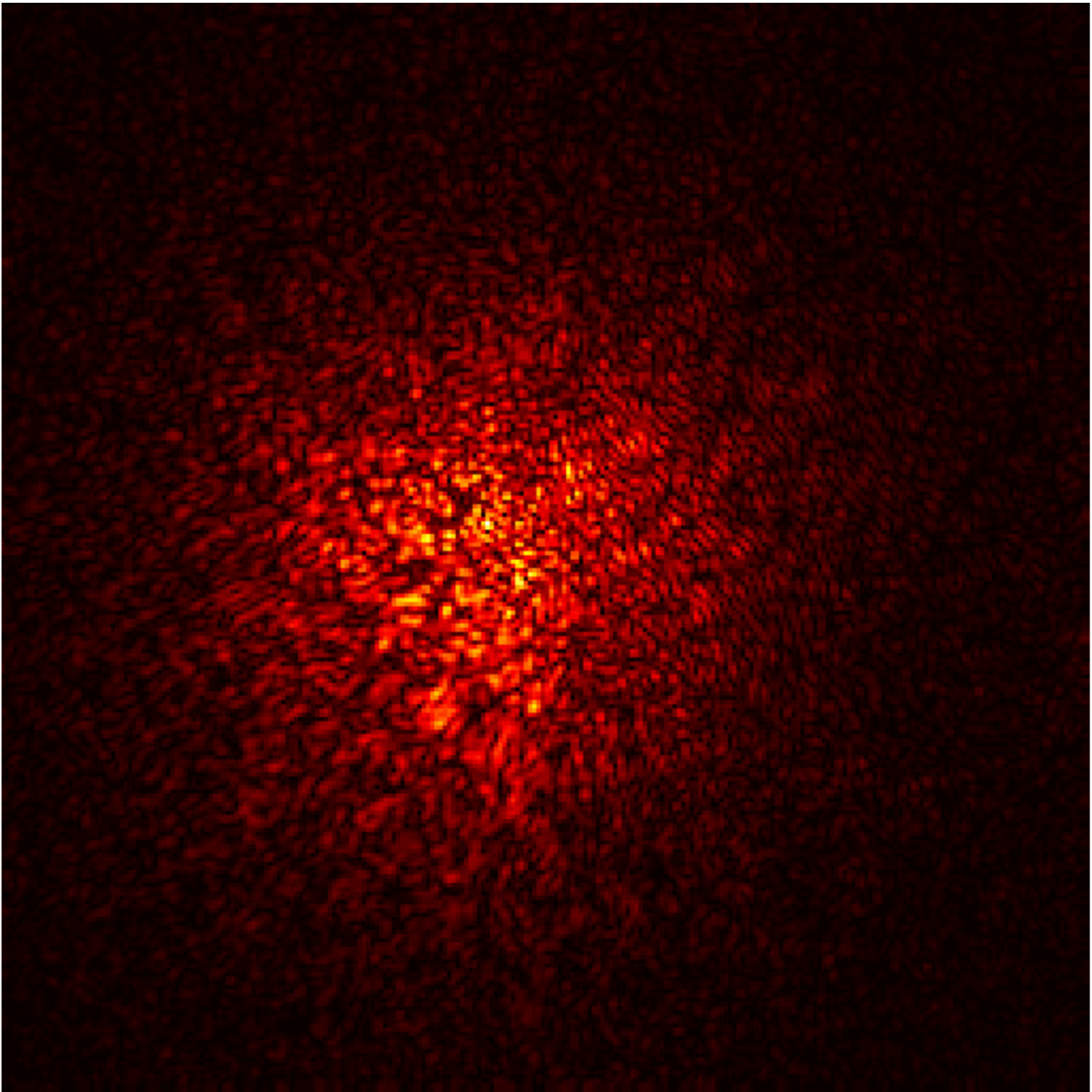}
%{\footnotesize True image, {\tt x}} & {\footnotesize Blurred image, {\tt b}} & {\footnotesize PSF}
\end{tabular}
\end{center}
\caption{\label{fig:SpeckleBlurExample} Data produced by the test problem
{\tt PRblurspeckle}, with the true image scene shown on the left,
the blurred image scene shown in the middle, and the PSF (which
defines the matrix {\tt A}) shown on the right. The PSF is displayed on
a square root scale.}
\end{figure}

We close this subsection with an example, where we generate a speckle blur test problem,
choosing the (non-default) test image
{\tt 'satellite'}, and reset the blur level to {\tt 'severe'}:
\begin{verbatim}
   options = PRset('trueImage', 'satellite', 'BlurLevel', 'severe');
   [A, b, x, ProbInfo] = PRblurspeckle(options);
\end{verbatim}
The {vectors \texttt{b} and \texttt{x}} produced by this test problem
can be displayed using {\tt PRshowb} and
{\tt PRshowx},
\begin{verbatim}
   PRshowb(b, ProbInfo)
   PRshowx(x, ProbInfo)
\end{verbatim}
We could also display the PSF using either of these ``show" functions, or
by using MATLAB's standard {\tt mesh} command:
\begin{verbatim}
   PRshowx(ProbInfo.psf, ProbInfo)
   mesh(ProbInfo.psf)
\end{verbatim}
In each of the ``show" cases, we change the colormap to {\tt hot}, and in the PSF case {we} use
a square root scale to display {the image intensity}.
The results are shown in Figure~\ref{fig:SpeckleBlurExample}.

\subsection{Inverse Interpolation}
\label{sec:inverseinterpolation}

Inverse interpolation -- also known as gridding --
is the problem of computing the values of a function
on a regular grid, given function values on arbitrarily located points,
in such a way that interpolation of the gridded function values (the unknowns)
produces the given values (the data) \cite{GMS12}.
This is obviously an inverse problem whose specifics depend
on the type of interpolation being used.
A different algorithm than ours is implemented in MATLAB's
\texttt{griddata} function.

As a simple example, consider linear interpolation on a 1D grid
with grid points $t_j^{\mathrm{G}}$, $j=1,2,\ldots$ and data
$(t_i,b_i)$ on the arbitrary points $t_1<t_2<t_3<\ldots$;
then the unknown function values $x_j$ at the grid points
must satisfy the interpolation relations (for all $i$):
  \[
    b_i = x_j + \frac{t_i - t_j^{\mathrm{G}}}{t_{j+1}^{\mathrm{G}}-t_j^{\mathrm{G}}}
    \, (x_{j+1}-x_j) , \quad \hbox{where} \quad t_i \in
    \left[ \,t_j^{\mathrm{G}} \, , \, t_{j+1}^{\mathrm{G}} \, \right] .
  \]
This gives a simple linear system of equations $A\, x=b$ with a sparse
coefficient matrix $A$ (two nonzeros per row).
Note that $A$ is rank deficient if there are consecutive grid
intervals with no data points.

\begin{figure}
\begin{center}
\includegraphics[width=0.45\textwidth]{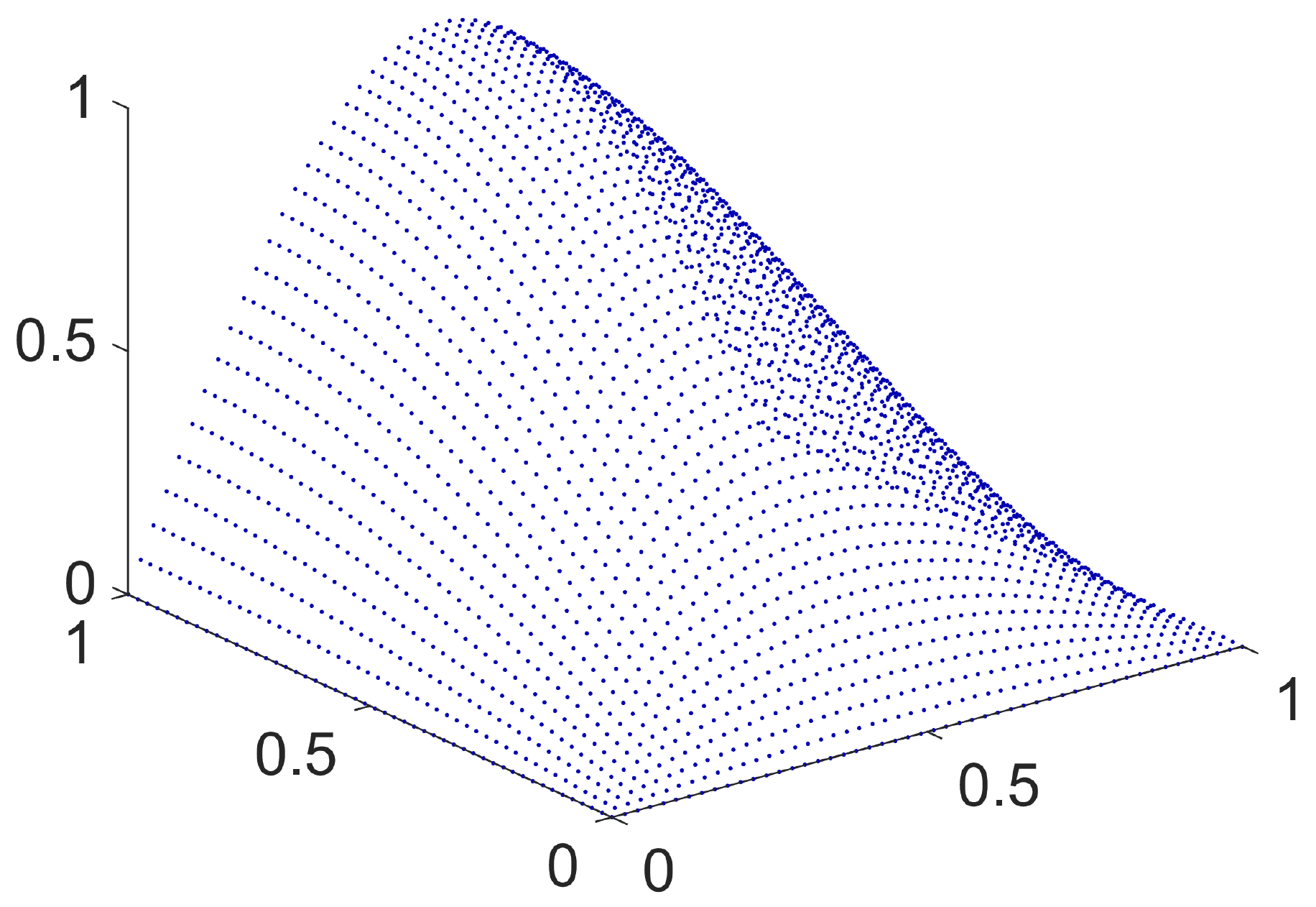} \hspace{4mm}
\includegraphics[width=0.45\textwidth]{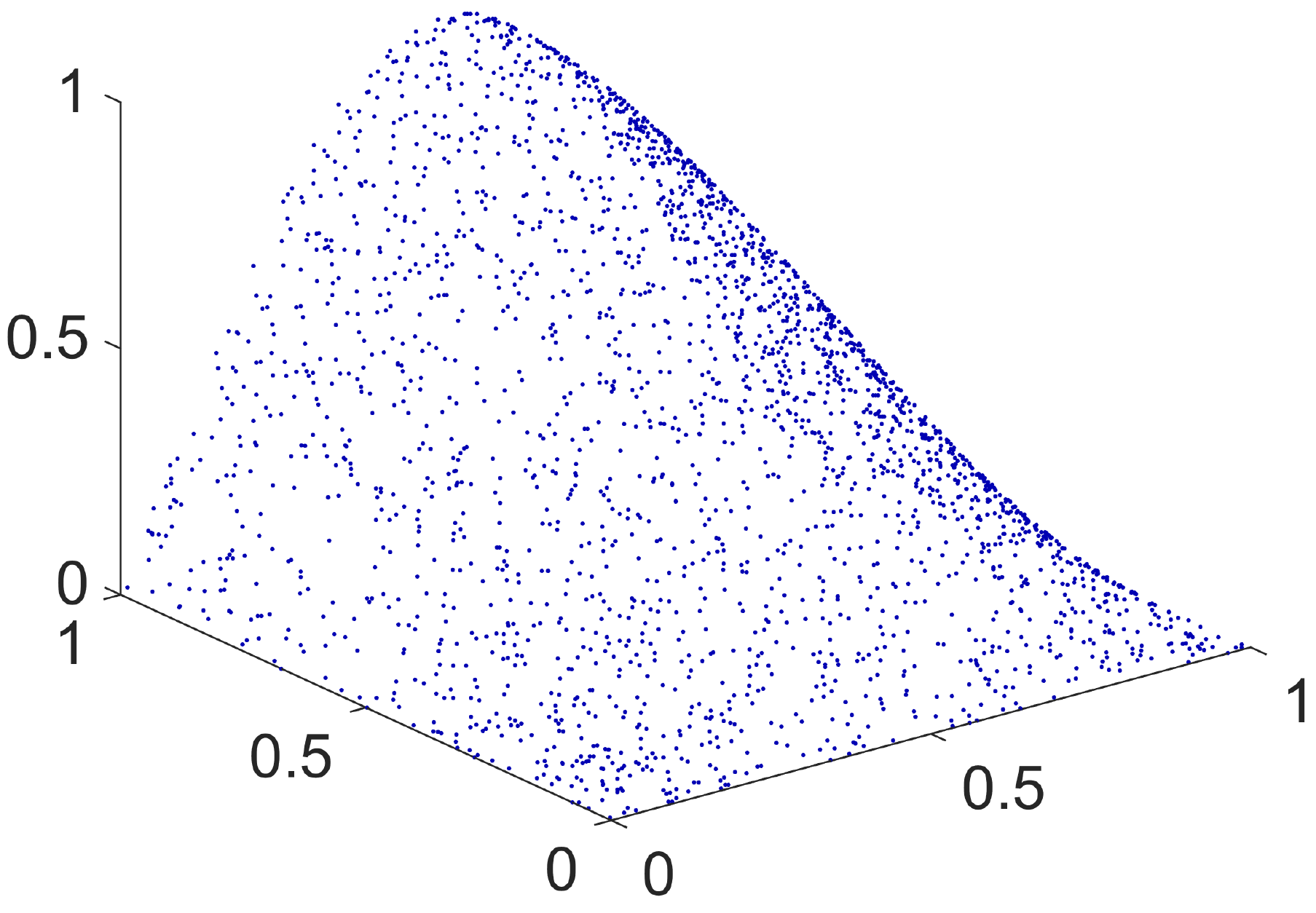}
\end{center}
\caption{\label{fig:invinterp2} Illustration of the 2D inverse interpolation
problem \texttt{PRinvinterp2} with \texttt{n} = 50.
Left:\ the true solution $x$ defined on a regular grid.
Right:\ the data $b$ defined on randomly scattered points.}
\end{figure}

Our test problem \texttt{PRinvinterp2} involves a regular 2D grid with
$N=\mathtt{n}^2$ grid points $(s_j^{\mathrm{G}},t_j^{\mathrm{G}})$ generated
by \texttt{meshgrid(linspace(0,1,n))}, and there are $M=N$
data points $(s_i,t_i)$ randomly distributed in $[0,1]\times [0,1]$.
The data values $b_i$ at the random points, as well as the true solution
$x_j$ at the grid points, are samples of the smooth
function $\phi(s,t) = \sin(\pi s) \, \sin(\pi/2\, t)$.
See Figure~\ref{fig:invinterp2} for an illustration; the figures are
generated with \texttt{PRshowx} and \texttt{PRshowb} after constructing
the default test data, that is with the following lines of MATLAB code:
\begin{verbatim}
   [A, b, x, ProbInfo] = PRinvinterp2;
   PRshowx(x, ProbInfo)
   PRshowb(b, ProbInfo)
\end{verbatim}
The default value of {\texttt{n}} is 128, but test data with other values of {\texttt{n}} can be easily
generated by specifying this value directly as an input to the function, e.g.,
\begin{verbatim}
   [A, b, x, ProbInfo] = PRinvinterp2(256);
\end{verbatim}
The {\tt options} structure has only one field, {\tt options.InterpMethod}, which
can be used to choose one of four different types of interpolation:\ nearest-neighbour,
linear (default), cubic, and spline.  For example, to use the default value
of {\texttt{n}} $= 128$, but the optional cubic interpolation, use the following code:
\begin{verbatim}
   options = PRset('InterpMethod', 'cubic');
   [A, b, x, ProbInfo] = PRinvinterp2(options);
\end{verbatim}
The forward computation (corresponding to multiplication with $A$)
is always done by means of MATLAB's \texttt{interp2} function, and thus
{\tt A} is not constructed explicitly as a matrix, but instead is represented
by a function handle. In this test problem we also provide the adjoint
operation, so the function handle {\tt A} can be used to compute
matrix-vector multiplications with both $A$ and $A^T$.
As was mentioned in the previous test problem,
all of the iterative
methods that we provide in the package do not need any additional
information from the user.
In case users would like to test their own iterative algorithms,
we recall that matrix-vector multiplications with $A$ and $A^T$ can be computed
using simple function calls.
For example,
$r = A^T(b-Ax)$ can be computed as
\begin{verbatim}
   r = A(b - A(x, 'notransp'), 'transp');
\end{verbatim}
In addition, to get the effective size of {\tt A} (that is,
if it could be constructed explicitly), use the MATLAB statement
\begin{verbatim}
   dims = A([], 'size');
\end{verbatim}
With the default value \texttt{n} = 128, the result is {\tt dims = [16384, 16384]}.
Because the adjoint operation (corresponding to $A^T$) is coded by us,
it is not necessary to use transpose-free methods with this test problem.

\subsection{Inverse Diffusion}
\label{sec:inversediffusion}

Many inverse PDE problems, such as parameter identification in electrical
impedance tomography, are nonlinear.
For this package we provide a simple \textit{linear} PDE test problem where
the solution is represented on a finite-element mesh and the
forward computation involves the solution of a time-dependent PDE\@.
The underlying problem is a 2D diffusion problem in the domain
$[0,T] \times [0,1] \times [0,1]$ in which the solution $u$ satisfies
  \begin{equation}
  \label{eq:PDE}
    \frac{\partial u}{\partial t} = \nabla^2 u
  \end{equation}
with homogeneous Neumann boundary conditions and a smooth function $u_0$
as initial condition at time $t=0$.
The forward problem maps $u_0$ to the solution $u_T$ at time $t=T$,
and the inverse problem is then to reconstruct the initial condition
from~$u_T$, cf.~\cite{MinGengRen13}.

\begin{figure}
\begin{center}
\includegraphics[width=0.45\textwidth]{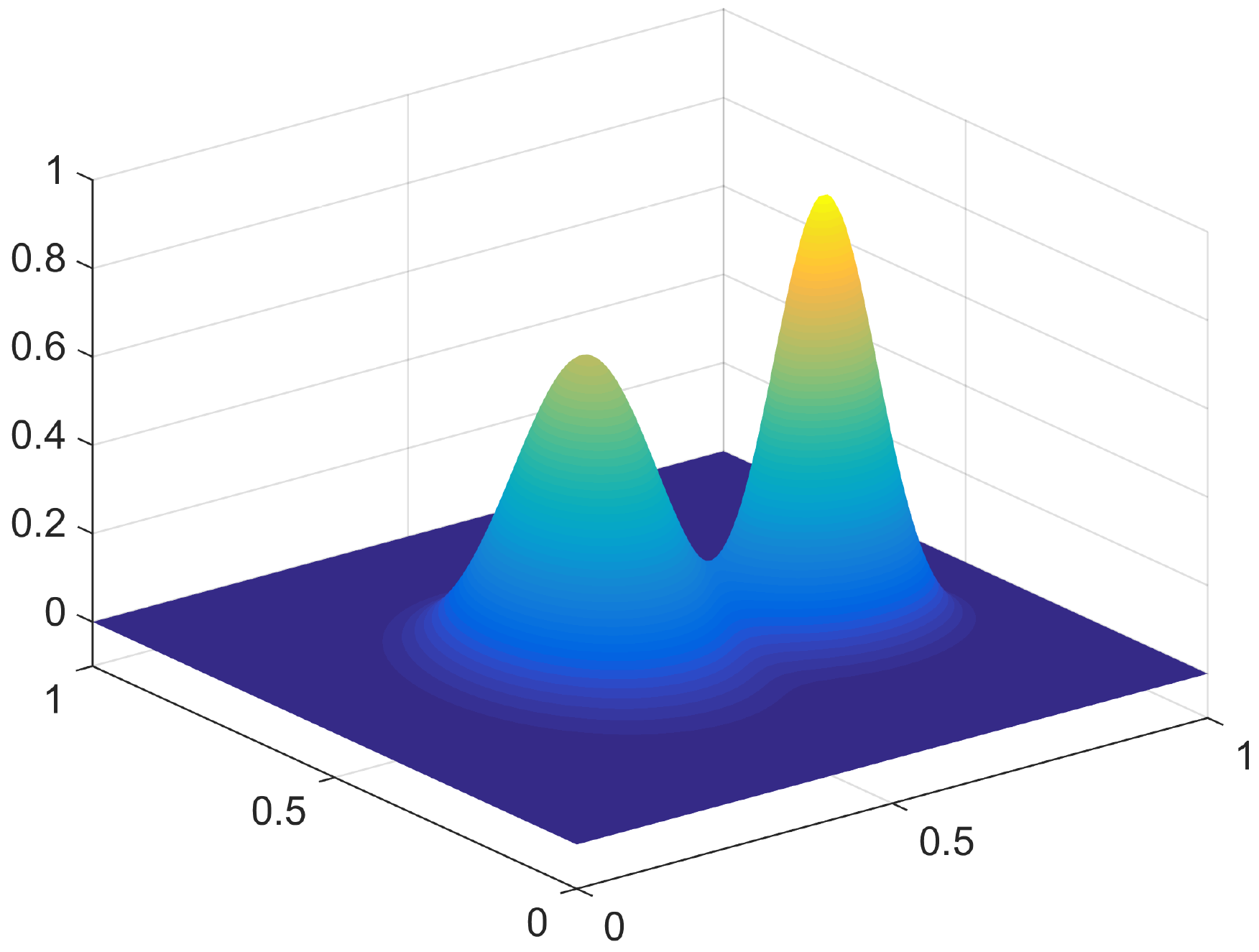} \hspace{4mm}
\includegraphics[width=0.45\textwidth]{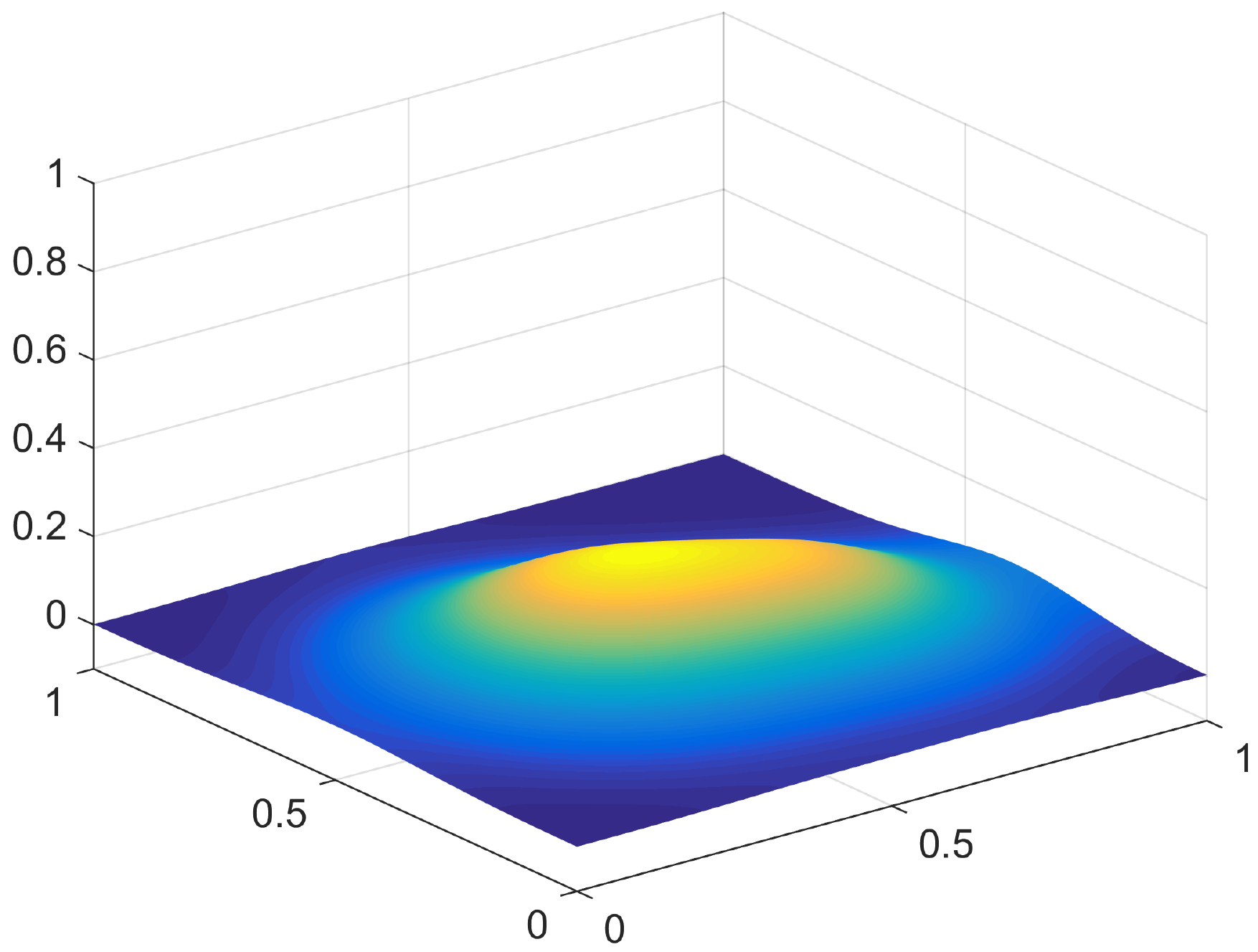}
\end{center}
\caption{\label{fig:diffusion} Illustration of the 2D inverse diffusion
problem \texttt{PRdiffusion} with \texttt{n} = 64.
Left:\ the true solution $x$ corresponding to the initial function~$u_0$.
Right:\ the data $b$ corresponding to the function~$u_T$ at time $T=0.01$.}
\end{figure}

We discretize the function $u$ on a uniform finite-element mesh with
$2(n-1)^2$ triangular elements; think of the domain as an
$(n-1) \times (n-1)$ pixel grid with two triangular elements in each pixel.
Then $u$ is represented by the $N = n^2$ values at the corners of the elements.
The forward computation -- represented by the function handle \texttt{A} --
is the numerical solution of the PDE (\ref{eq:PDE}), and it
is implemented by the Crank-Nicolson-Galerkin finite-element method.

The true solution $x$ and the right-hand side $b$ consist of the $N$
values of $u_0$ and $u_T$, respectively, at the corners of the elements;
see Figure~\ref{fig:diffusion} for an example, which was generated using
the statement:
\begin{verbatim}
   [A, b, x, ProbInfo] = PRdiffusion;
\end{verbatim}
This basic call uses the default input value of {\tt n} =128,
and sets default values for {\tt options}, which includes
\begin{itemize}
\item
{\tt options.TFinal} is the diffusion time (default is 0.01).
\item
{\tt options.Tsteps} is the number of time steps (default is 100).
\end{itemize}
As previously stated, {\tt A} is returned as a function handle that
can be used to perform matrix-vector multiplications.
As will be illustrated in Section~\ref{sec:Examples},
all of the iterative
methods that we provide in the package do not require any additional
information from the user.
In case users would like to use this test problem in their own iterative algorithms,
we recall that matrix-vector multiplications can be computed
using a simple function call.
For example,
to compute $r = b - A\, x$, use the MATLAB statement
\begin{verbatim}
   r = b - A(x, 'notransp');
\end{verbatim}
Note that if {\tt A} could be constructed explicitly, it would be an
$N \times N$ matrix, where $N = n^2$. The following MATLAB statement can be used to determine this information:
\begin{verbatim}
   dims = A([], 'size');
\end{verbatim}
thus, with the default value of {\tt n} = 128, {\tt N}= 16384.

As with other test problems in this package, an alternate value for {\tt n} can
be directly specified as an input to {\tt PRdiffusion}.  For example, if we
want to use {\tt n} = 256, and default values for {\tt options}, then
we can simply call {\tt PRdiffusion} as follows:
\begin{verbatim}
   [A, b, x, ProbInfo] = PRdiffusion(256);
\end{verbatim}
In addition, the default values for {\tt options}
can easily be changed using {\tt PRset}.
For example, if we want to use {\tt n} = 256, a
diffusion time of 0.3, and 50 time steps, then we type:
\begin{verbatim}
   options = PRset('TFinal', 0.3, 'Tsteps', 50);
   [A, b, x, ProbInfo] = PRdiffusion(256, options);
\end{verbatim}

\subsection{NMR Relaxometry}
\label{sec:nmr}

Nuclear Magnetic Resonance (NMR) relaxometry consists of reconstructing
a distribution of relaxation times associated with the probed material,
starting from a signal measured at given times.
Two-dimensional (2D) NMR relaxometry can be performed using particular
excitation sequences and acquisition strategies, so that the joint distribution
of the longitudinal and transverse relaxation times $T^1$ and $T^2$ can be recovered,
providing more chemical information about the probed material than its
one-dimensional analogous \cite{MCG12}.
2D NMR relaxometry is mathematically modeled using the following
Fredholm integral equation of the first kind
  \begin{equation}
  \label{eq:nmr}
    \int_{0}^{\widehat{T}^2} \int_{0}^{\widehat{T}^1}
    k(\tau^1, \tau^2, T^1, T^2) \, f(T^1,T^2) \, dT^1 \, dT^2
    = g(\tau^1, \tau^2)\,,
  \end{equation}
where $g(\tau^1, \tau^2)$ is the noiseless signal as a function of experiment
times $(\tau^1, \tau^2)$, and $f(T^1,T^2)$ is the density distribution function.
The kernel $k(\tau^1, \tau^2, T^1, T^2)$ in equation (\ref{eq:nmr}) is separable
{and given by}
  \[
    k(\tau^1, \tau^2, T^1, T^2)=\left(1-2\exp(-\tau^1/T^1)\right)\exp(-\tau^2/T^2)\,,
  \]
and, upon variable transformation, it can be regarded as a Laplace kernel.
Perturbations arising in 2D NMR relaxometry measurements are typically modeled
as Gaussian white noise.
Common techniques to regularize the inversion process include the incorporation
of box constraints and smoothness priors on the solution \cite{BBFLZ16}.

We discretize the integral in (\ref{eq:nmr}) using the the midpoint quadrature
rule with logarithmically equispaced nodes
  \[
    T^1_{1},T^1_{2},\dots,T^1_{n_1}\quad\mbox{and}\quad T^2_1,T^2_{2},\dots,T^2_{n_2}\,,
  \]
and considering a corresponding change of variables.
We then enforce collocation on the logarithmically equispaced sampled values
  \[
    \tau^1_{1},\tau^1_{2},\dots,\tau^1_{m_1}\quad\mbox{and}\quad \tau^2_1,\tau^2_{2},\dots,\tau^2_{m_2}\,,
  \]
so that equation (\ref{eq:nmr}) can be discretized as
  \begin{equation}
  \label{separable}
    A^1\,X\,(A^2)^T=B\,,
    \end{equation}
where
  \[
    \begin{array}{llll}
    A^1_{\ell_1, k_1}&=&1-2\exp(-\tau^1_{\ell_1}/T^1_{k_1})\,,&
      \quad \ell_1=1,\dots,m_1\,,\;k_1=1,\dots,n_1\,,\\
    A^2_{\ell_2, k_2}&=&\exp(-\tau^2_{\ell_2}/T^2_{k_2})\,,&
      \quad \ell_2=1,\dots,m_2\,,\;k_2=1,\dots,n_2\,,\\
    B_{\ell_1, \ell_2}&=&g(\tau^1_{\ell_1},\tau^2_{\ell_2})\,,&
      \quad \ell_1=1,\dots,m_1\,,\;\ell_2=1,\dots,m_2\,,\\
    X_{k_1, k_2}&=&f(T^1_{k_1},T^2_{k_2})\,,&\quad k_1=1,\dots,n_1\,,\;k_2=1,\dots,n_2\,.
    \end{array}
  \]
Equation (\ref{separable}) is a consequence of the fact that the kernel in
(\ref{eq:nmr}) is separable.
Taking $M=m_1m_2$ and $N=n_1n_2$, and defining $x\in\mathbb{R}^N$ and $b\in\mathbb{R}^M$
as the vectors obtained by stacking the columns of $X\in\mathbb{R}^{n_1\times n_2}$
and $B\in\mathbb{R}^{m_1\times m_2}$, respectively,
we obtain the linear system (\ref{eq:Axb}).
Due to the separability of {the kernel} $k$ the matrix $A$ has Kronecker structure,
$A=A_2\otimes A_1$;
we do not construct $A$ explicitly, but instead use a function handle
that implements matrix-vector multiplications with $A$ and $A^T$
{through the relation
$(A_2\otimes A_1)\, x = \mathrm{vec} \bigl( A^1\,X\,(A^2)^T \bigr)$.}
The function to construct this example is {\tt PRnmr} and, like other
applications in our package, is called using:
\begin{verbatim}
   [A, b, x, ProbInfo] = PRnmr(n, options);
\end{verbatim}
The input parameter {\tt n} specifies the size of the relaxation time distribution to be recovered,
and can either be an integer scalar, in which case
it is assumed that {\tt n} $ = n_1 = n_2$, or a vector {\tt n} $= [n_1,\; n_2]$.  The default
value is {\tt n} = 128.

As with the previous example, since {\tt A} is a function handle, to
compute \linebreak[4]$r = A^T(b-A\, x)$
(e.g., for users who want to use this problem to test their own iterative methods)
we can use the MATLAB statement
\begin{verbatim}
   r = A(b - A(x, 'notransp'), 'transp');
\end{verbatim}
and the effective size of {\tt A} can be obtained as
\begin{verbatim}
   dims = A([], 'size');
\end{verbatim}
With the default value of {\tt n} $=128$,
this returns {\tt dims = [65536, 16384]}.

\begin{figure}
\begin{center}
\includegraphics[width=0.45\textwidth]{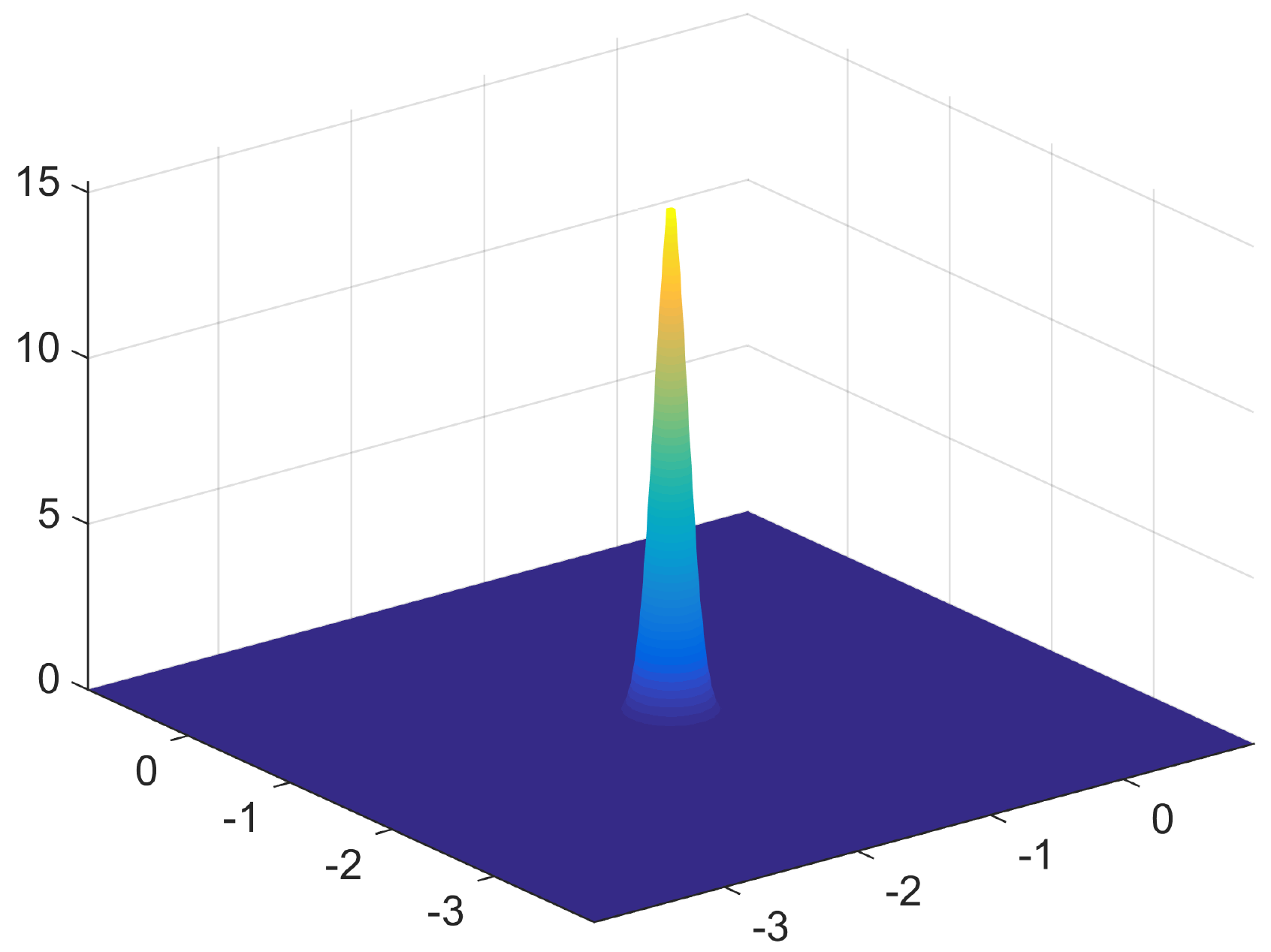} \hspace{4mm}
\includegraphics[width=0.45\textwidth]{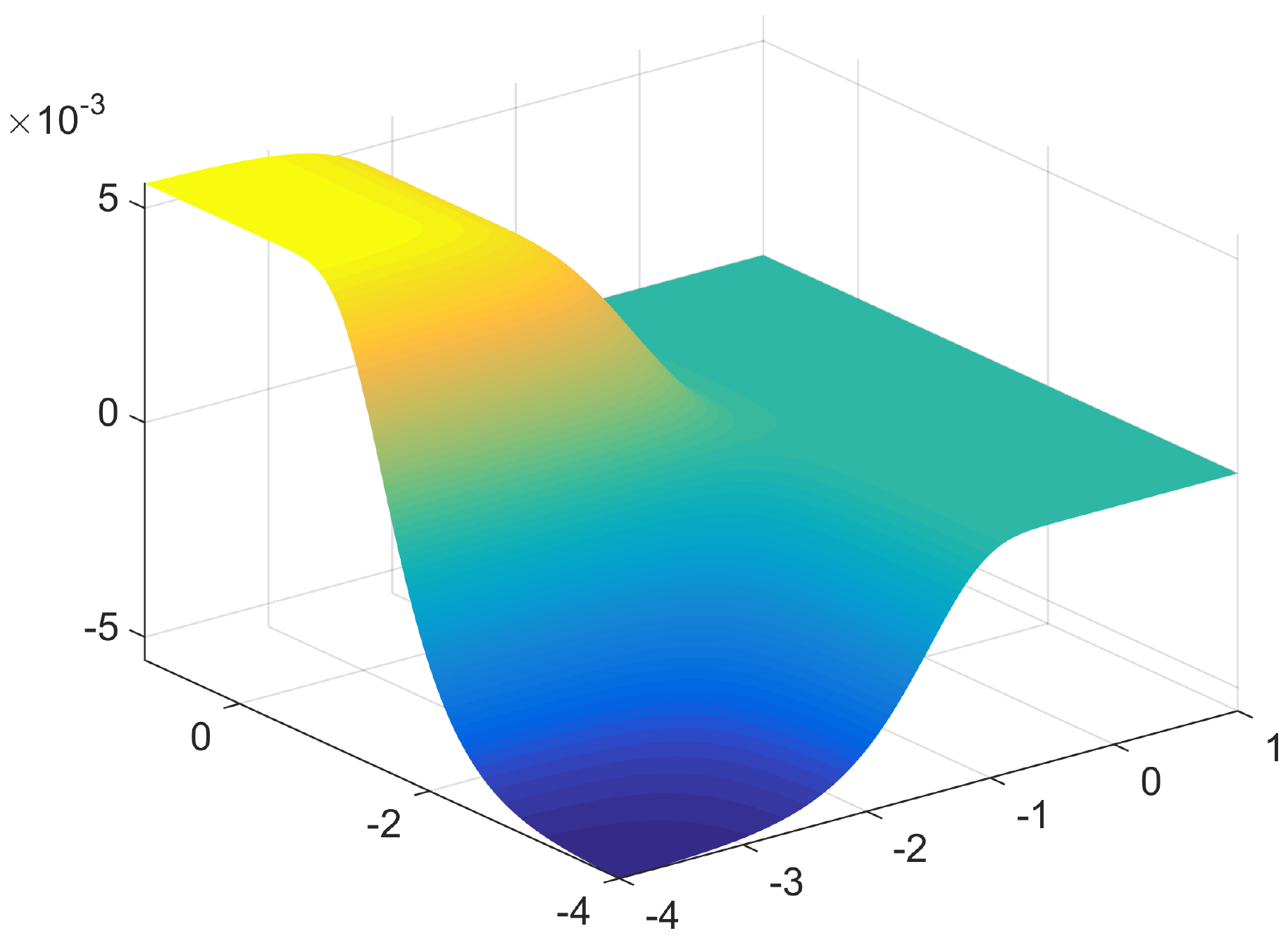}
\end{center}
\caption{\label{fig:nmr} Illustration of the NMR relaxometry
problem \texttt{PRnmr} with problem size {\texttt{n}} = 128.
Left:\ the true solution $x$ as a function of
$(\log_{10}(T^1),\log_{10}(T^2))$.
Right:\ the data $b$ as a function of
$(\log_{10}(\tau^1),\log_{10}(\tau^2))$.}
\end{figure}

The {\tt options} structure can be used to change other default parameters, including:
\begin{itemize}
\item
{\tt options.numData} is the number of acquired measurements, $m_1$ and $m_2$.
Default is $m_1 = 2n_1$ and $m_2 = 2n_2$.
\item
{\tt options.material} specifies the phantom for the relaxation time distribution,
which can be set to be {\tt 'carbonate'} (default),  {\tt 'methane', 'organic'} or {\tt 'hydroxyl'}.
The chosen {phantom} is returned as the vector {\tt x}.
\item
{\tt options.Tloglimits} is an array with two values, {\tt [Tlogleft, Tlogright]},
that define the limits for the logarithm of the relaxation times $T$, where
  \begin{center}
   \begin{tabular}{l}
   $T^1$ = {\tt logspace(Tlogleft, Tlogright, $n_1$),} \\
   $T^2$ = {\tt logspace(Tlogleft, Tlogright, $n_2$),}
    \end{tabular}
  \end{center}
The default is
$[-4, \; 1]$.
\item
{\tt options.tauloglimits} is an array with two values, {\tt [taulogleft, taulogright]},
that define the limits for the logarithm of the relaxation times $T$, where
  \begin{center}
   \begin{tabular}{l}
   $\tau^1$ = {\tt logspace(taulogleft, taulogright, $m_1$),} \\
   $\tau^2$ = {\tt logspace(taulogleft, taulogright, $m_2$),}
    \end{tabular}
  \end{center}
The default is
$[-4, \; 1]$.
\end{itemize}
The plots shown in Figure~\ref{fig:nmr} were obtained by using
{\tt PRshowx} and {\tt PRshowb} to display the data produced from
the most basic call to {\tt PRnmr} with default choices for {\tt n} and {\tt options}; that is,
\begin{verbatim}
   [A, b, x, ProbInfo] = PRnmr;
   PRshowx(x, ProbInfo)
   PRshowb(b, ProbInfo)
\end{verbatim}

\subsection{Tomography}
\label{sec:tomo}

Tomographic reconstruction problems come in many different forms,
and we provide three different types of such problems, which can
generate data using one of the following three statements:
\begin{verbatim}
   [A, b, x, ProbInfo] = PRtomo(n, options);
   [A, b, x, ProbInfo] = PRspherical(n, options);
   [A, b, x, ProbInfo] = PRseismic(n, options);
\end{verbatim}
All three problems are from the \textsc{AIR Tools II} package and we
refer to \cite{AIRTools} for more details and pictures of the
test images.  In each case the default value of {\tt n} determines
the size of {\tt x} (specifically, $x$ represents an $n \times n$ image).
The fields that can be specified in {\tt options} depend on the kind
of tomography taken into account.

\texttt{PRtomo} is used to generate test problems that model
\textit{X-ray attenuation tomography}, often referred to as computed
tomography (CT)\@.
This kind of tomography plays a large role in medical imaging
and materials science.
The data consists of measurements of the damping of X-rays that
penetrate the object and, to a good approximation, can be assumed to travel
along straight lines; see \cite{Buzug} for details and mathematical models.
The goal is then, from the data, to reconstruct an image of the
object's spatially varying attenuation coefficient.

Since each ray only traverses a small number of the total amount
of image pixels,
the matrix $A$ will be very sparse (for an $n \times n$ image there
are at most $2n$ nonzero elements per row of~$A$).
We provide two different measurement geometries (there are many
more in practice):
\begin{itemize}
\item
\texttt{options.CTtype} = \texttt{'parallel'} (default) gives a parallel-beam
tomography where, for each source-detector position angle, there are
a number of equidistantly spaced parallel X-rays.
This is the typical geometry in synchrotron X-ray measurements,
and it corresponds to the well-known Radon transform.
\item
\texttt{options.CTtype} = \texttt{'fancurved'} gives, for each
source-detector position angle, a fan of X-rays from a single
source to a curved detector, with an identical angular span between
all the rays.
This is often the case in large medical X-ray scanners.
\end{itemize}
The data is usually organized as an image called the sinogram, in which
each column consists of the data for one source-detector position angle.
The user can choose the number of angles, the number of rays per angle, etc.

\texttt{PRspherical} is used to generate test problems that model
\textit{spherical means tomography}.
This kind of tomography arises, e.g., in photo-acoustic imaging
based on the spherical Radon transform, where the data consists of integrals
along circles whose centers are located outside the object.
The goal is to reconstruct an image of the initial pressure distribution
inside the object (caused by a laser stimulation).
Since each circle only intersects a small number of image pixels,
the matrix $A$ is sparse.
The data is organized in a sinogram-like image whose
columns are the data for each circle center.
The user can choose the number of circle centers and
the number of concentric integration circles per center.

\texttt{PRseismic} is used to generate test problems that model
\textit{seismic travel-time tomography}.
This type of tomography uses measurements of
the delay of seismic waves to reconstruct an image of the slowness
(the reciprocal of the sound speed) in the domain of interest.
In our model problem, the sensors are located along two edges
of the image (corresponding to the surface and one bore hole)
while the wave sources are located along a third edge
(corresponding to another bore hole).
We provide two different models of the seismic wave:
\begin{itemize}
\item
\texttt{options.wavemodel} = \texttt{'ray'} (default)
corresponds to an assumption
that the wave frequency is infinite, such that the waves can be
well represented by straight lines (similar{ly} to X-ray tomography).
\item
\texttt{options.wavemodel} = \texttt{'fresnel'} corresponds to a model with
a finite wave frequency, where it is assumed that the wave is confined
to its first Fresnel zone -- a cigar-shaped domain with its endpoints
at the source and the detector.
\end{itemize}
Similar{ly} to X-ray tomography, in both
cases we obtain a sparse matrix, which is more sparse for the line model.
We organize the data in an image where each column contains
all the measurements from one source.
Since {\tt A} is a sparse matrix for all the tomography test problems,
standard MATLAB operators for multiplication, transpose, etc., can be used.

\subsection{The Severity of the Test Problems}

\begin{figure}
\begin{center}
\begin{tabular}{cc}
\includegraphics[width=0.47\textwidth]{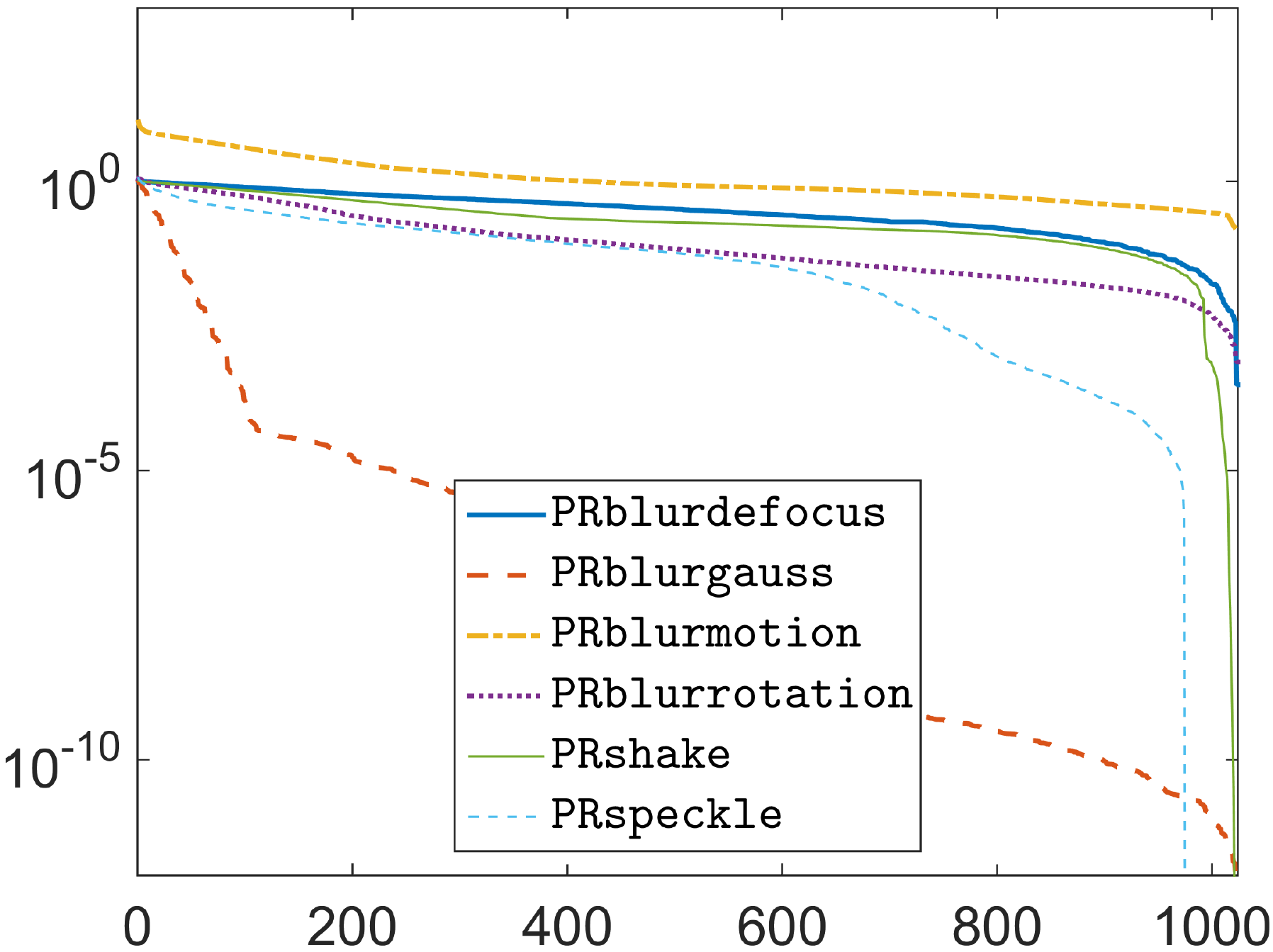} &
\includegraphics[width=0.47\textwidth]{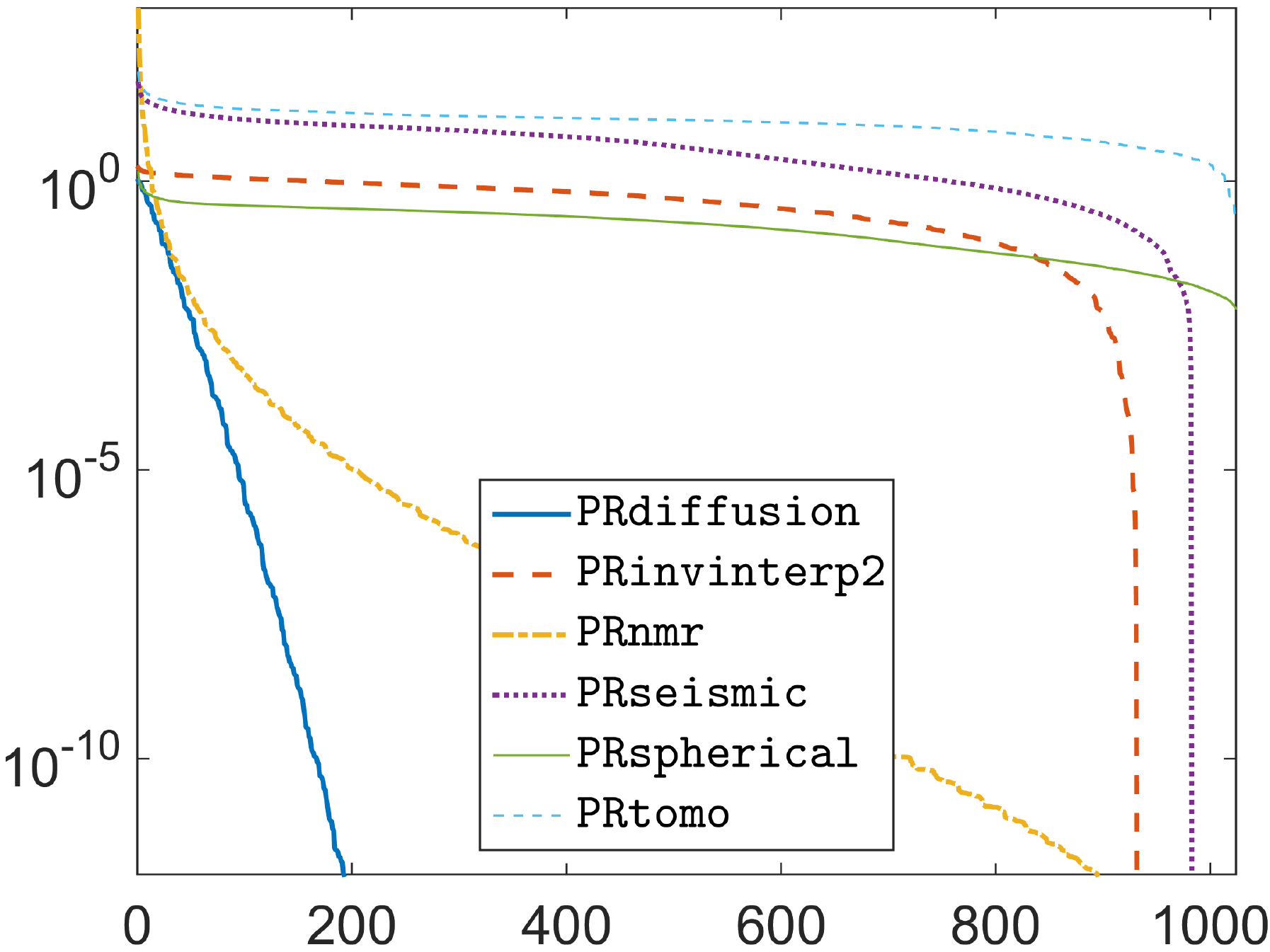}
\end{tabular}
\end{center}
\caption{\label{fig:decay} The singular values of the matrix $A$
for all 12 test problems in this package, for $n=32$ and using
default options.  \texttt{PRblurgauss}, \texttt{PRdiffusion} and
\texttt{PRnmr} are severely ill posed, the remaining problems
are mildly ill posed.}
\end{figure}

It is convenient to have a measure of the severity of the test problems
included in this package.
In linear algebra this is often measured by the condition number of the
matrix $A$, but the decay of the singular values of $A$ is a much better
measure of the severity of the underlying problem:\ the faster the decay
the severer the problem (and hence the larger the condition number);
see, e.g., \cite{Hansen}.

Figure \ref{fig:decay} shows the singular values of $A$ for all 12 test problems,
for the particular choice $n = 32$ and using default options.
Note that the fast decay of the singular values towards $N = n^2 = 1024$,
observed for most problems, is a discretization artifact.
The severely ill-posed problems are image deblurring with a Gaussian PSF,
the inverse diffusion problem, and the NMR relaxometry problem; the
remaining problems are mildly ill posed.
The help lines in the \texttt{PR\_\_\_} functions describe which problem
parameters affect the problem's severity.

\subsection{Adding Noise to the Data}

We also provide a function in order to make it easy to add noise to the data $b$:
\begin{verbatim}
   [bn, NoiseInfo] = PRnoise(b, NoiseLevel, kind);
\end{verbatim}
where the output \texttt{bn} = \texttt{b} + \texttt{noise} is the noisy data;
\texttt{noise} is the vector of perturbations, and it is available
within the output structure \texttt{NoiseInfo}.
As is the case with other \texttt{PR\_\_\_} functions, \texttt{PRnoise} can be called without specifying any input, in which case default values are used. The noise is scaled such that
  \begin{equation}
  \label{eq:RelNoiseLevel}
    \| \mathtt{noise} \|_2 / \| \mathtt{b} \|_2 =
    \mathtt{NoiseLevel}
  \end{equation}
with the default \texttt{NoiseLevel} = 0.01.
We provide three different kinds of noise {that can be easily obtained by setting the following options:}
\begin{itemize}
\item
\texttt{kind = 'gauss'} (default) gives Gaussian white noise,
$\mathtt{noise(i)} \sim \mathcal{N}(0,\sigma^2)$, with zero mean
and  with the standard deviation $\sigma$ chosen to satisfy
(\ref{eq:RelNoiseLevel}):
This noise is easily generated by means of:
\begin{verbatim}
   [bn, NoiseInfo] = PRnoise(b, NoiseLevel, 'gauss');
\end{verbatim}
%\begin{verbatim}
%  r = randn(M,1)
%   noise = ((RelNoiseLevel*norm(b(:)))/norm(r))*r
%   noise = reshape(noise,m,n)
%   bn = b + noise
%\end{verbatim}
\item
\texttt{kind = 'laplace'} gives Laplacian noise,
$\mathtt{noise(i)} \sim \mathcal{L}(0,\beta)$, with zero mean, and
the scale parameter $\beta$ is chosen to satisfy
(\ref{eq:RelNoiseLevel}).
This noise is easily generated by means of:
\begin{verbatim}
   [bn, NoiseInfo] = PRnoise(b, NoiseLevel, 'laplace');
\end{verbatim}
%\begin{verbatim}
%   r = rand(M,1)
%   r = sign(0.5-r).*(1/sqrt(2)).*log(2*min(r,1-r))
%   noise = ((RelNoiseLevel*norm(b(:)))/norm(r))*r
%   noise = reshape(noise,m,n)
%   bn = b + noise
%\end{verbatim}
\item
\texttt{kind = 'multiplicative'} gives a specific type of
multiplicative noise (often encountered in radar and ultrasound
imaging \cite{DoZe}) where each element
\texttt{bn(i)} equals \texttt{b(i)} times a random variable following
a Gamma distribution $\mathrm{\Gamma}(\kappa,\theta)$ with mean
$\kappa\,\theta = 1$ and the parameter $\kappa$ chosen such that
(\ref{eq:RelNoiseLevel}) is approximately satisfied:
\begin{verbatim}
   [bn, NoiseInfo] = PRnoise(b, NoiseLevel, 'multiplicative');
\end{verbatim}
%\begin{verbatim}
%   dfun = @(k) norm(b(:).*gamrnd(k,1/k,M,1)-b(:))/norm(b(:)) ...
%               - RelNoiseLevel
%   kappa = fzero(dfun,1/RelNoiseLevel^2)
%   bn = b.*reshape(gamrnd(kappa,1/kappa,m,1),m,n)
%   noise = bn - b
%\end{verbatim}
\end{itemize}
Information about the kind of generated noise and its level are available
within the output structure \texttt{NoiseInfo}.
Note that {\tt PRnoise} makes use of MATLAB's random number
generator functions, and thus to construct repeatable experiments
(i.e., to generate the same {\tt bn} for multiple experiments),
users should use MATLAB's {\tt rng} function to control the
seed of the random number generator before calling {\tt PRnoise}.

Other types of noise can be added by means of the function \texttt{imnoise}
from MATLAB's Image Processing Toolbox.

While Poisson noise is also a common type of noise in imaging,
it is not included in this package because it does not
conform to the use of \texttt{PRnoise}.
Specifically, in the presence of Poisson noise each noisy data element
\texttt{bn(i)} is an integer random variable following the
Poisson distribution $\mathcal{P}(\mathtt{b(i)})$, i.e.,
the noise-free data element
\texttt{b(i)} is both the mean and the variance of~\texttt{bn(i)}.
Hence, if we want to scale the ``noise'' vector
$\mathtt{noise} = \mathtt{bn}-\mathtt{b}$ then we can only do this
by scaling the noise-free data vector~$\mathtt{b}$ and the solution
vector \texttt{x} accordingly (the scaling factor can be found, e.g.,
by a simple fixed-point scheme.
Poisson noise can also be incorporated by means of \texttt{poissrnd}
from the Statistics Toolbox.

Another important type of noise, which arises in X-ray computed tomography,
can be referred to as ``log-Poisson'' (not a standard name).
Here the noisy elements of the right-hand side
in the linear model (\ref{eq:Axb}) are given by $b_i = \log(\tilde{d}_i)$
with $\tilde{d}_i \sim \mathcal{P}(d_i)$, where $d_i$ is the expected
photon count for the $i$th measurement.
It can be shown that $\log(\tilde{d}_i)$ approximately follows the
normal distribution $\mathcal{N}(\log(d_i),d_i^{-1})$
(corresponding to a quadratic approximation of the associated
likelihood function, cf.\ \cite{SaBo93}).
This provides a simple way to generate reasonably realistic noise
for X-ray tomography problems of the form (\ref{eq:Axb}) with the code:
\begin{verbatim}
   noise = randn(size(b))./sqrt(b);
   bn    = b + noise;
\end{verbatim}
Again, note that one must scale the noise-free $\mathtt{b}$ in order to
scale the relative noise level.

\section{Examples and Demonstrations}
\label{sec:Examples}

In this section we demonstrate the use of the iterative reconstruction methods
and the test problems by means of some numerical examples. Scripts to run these examples are
%These {\color{red}\sout{demonstrative}demonstration} examples are
available in \IRT\ with the naming convention \texttt{EX\_\_\_}.

\subsection{Solving 2D Image Deblurring Problems with CGLS and Hybrid Methods}
\label{sec:EXblur_cgls}

\begin{figure}[htbp]
\begin{center}
\begin{tabular}{cc}
\includegraphics[width=3.75cm]{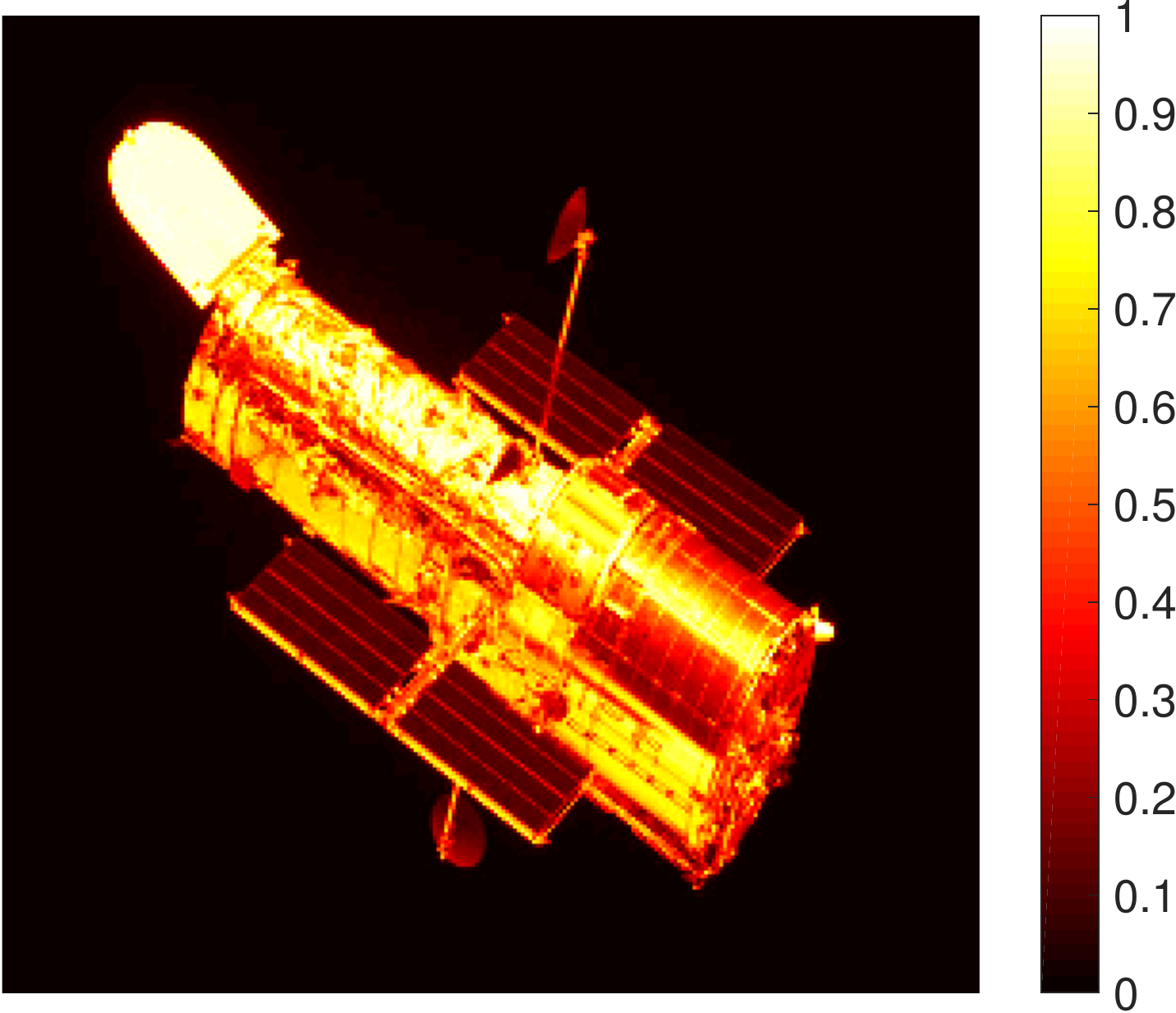} &
\includegraphics[width=3.75cm]{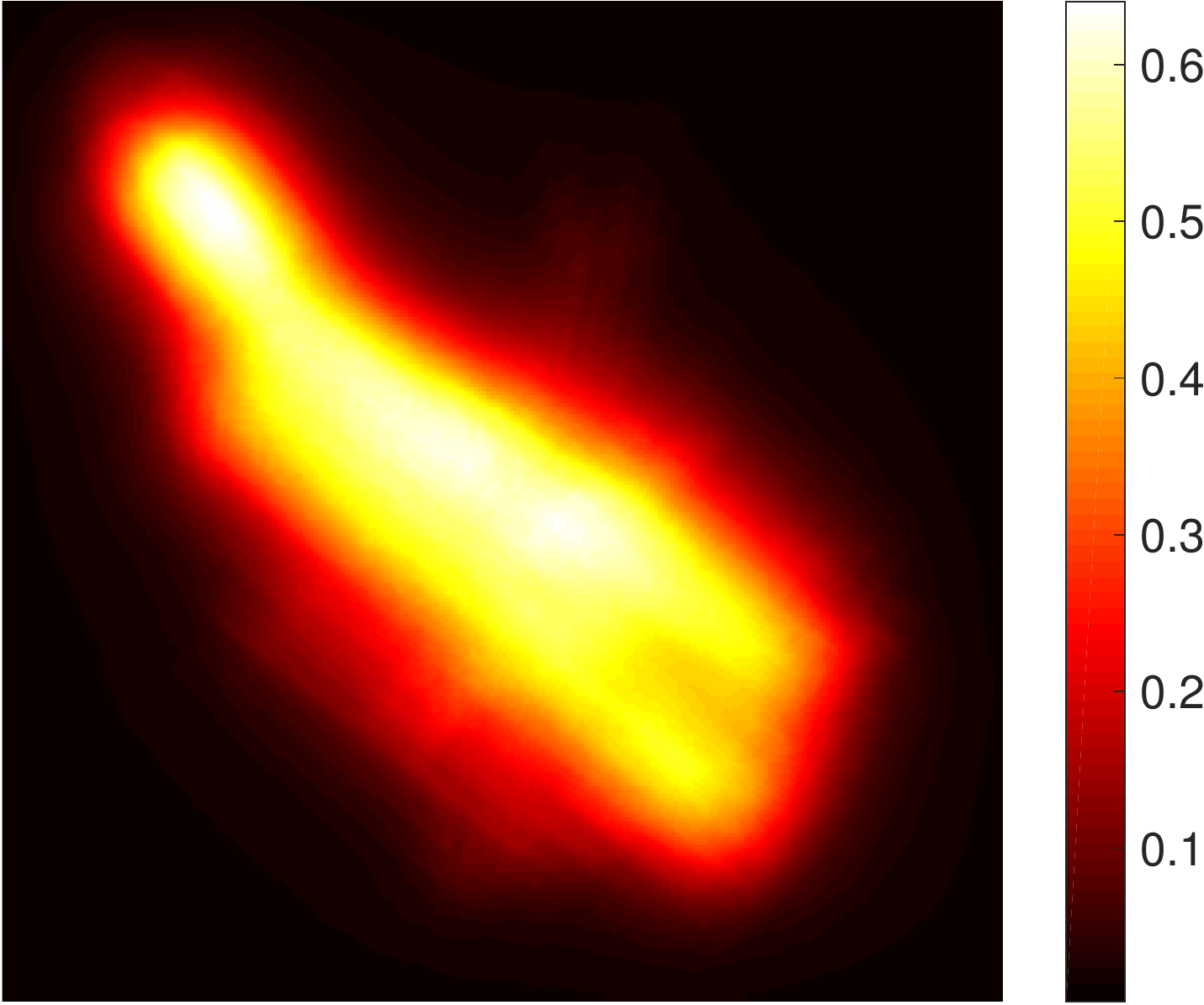}  \\
(a) & (b)
\end{tabular}
\caption{Image blurring test data for the example in Section~\ref{sec:EXblur_cgls}:
(a) true image, (b) blurred and noisy image.}
\label{fig:EXblur_cgls_images}
\end{center}
\end{figure}

\begin{figure}[htbp]
\begin{center}
\includegraphics[width=6cm]{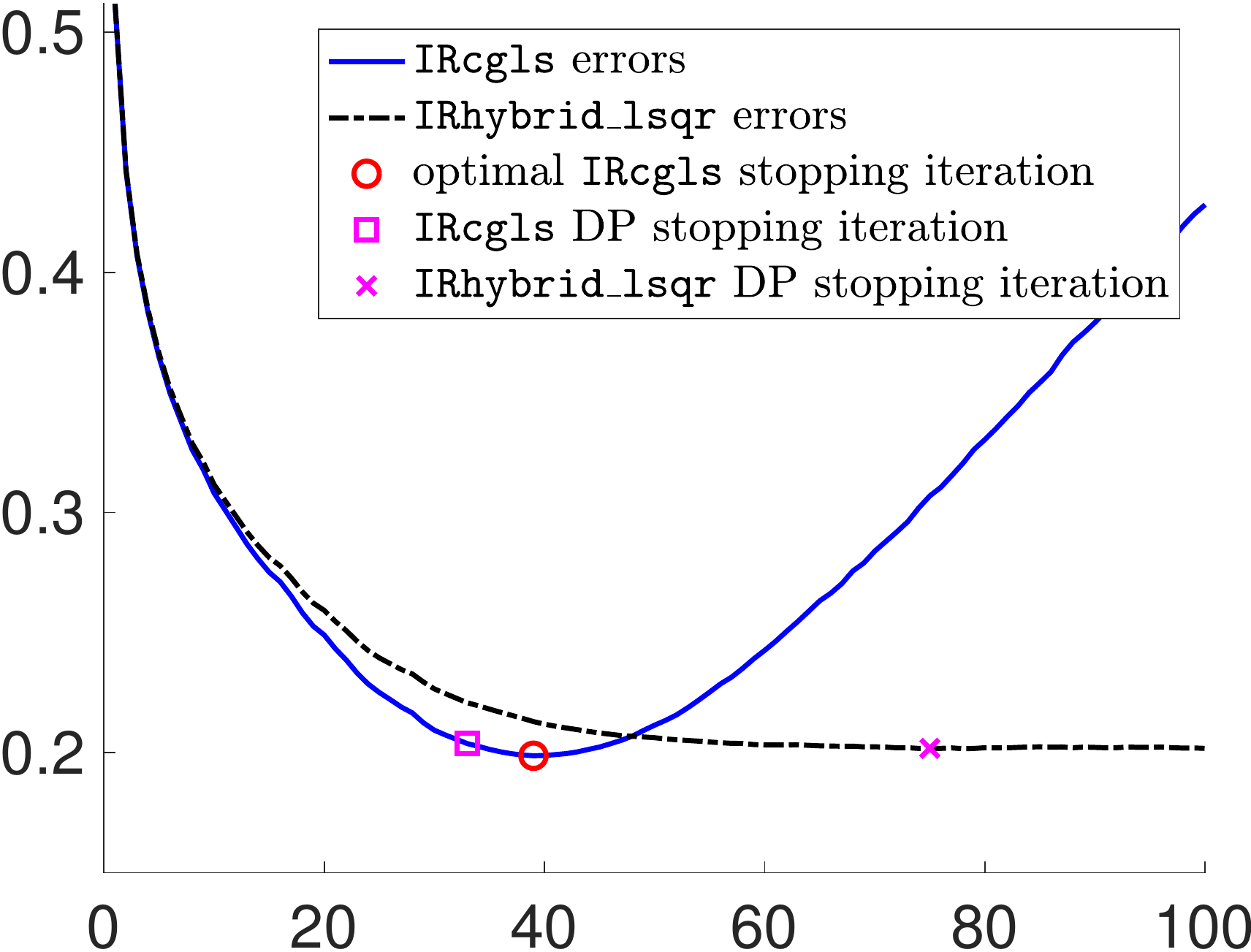} \\
\caption{Relative error plot for the image deblurring test problem from the example
in Section~\ref{sec:EXblur_cgls}.
The blue solid curve displays the iteration history (relative errors) of {\tt IRcgls},
the red circle marks iteration 39 where
the {\tt IRcgls} relative error is at its minimum value,
and the magenta square marks iteration 33 which is the
{\tt IRcgls} stopping iteration chosen by the discrepancy principle.
The black dashed curve is the iteration history
(relative errors) of {\tt IRhybrid\_lsqr}, and the magenta $\times$ marks
iteration 75 which is the {\tt IRhybrid\_lsqr} stopping iteration chose by
when using the weighted GCV {\tt 'wgcv'} parameter-choice method for the
projeted problem.}
\label{fig:EXblur_cgls_Enrm}
\end{center}
\end{figure}

\begin{figure}[htbp]
\begin{center}
\begin{tabular}{ccc}
\includegraphics[width=3.75cm]{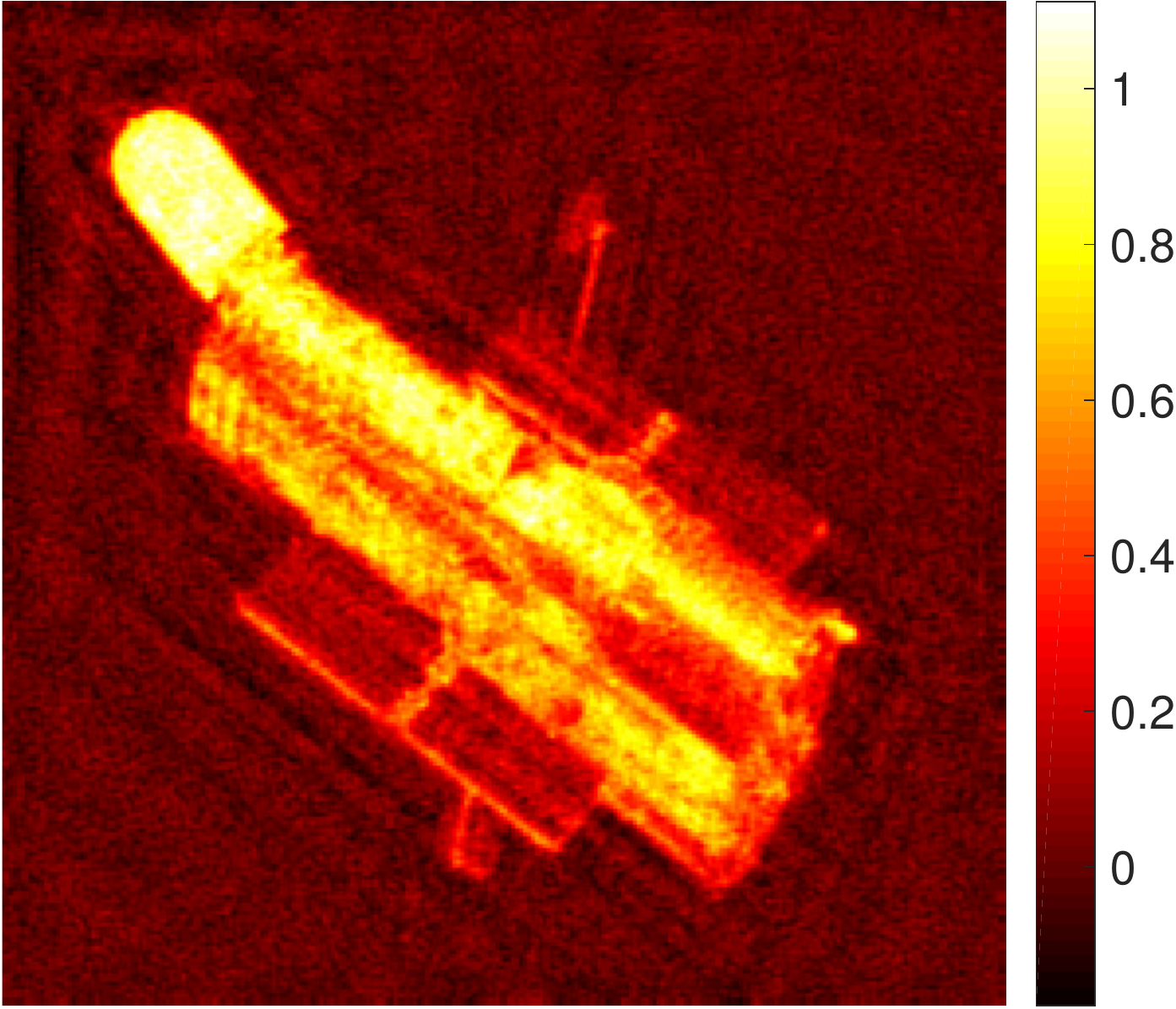} &
\includegraphics[width=3.75cm]{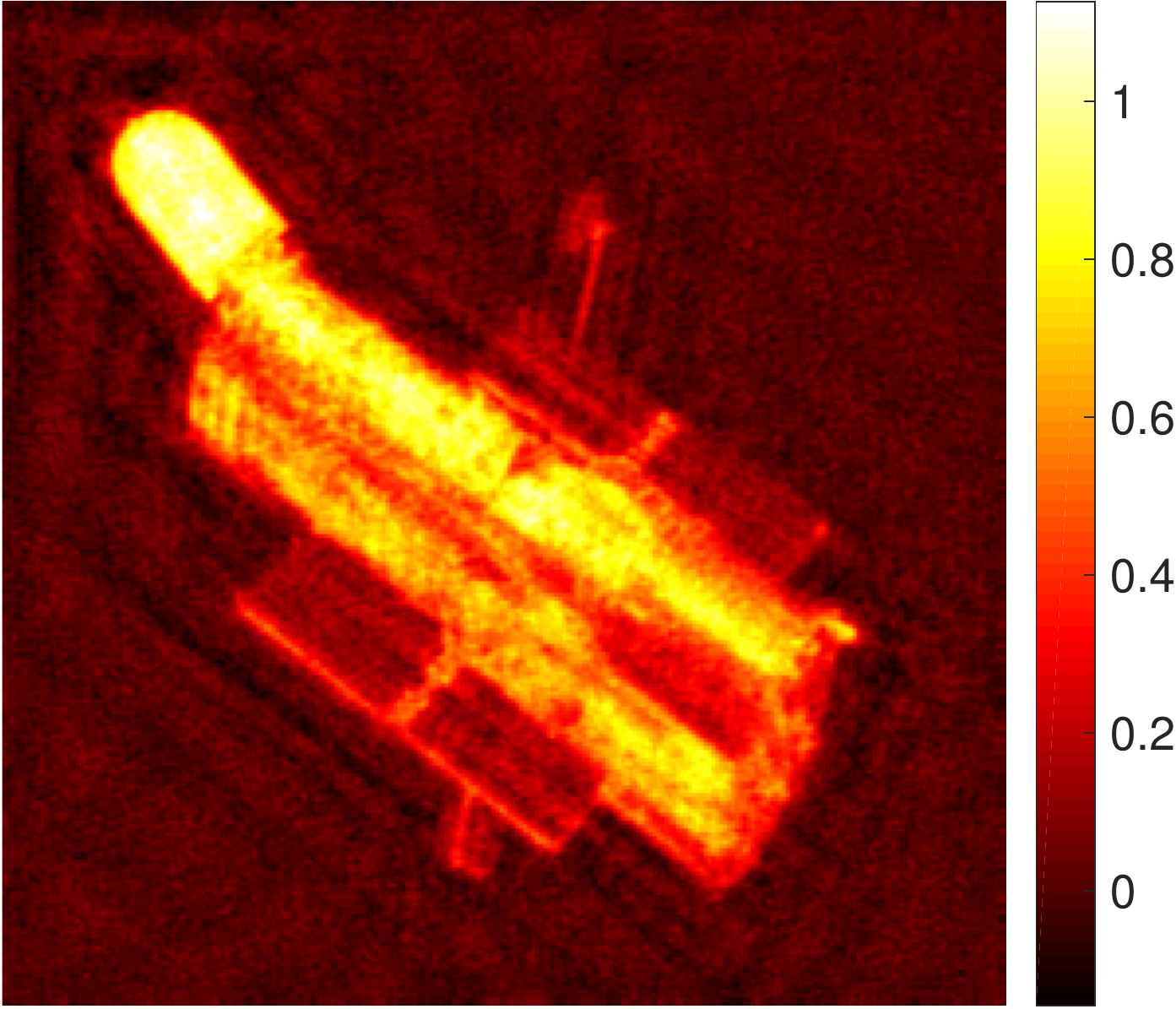} &
\includegraphics[width=3.75cm]{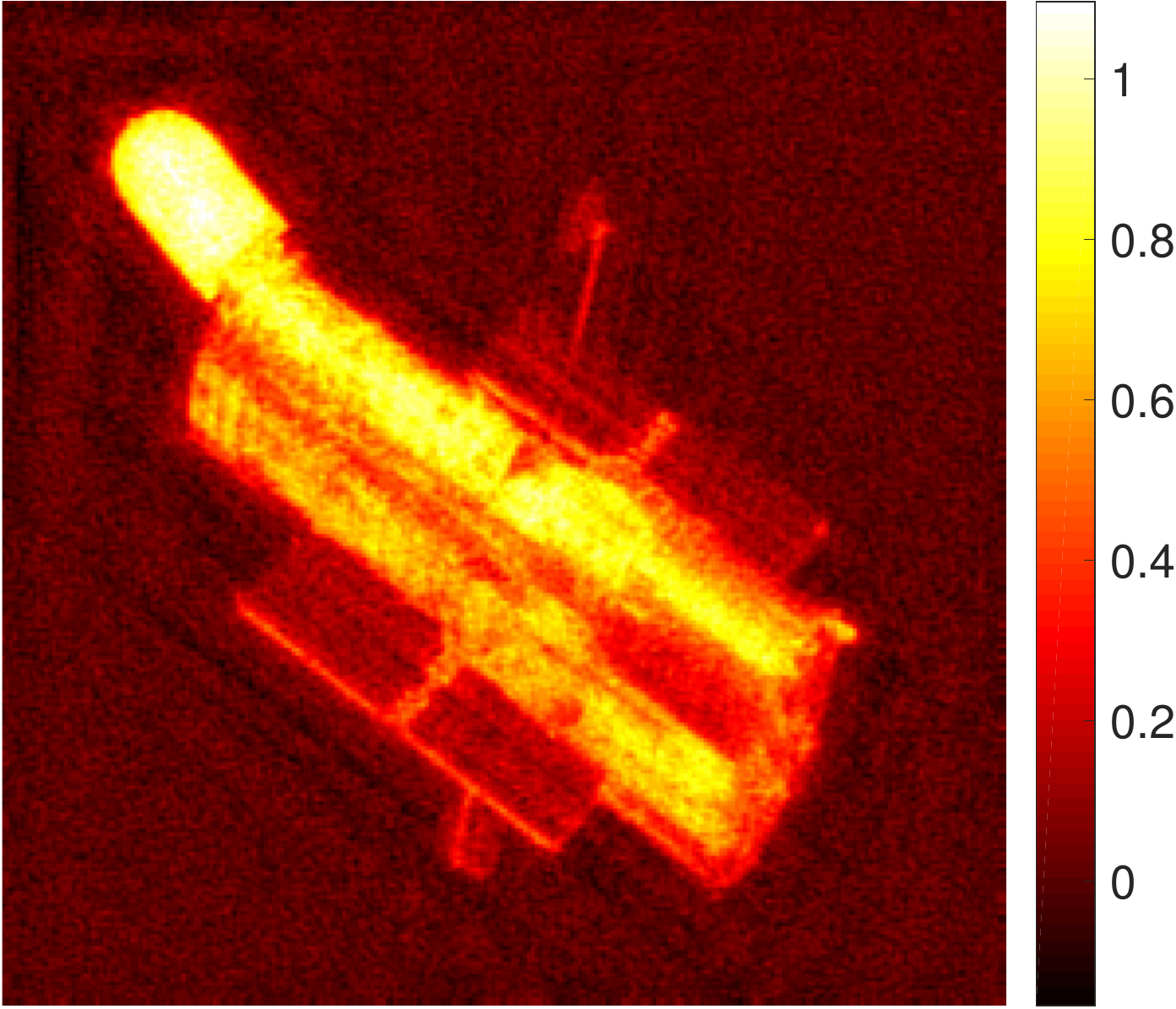}\\
(a) & (b) & (c)
\end{tabular}
\caption{Restored images from the example in Section~\ref{sec:EXblur_cgls}:
(a) restored image using 39 iterations of {\tt IRcgls}, (b) restored
image after 33 iterations of {\tt IRcgls}, (c) restored image after 75
iterations of {\tt IRybrid\_lsqr}.}
\label{fig:EXblur_cgls_results}
\end{center}
\end{figure}

Here we illustrate the use of {\tt IRcgls} and {\tt IRhybrid\_lsqr}
using the speckle image deblurring example
{\tt PRblurspeckle} described in Section~\ref{sec:blur}.
We use the default image size of \texttt{n} = 256 (i.e., the true and
blurred images have $256 \times 256$ pixels), the default true image
(Hubble Space Telescope), the default level of blurring (moderate), and
we add 1\% Gaussian noise; specifically, we generate the data using the following
lines of MATLAB code:
\begin{verbatim}
   NoiseLevel = 0.01;
   [A, b, x, ProbInfo] = PRblurspeckle;
   [bn, NoiseInfo] = PRnoise(b, 'gauss', NoiseLevel);
\end{verbatim}
Figure~\ref{fig:EXblur_cgls_images} shows the resulting true image {\tt x}
(Fig.~\ref{fig:EXblur_cgls_images}a) and
the blurred and noisy image {\tt bn} (Fig.~\ref{fig:EXblur_cgls_images}b).
We begin by running {\tt IRcgls} for
100 iterations, saving {only the iteration satisfying the stopping criterion}; we use the input {\tt options} to provide the
true solution to {\tt IRcgls}, so that we can investigate how the relative errors behave
at each iteration.  Specifically:
\begin{verbatim}
   options = IRset('x_true', x);
   [X, IterInfo_cgls] = IRcgls(A, bn, options);
\end{verbatim}
Note that in this example we do not specify a maximum number of
iterations, so the method uses the default value 100.
The output {\tt X} contains the solution at the final iteration;
in this example, convergence criteria are not satisfied,
so the method runs the full 100 iterations, and thus {\tt X} is the solution
at iteration 100.
Also note that because we specified the true solution \texttt{x} in {\tt options},
the relative errors $\| x^{(k)} - \mathtt{x} \|_2 / \| \mathtt{x} \|_2$
at each iteration are saved in the output structure, {\tt IterInfo\_cgls.Enrm}.
A plot of the relative errors can then be easily displayed as
\begin{verbatim}
   plot(IterInfo_cgls.Enrm)
\end{verbatim}
From this plot, which is shown by the blue solid curve in
Figure~\ref{fig:EXblur_cgls_Enrm},
we observe the well-known semi-convergence behavior of CGLS, and we
can also observe that the
smallest relative error occurs at iteration 39 (denoted by the
red circle in the plot). We refer to the solution where the relative error
is minimized as the ``best regularized solution.''
One feature of our iterative methods is that if the true solution is provided through
the options structure, then in addition to computing error norms, this best regularized
solution is also saved in {\tt IterInfo\_cgls.BestReg.X}, and the iteration where the error is smallest can be found
in {\tt IterInfo\_cgls.BestReg.It}.
Thus, if we want to
display this solution where the error is minimized, we can use {\tt PRshowx} as follows:
\begin{verbatim}
   PRshowx(IterInfo_cgls.BestReg.X, ProbInfo)
\end{verbatim}
This solution is shown in Figure~\ref{fig:EXblur_cgls_results}a.

Using the true solution to determine a stopping iteration is cheating,
but our implementations
can use other schemes that do not require knowing the true image.
For example, if we know the noise level in the data, then that information
 can be used along with the discrepancy
principle to determine a stopping iteration.
To do this, we simply need to change the options,
and run {\tt IRcgls}; specifically,
\begin{verbatim}
   options = IRset(options, 'NoiseLevel', NoiseLevel);
   [X, IterInfo_cgls_dp] = IRcgls(A, bn, options);
\end{verbatim}
We emphasize that the previously set options remain unchanged -- only the
one that is specified (in this example \texttt{'NoiseLevel'}) is changed.
Now the discrepancy principle terminates {\tt IRcgls} at iteration 33.
The relative error at this iteration is shown by the magenta square in
Figure~\ref{fig:EXblur_cgls_Enrm}.
It is well-known that the discrepancy principle tends to compute overly smooth
solutions, but this is not the case here where we know the exact error norm,
and we are able to compute the good restoration shown in
Figure~\ref{fig:EXblur_cgls_results}b.

We conclude this subsection by illustrating the use of one of the hybrid methods,
namely {\tt IRhybrid\_lsqr}.
This scheme enforces regularization at each iteration, and thus
avoids the semi-convergence behavior seen in {\tt IRcgls}.
In order to illustrate an approach that does not require an estimate
of the error norm, we use the weighted GCV {\tt 'wgcv'} parameter-choice
method (which is default) for the projected problem.
Specifically, if we use {\tt IRhybrid\_lsqr} and if we properly modify
the information about the regularization parameter choice in the previously
defined options,
\begin{verbatim}
    options = IRset(options, 'RegParam', 'wgcv');
    [X, IterInfo_hybrid] = IRhybrid_lsqr(A, bn, options);
\end{verbatim}
the method terminates at iteration 75 (which can be found from the output structure
{\tt IterInfo\_hybrid.its}).

If we want to show that {\tt IRhybrid\_lsqr} avoids the semi-convergence
behavior, we need to force the method to run more iterations, past the recommended
stopping iteration. We can do this by using an additional {\tt NoStop}
specification in the options.  That is,
\begin{verbatim}
   options = IRset(options, 'NoStop', 'on');
   [X_hybrid, IterInfo_hybrid] = IRhybrid_lsqr(A, bn, options);
\end{verbatim}
With the {\tt NoStop} option turned {\tt 'on'}, the iterations continue to the
default maximum of 100.
In this case, the vector {\tt X\_hybrid} is the solution at iteration 100,
but we also save the
solution at the recommended stopping iteration in the output structure,
{\tt IterInfo\_hybrid.StopReg.X}, and the iteration where the stopping criterion is
satisfied is saved in {\tt IterInfo\_hybrid.StopReg.It}.
Note that the field {\tt StopReg} is different than the field {\tt BestReg}:
the former stores information about the solution that satisfies the
stopping criterion; the latter stores information about the best computed
solution (and requires \texttt{x\_true} to be specified among the input options).
The relative errors for 100 iterations of {\tt IRhybrid\_lsqr} are shown
in the black dashed curve of Figure~\ref{fig:EXblur_cgls_Enrm},
with the recommended stopping iteration denoted by the magenta~$\times$.
The solution at this
recommended stopping iteration is shown in Figure~\ref{fig:EXblur_cgls_results}c.

The code used to generate the test problem and results described in this example is
provided in our package in the script {\tt EXblur\_cgls\_hybrid.m}.

\subsection{Solving the 2D Inverse Interpolation Problem with Priorconditioned CGLS}

\begin{figure}[htbp]
\begin{center}
\begin{tabular}{ccc}
\includegraphics[width=3.78cm]{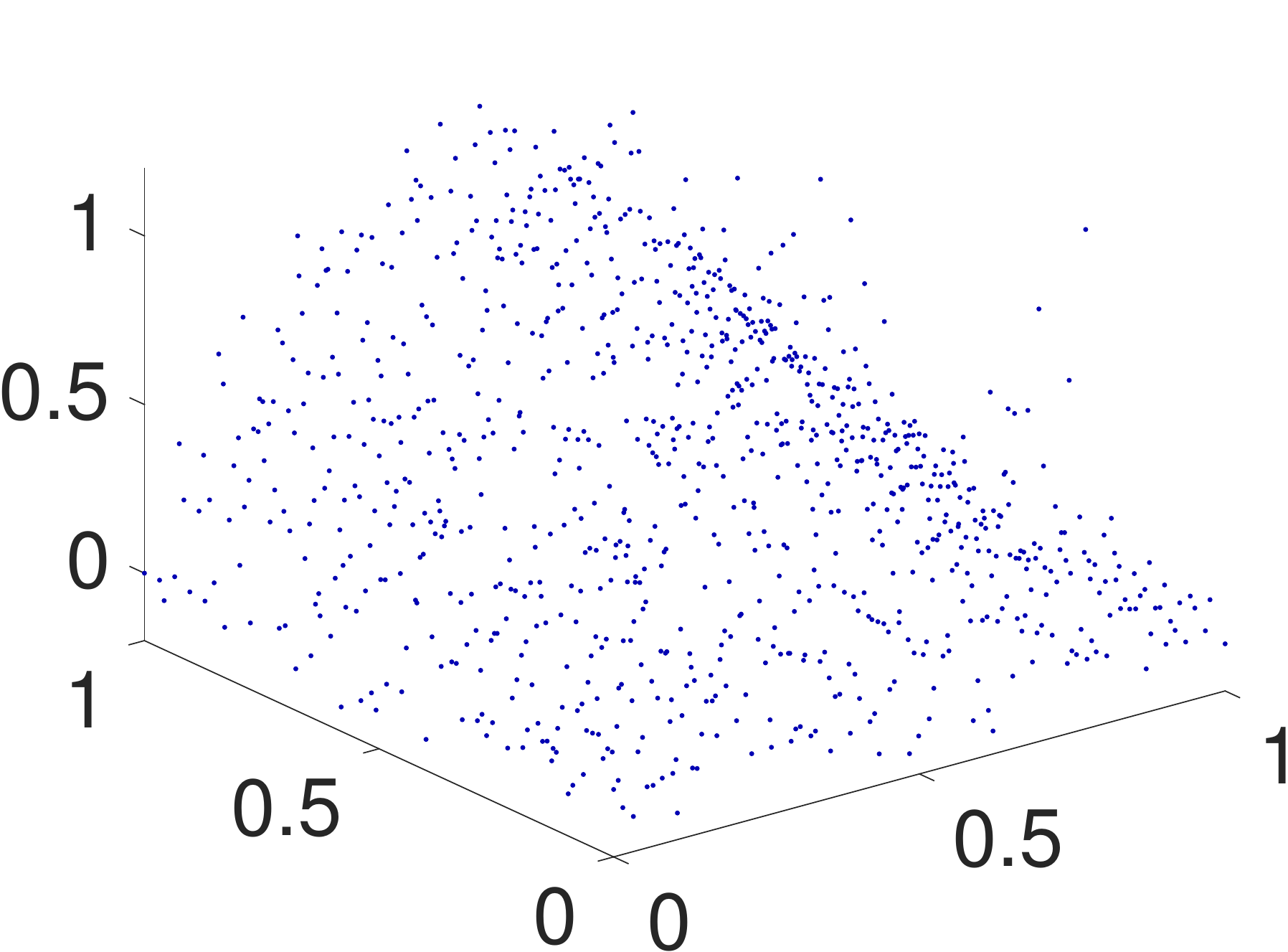} &
\includegraphics[width=3.78cm]{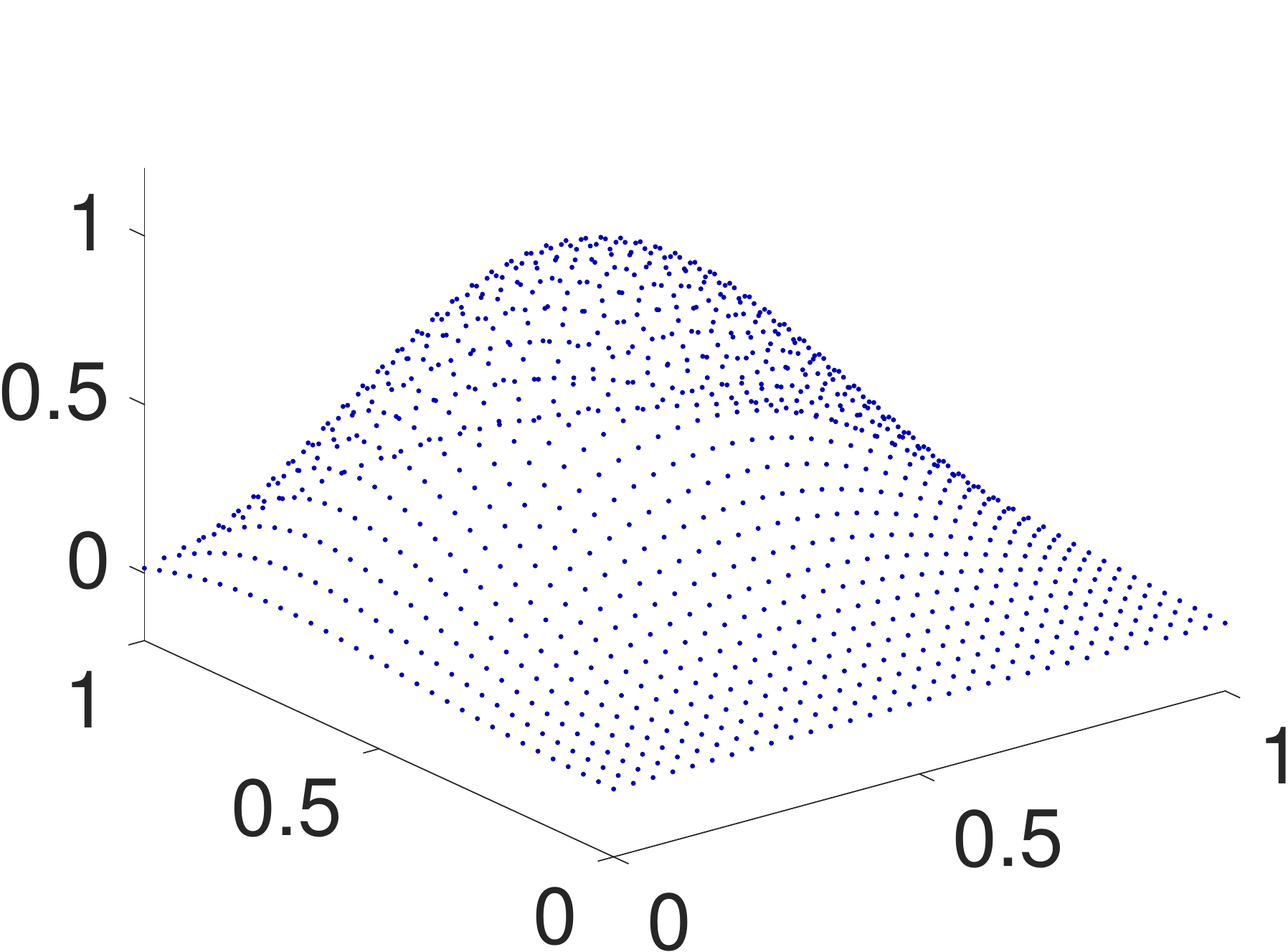} &
\includegraphics[width=3.78cm]{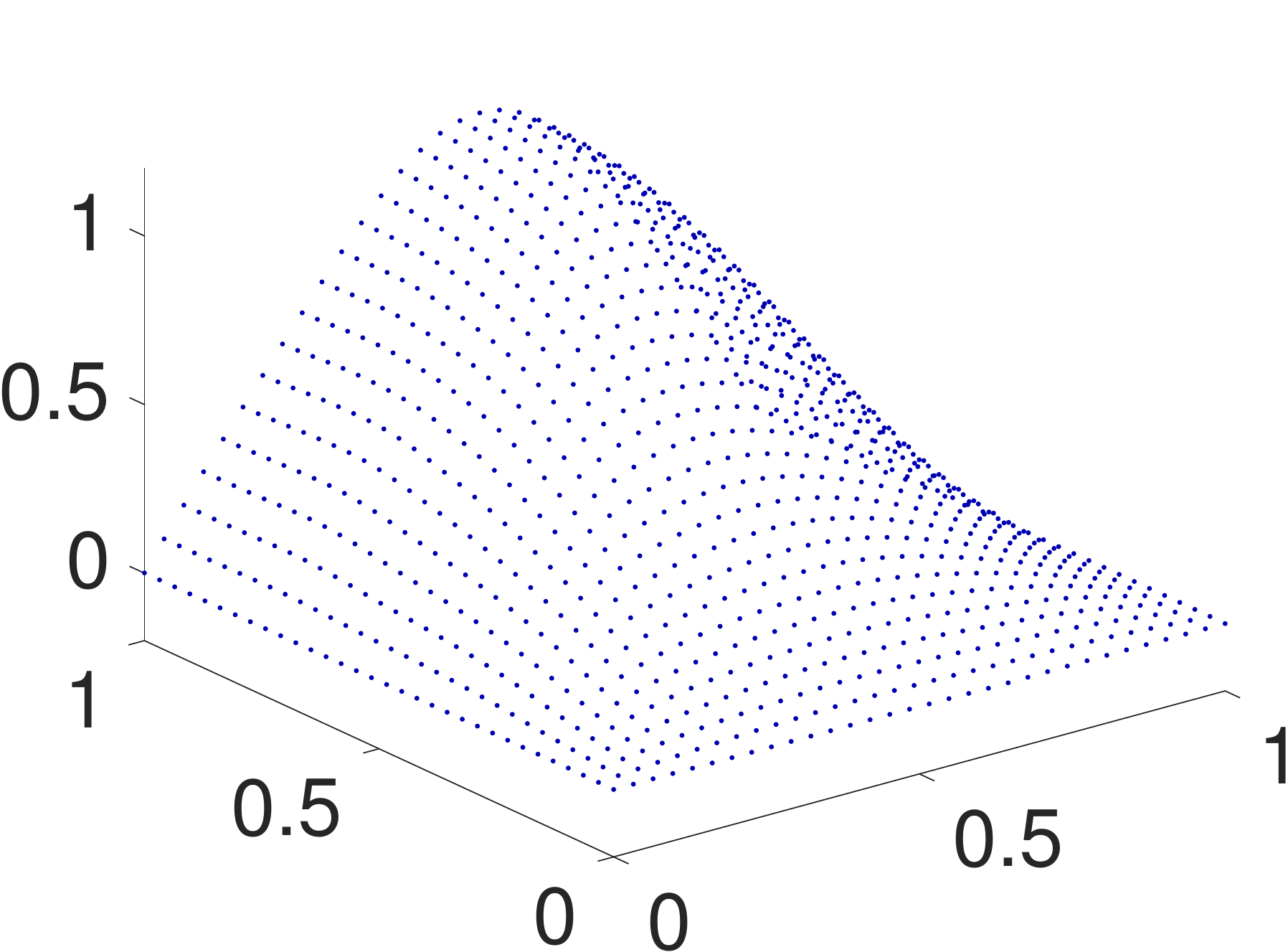}\\
(a) & (b) & (c)
% width=3.25cm % width=3.75cm
\end{tabular}
\caption{\label{fig:invinterp2Demo} Illustration of the solution of the
2D inverse interpolation problem \texttt{PRinvinterp2} with \texttt{n} = {32}
by means of \texttt{IRcgls} and with three different regularization
matrices: (a) the identity gives a very noisy solution, (b) the
2D Laplacian with zero boundary conditions everywhere gives a smooth
but incorrect solution, (c) enforcing instead a zero derivative on the
boundary where the solution is nonzero gives a good approximate solution.}
\end{center}
\end{figure}

Here we use the 2D inverse interpolation test problem {\tt PRinvinterp2}
to illustrate how to use \textit{prior-conditioning} in {\tt IRcgls}.
We begin by generating the test problem using {\tt n} $= 32$, and add 5\% Gaussian noise:
\begin{verbatim}
   n = 32;
   [A, b, x, ProbInfo] = PRinvinterp2(n);
   bn = PRnoise(b, 0.05);
\end{verbatim}
The true solution {\tt x} and data {\tt b} were already shown in Section~\ref{sec:inverseinterpolation},
Figure~\ref{fig:invinterp2}; the data {\tt bn} looks very similar to {\tt b}.
We attempt to solve this problem with three different versions of CGLS:
\begin{itemize}
\item
Standard CGLS, using the statement:
\begin{verbatim}
   K = 1:200;
   [X1, IterInfo1] = IRcgls(A, bn, K);
\end{verbatim}
Unfortunately, without providing additional information, CGLS cannot recognize
an appropriate stopping iteration, and the final computed solution
is a poor approximation; see Figure~\ref{fig:invinterp2Demo}a.
\item
Priorconditioned CGLS with \texttt{options.RegMatrix = 'Laplacian2D'}
which enforces zero boundary conditions everywhere.  This can be computed using
\begin{verbatim}
   options.RegMatrix = 'Laplacian2D';
   [X2, IterInfo2] = IRcgls(A, bn, K, options);
\end{verbatim}
In this case, CGLS finds a smoother solution that somewhat resembles the
exact solution in half of the domain.
But in the other half, the solution (while still smooth) is incorrect due to
the zero boundary condition at that edge where the exact solution is nonzero.
This is clearly seen in Figure~\ref{fig:invinterp2Demo}b.
\item
We can also create our own prior-conditioning matrix $L$.
Specifically we construct a matrix $L$ that is
similar to the 2D Laplacian, except we enforce a \textit{zero derivative} on one
of the boundaries;
\begin{verbatim}
   L1 = spdiags([ones(n,1),-2*ones(n,1),ones(n,1)],[-1,0,1],n,n);
   L1(1,1:2) = [1,0]; L1(n,n-1:n) = [0,1];
   L2 = L1; L2(n,n-1:n) = [-1,1];
   L = [ kron(speye(n),L2) ; kron(L1,speye(n)) ];
   L = qr(L,0);
   options.RegMatrix = L;
   [X3, IterInfo3] = IRcgls(A, bn, K, options);
\end{verbatim}
In this case, the iteration terminates after $k=10$ iterations, and
because the prior-conditioner enforces correct boundary conditions,
we obtain a very good computed approximation; see Figure~\ref{fig:invinterp2Demo}c.
\end{itemize}

In this example, we rely on the default normal equations
residual for the stopping rule,
$\| A^T (b-A\ x^{(k)}) \|_2 / {\| A^Tb \|_2} \leq \mathtt{options.NE\_Rtol}
= 10^{-12}$, where $x^{(k)}$ is the computed approximate solution at iteration $k$.
We also remark that in each of the calls to {\tt IRcgls},
the third input argument {\tt K = 1:{200}} is used to
request that the methods return all solution iterates in {\tt X1},
{\tt X2} and {\tt X3}.
For example, since the first call to {\tt IRcgls} runs all {200} iterations,
{\tt X1} is an array of size ${1024} \times 200$, but since the other
two calls only needed 5 iterations, {\tt X2} and {\tt X3} are arrays of
size ${1024} \times 5$.  This can be very useful if one wants to
view solutions at earlier iterations. For example, it would be
very easy to see how the solution at iteration 5
of the first call to {\tt IRcgls} compares with
second two calls, e.g.,
\begin{verbatim}
   PRshowx(X1(:,5), ProbInfo)
\end{verbatim}
However, requesting all the solution
iterates can lead to a large amount of storage,
especially when solving very large problems, so we caution
users to use this capability wisely.  For example, {\tt K} can be
any set of integers, such as {\tt K = [1, 10:{10:200}]}, which would
return solutions at iterations $1,\,10,\, 20,\, \ldots, \, {200}$.

The code used to generate the test problem and results described in this example is
provided in our package in the script {\tt EXinvinterp2\_cgls.m}.

\subsection{Solving the 2D Inverse Diffusion Problem with RRGMRES}

\begin{figure}
\begin{center}
\begin{tabular}{ccc}
\includegraphics[width=3.75cm]{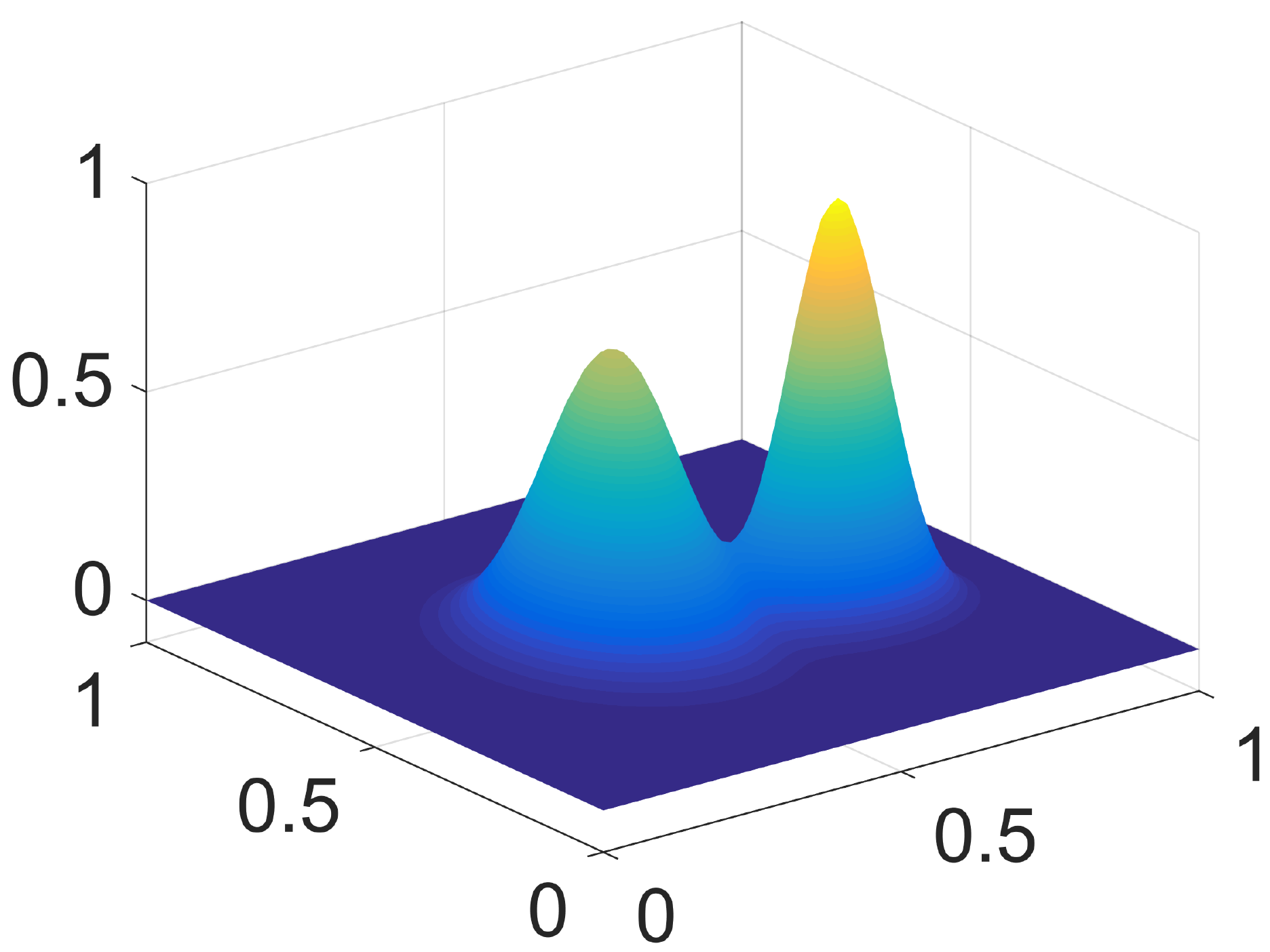} &
\includegraphics[width=3.75cm]{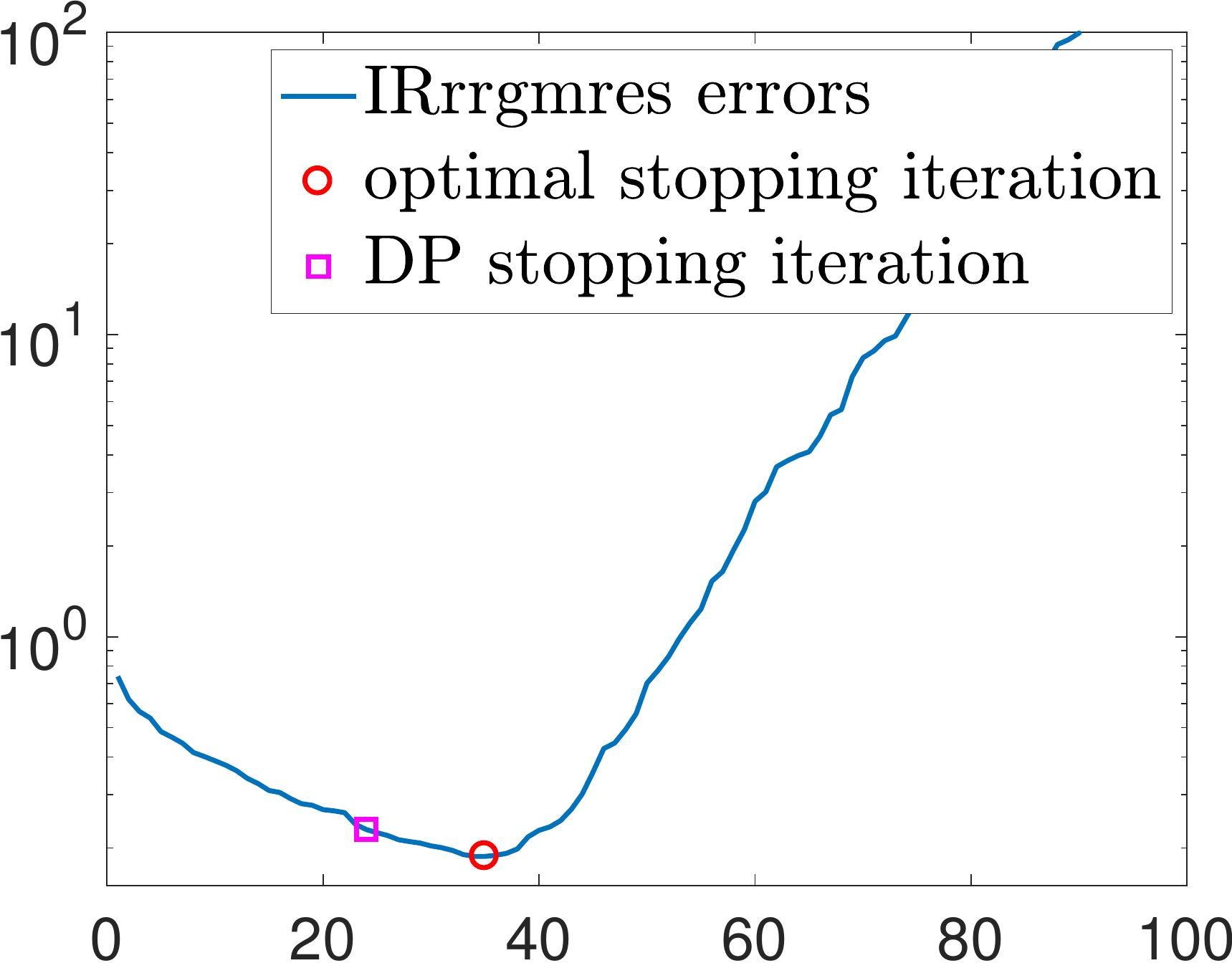} &
\includegraphics[width=3.75cm]{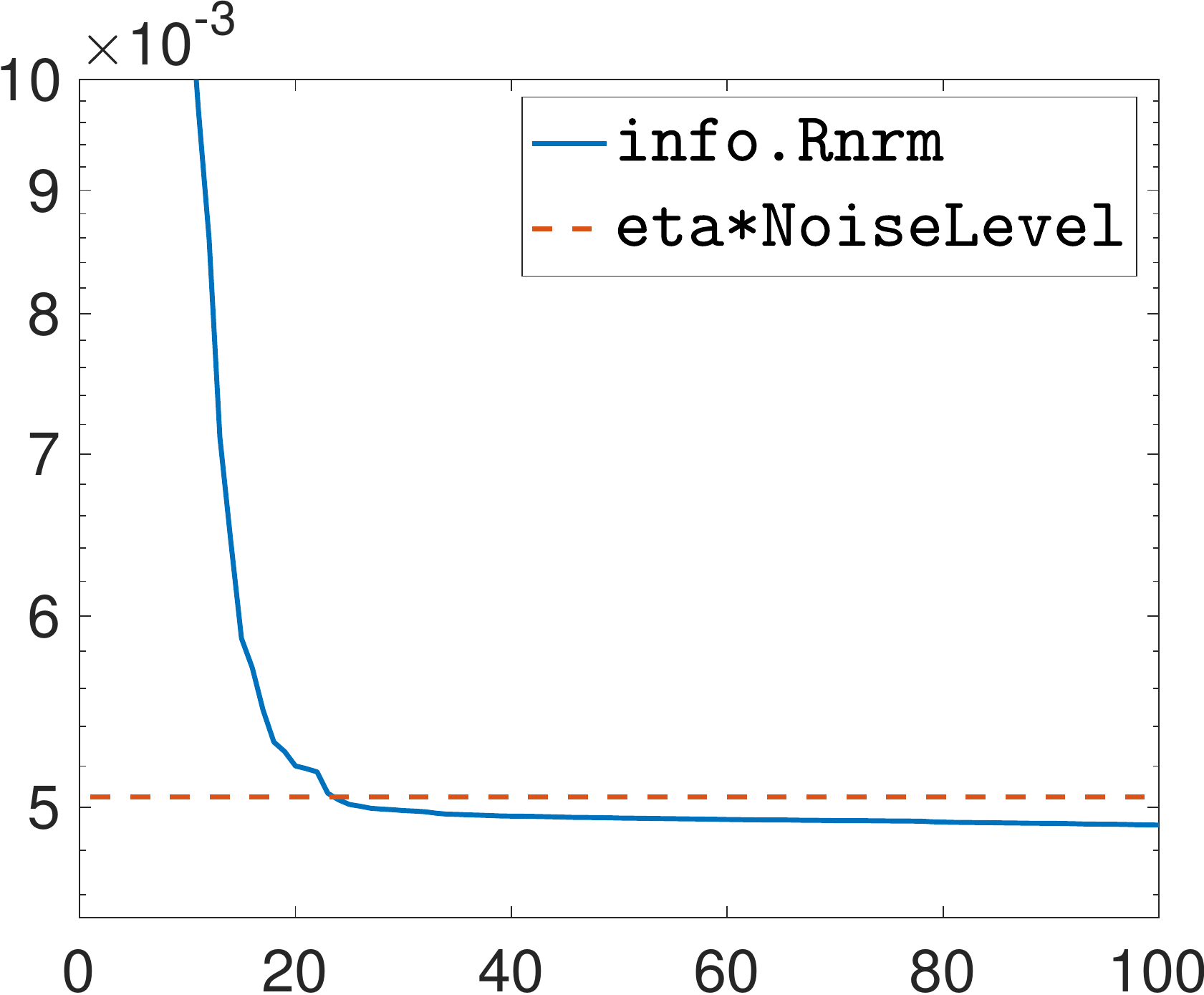}\\
(a) & (b) & (c)\\
\includegraphics[width=3.75cm]{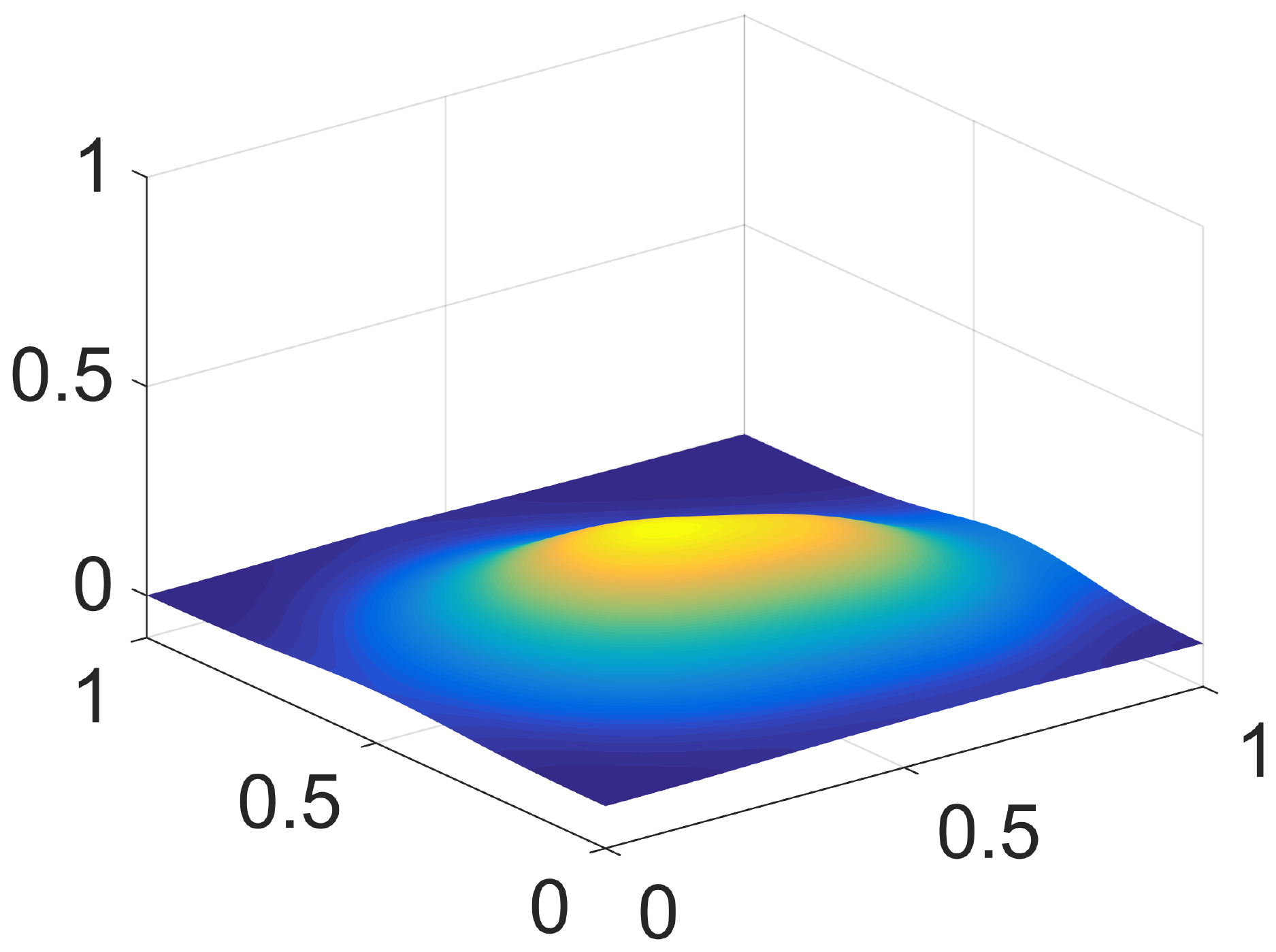} &
\includegraphics[width=3.75cm]{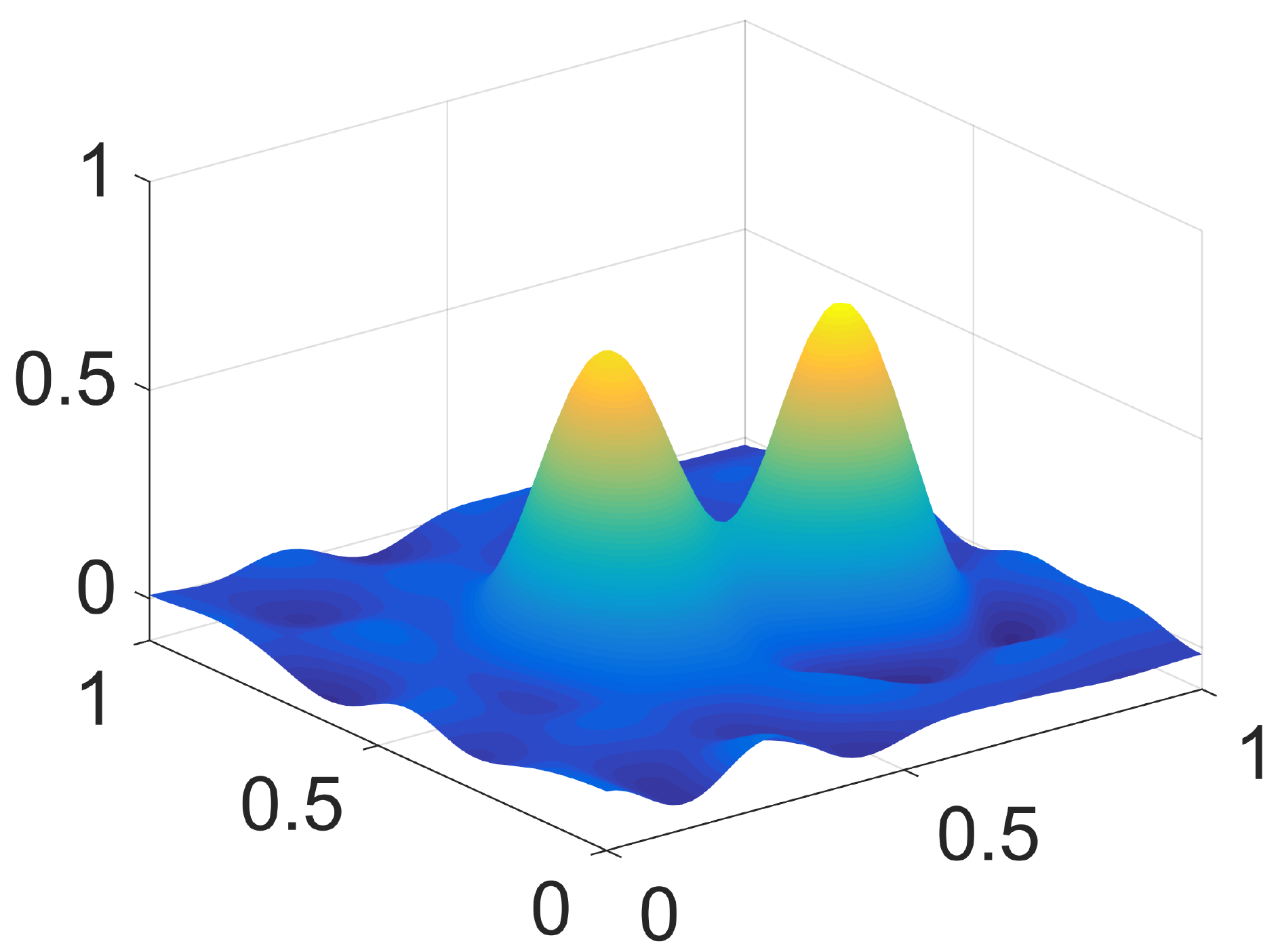} &
\includegraphics[width=3.75cm]{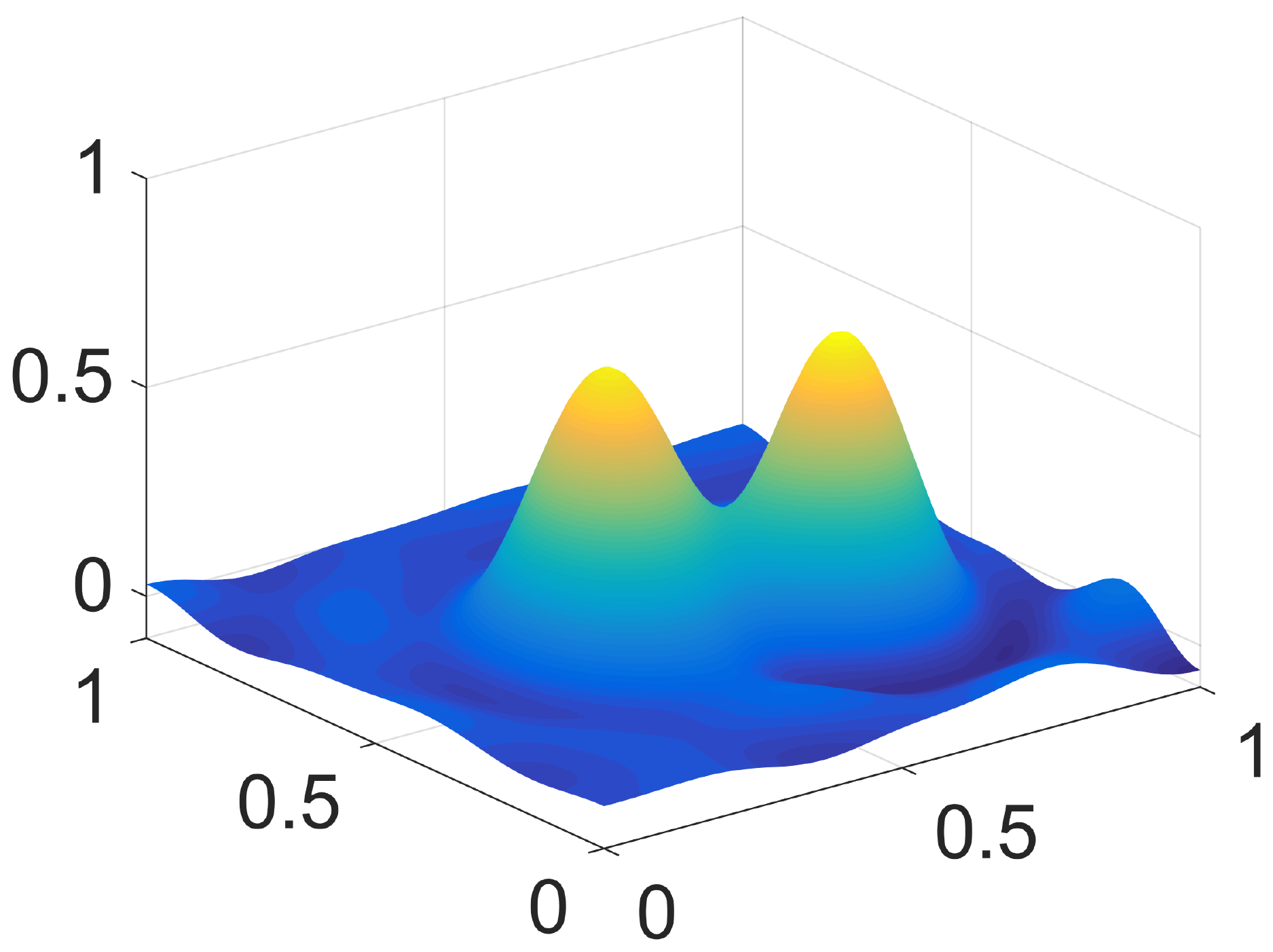}\\
(d) & (e) & (f)
\end{tabular}
\end{center}
\caption{\label{fig:diffusionDemo} Illustration of the solution of the
2D inverse diffusion problem \texttt{PRdiffusion} with \texttt{n} = 64
by means of \texttt{IRrrgmres}.  We set \texttt{options.NoStop = 'on'}
to force the iterations to continue beyond the number of iterations
selected by the discrepancy principle stopping rule.
(a) true solution \texttt{x}, (b) relative error history,
(c) relative residual norm history, (d) noisy data \texttt{bn},
(e) best reconstruction ($k = 35$), (f) reconstruction obtained when
the discrepancy principle is satisfied ($k=24$).}
\end{figure}

This example illustrates the use of \texttt{IRrrgmes} which does not require
operations with the adjoint operator (the matrix transpose).
We consider the 2D inverse diffusion problem
from \S\ref{sec:inversediffusion} in \texttt{PRdiffusion}.
As with previous examples in this section, we begin by setting up the
test problem:
\begin{verbatim}
   n = 64;
   NoiseLevel = 0.005;
   [A, b, x, ProbInfo] = PRdiffusion(n);
   [bn, NoiseInfo] = PRnoise(b, NoiseLevel);
\end{verbatim}
The true solution {\tt x} and data {\tt bn} are shown
in Figure~\ref{fig:diffusionDemo}a and Figure~\ref{fig:diffusionDemo}d, respectively.
We now use RRGMRES to solve the problem, but first we set a few
options {by modifying the appropriate fields in the option structure}:
\begin{itemize}
\item
First, we set {\tt options.x\_true} to the true solution {\tt x} which
allows the method to compute relative error norms.
%$\|x^{(k)} - \mathtt{x} \|_2/\| \mathtt{x} \|_2$, where $x^{(k)}$ is the
%computed approximate solution at iteration $k$.
\item
Next, we want to use the discrepancy principle as stopping rule,
so we need to set {\tt options.NoiseLevel}. We also change the default
safety parameter {\tt eta} to be 1.01.
\item
Finally, we turn on the option {\tt NoStop} so that the iteration will
proceed to the maximum number of iterations, even if a
stopping criterion is satisfied.
\end{itemize}
{As mentioned earlier, }all these parameters can be set in a single call to {\tt IRset},
\begin{verbatim}
   options = IRset('x_true', x, 'NoiseLevel', NoiseLevel, ...
         'eta', 1.01,'NoStop', 'on');
\end{verbatim}
We can then use RRGMRES as follows:
\begin{verbatim}
   [X, IterInfo] = IRrrgmres(A, bn, K, options);
\end{verbatim}
Once the iterations are completed, we can access several pieces of information
from the structure {\tt IterInfo}. Specifically,
\begin{itemize}
\item
{\tt IterInfo.Enrm} contains the {relative} error norms, which are displayed
in Figure~\ref{fig:diffusionDemo}b.
%the first row, second column of Figure~\ref{fig:diffusionDemo}.
\item
The ``best regularized solution" is saved in
{\tt IterInfo.BestReg.X}, and the iteration where the error
is smallest can be found in {\tt IterInfo.BestReg.It}.
This solution is shown in Figure~\ref{fig:diffusionDemo}e.
%This solution, and the corresponding iteration number, are shown
%in the second row, second column of Figure~\ref{fig:diffusionDemo}.
\item
{\tt IterInfo.Rnrm} contains the {relative} residual norms at each
iteration,\linebreak[4]$\|b - A\, x^{(k)}\|_2{/\|b\|_2}$,
which are displayed in Figure~\ref{fig:diffusionDemo}c,
%the first row, third column of Figure~\ref{fig:diffusionDemo},
along with a line marking the stopping point defined by the discrepancy
principle.  That is, once the residual norm reaches the red dashed line,
the convergence criterion defined by the discrepancy principle is considered
satisfied.
\item
The precise iteration {satisfying }the stopping criterion, along with
its corresponding solution, can be obtained from
{\tt IterInfo.StopReg.It} and\linebreak[4]{\tt IterInfo.StopReg.X}, respectively.
This solution is shown in Figure~\ref{fig:diffusionDemo}f.
%This solution, and the corresponding iteration number, are shown
%in the second row, third column of Figure~\ref{fig:diffusionDemo}.
\end{itemize}
We again emphasize that all the iterative methods implemented in our package
have a similar input and output structure.

The code used to generate the test problem and results described in this example is
provided in our package in the script {\tt EXdiffusion\_rrgmres.m}.

\subsection{Solving the 2D NMR Relaxometry Problem with MRNSD}

\begin{figure}
\begin{center}
\begin{tabular}{ccc}
\includegraphics[width=3.75cm]{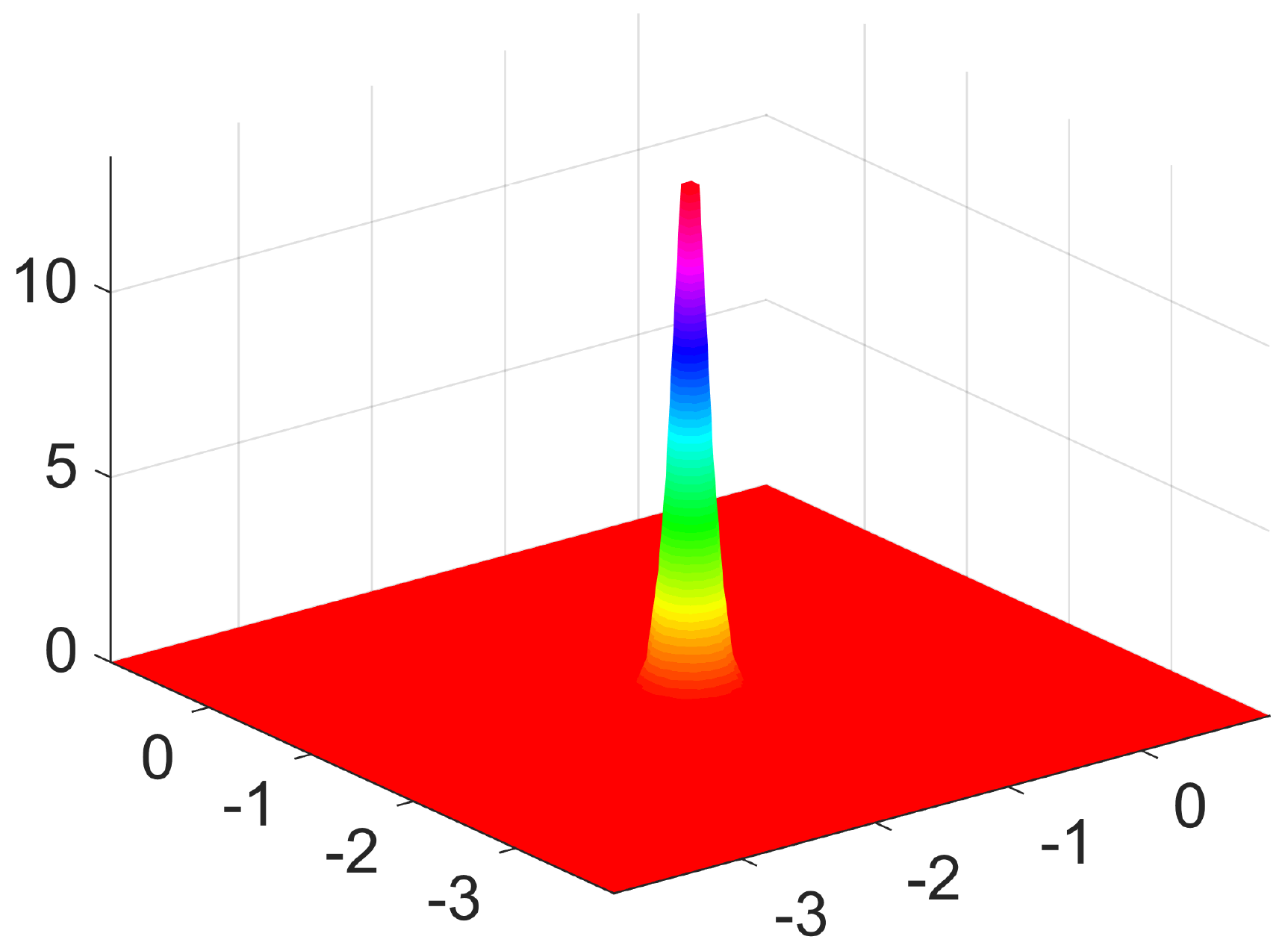} &
\includegraphics[width=3.75cm]{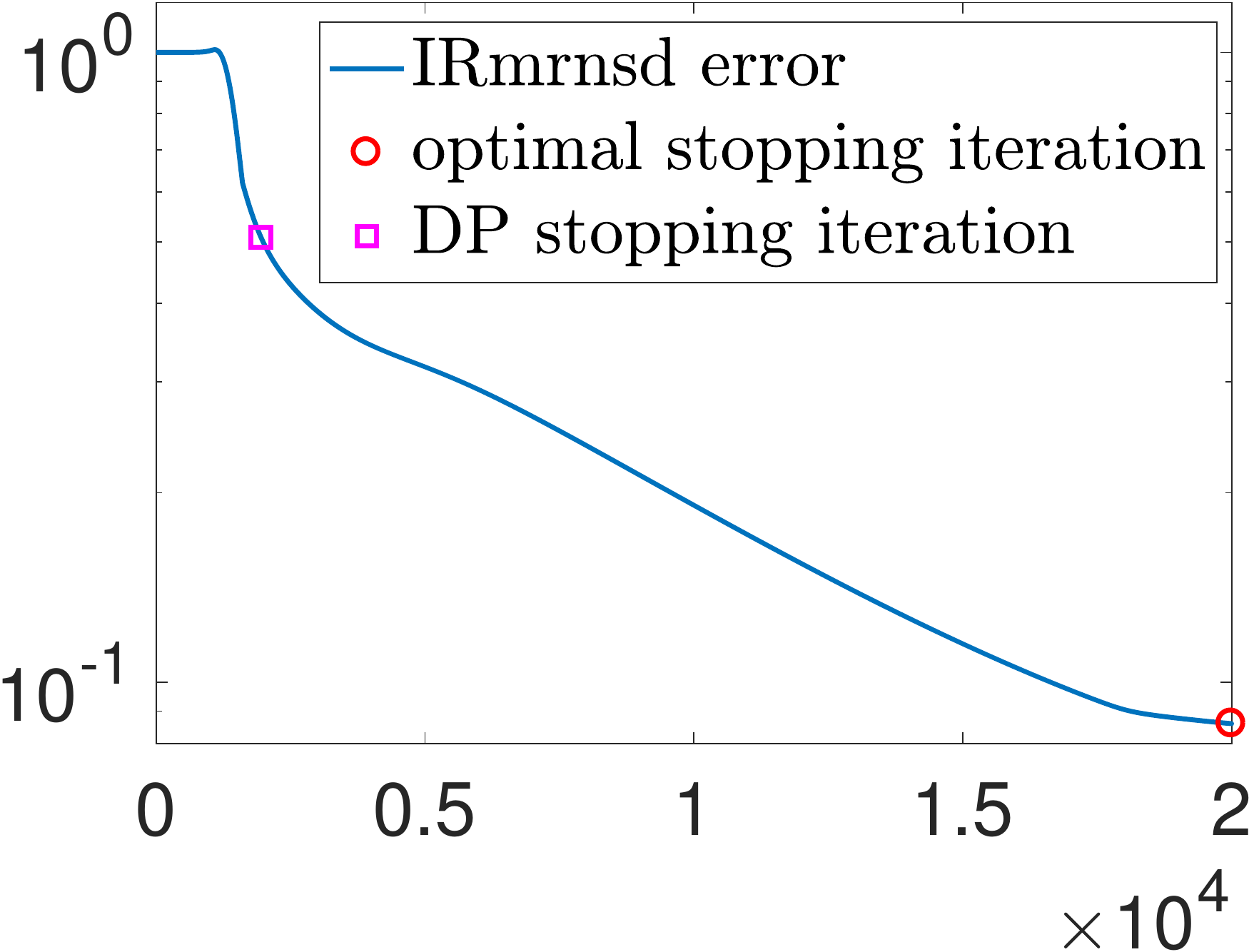} &
\includegraphics[width=3.75cm]{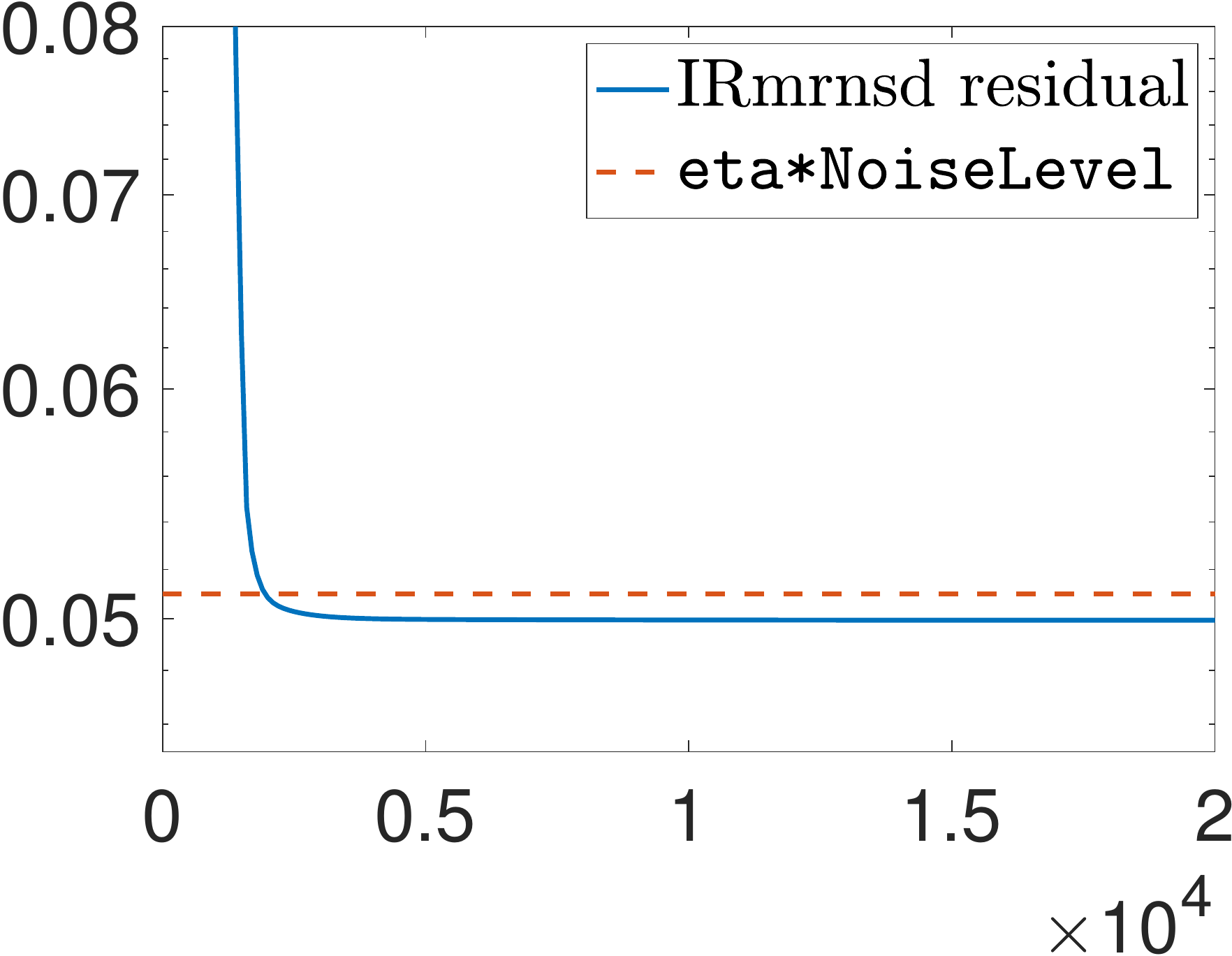}\\
(a) & (b) & (c)\\
\includegraphics[width=3.75cm]{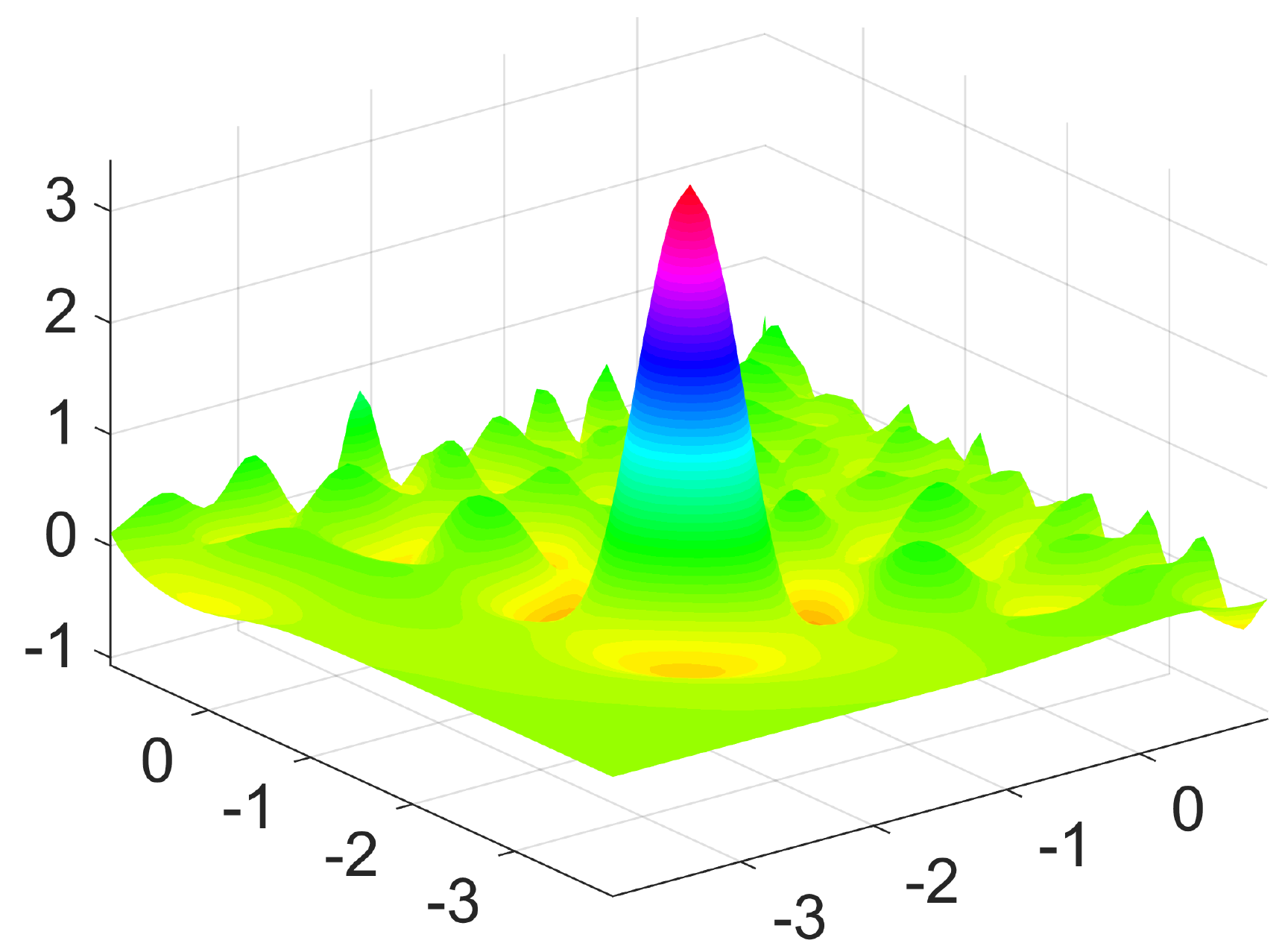} &
\includegraphics[width=3.75cm]{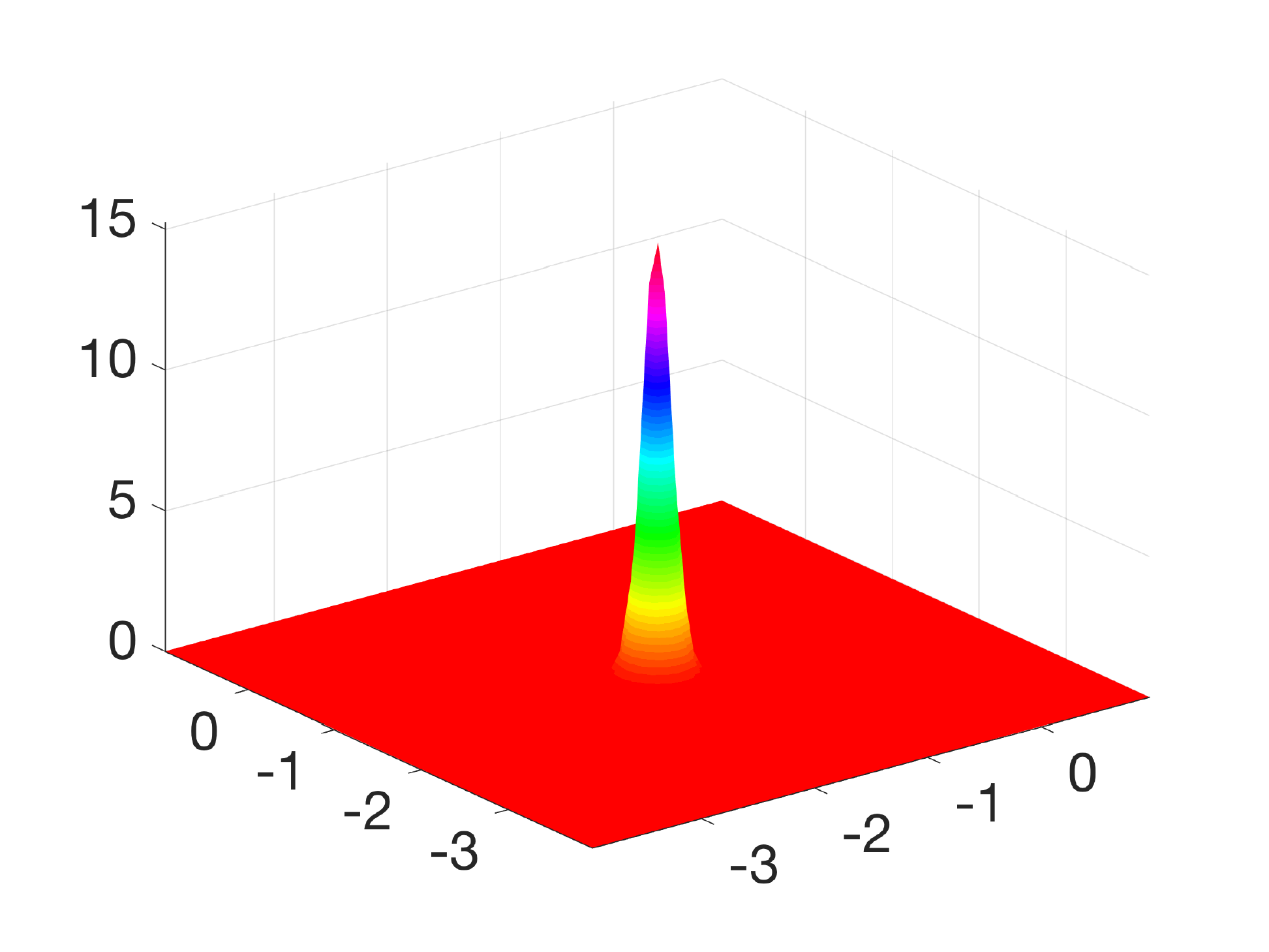} &
\includegraphics[width=3.75cm]{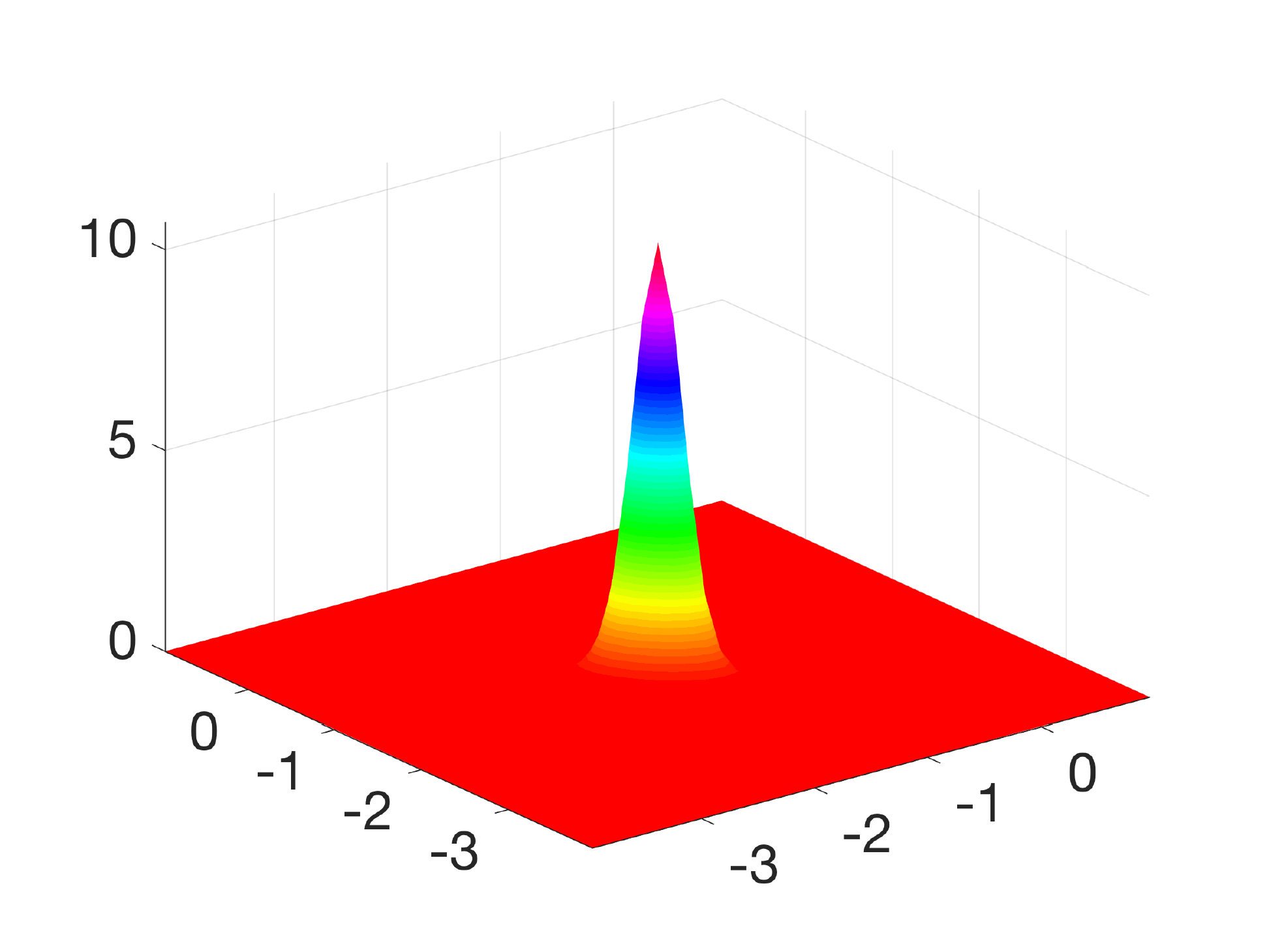}\\
(d) & (e) & (f)
\end{tabular}
\end{center}
\caption{\label{fig:nmrDemo} Illustration of the solution of the
2D NMR relaxometry problem \texttt{PRrmn} with \texttt{n} = 64
by means of \texttt{IRmrnsd}.  We use a color map that emphasizes the
behavior of the large flat region.  We set \texttt{options.NoStop = 'on'}
to force the iterations to continue beyond the number of iterations
selected by the discrepancy principle stopping rule.
(a) true solution \texttt{x}, (b) relative error history,
(c) relative residual norm history, (d) best reconstruction by CGLS
($k = 94$), (e) best reconstruction by MRNSD ($k = 19999$),
(f) reconstruction obtained by MRNSD when the discrepancy principle
is satisfied ($k=1950$).}
\end{figure}

To demonstrate the advantage of imposing nonnegativity constraints we
consider the 2D NRM relaxometry problem from \S\ref{sec:nmr}
implemented in \texttt{PRnmr}, with \texttt{n} = 64.
As done for the other test problems, we begin by setting
\begin{verbatim}
   n = 64;
   NoiseLevel = 0.05;
   [A, b, x, ProbInfo] = PRnmr(n);
   bn = PRnoise(b, NoiseLevel);
\end{verbatim}
The true solution {\tt x} is shown
in Figure~\ref{fig:nmrDemo}a.
This test problem is extremely hard to solve,
and every iterative method available in our package requires a large
amount of iterations to compute a meaningful approximation of~{\tt x}.
We allow 20000 iterations at most, and we can store one approximate
solution every 1000 iterations by setting
\linebreak[4]\texttt{K = [1, 1000:1000:20000]}.
We assign the following options by calling the \texttt{IRset} function:
\begin{verbatim}
   options = IRset('x_true', x, 'NoiseLevel', NoiseLevel, ...
         'eta', 1.01, 'NoStop', 'on');
\end{verbatim}
We now use CGLS as follows:
\begin{verbatim}
   [X_cgls, IterInfo_cgls] = IRcgls(A, bn, K, options);
\end{verbatim}
The solution computed by means of \texttt{IRcgls} is shown in Figure~\ref{fig:nmrDemo}d;
%in the bottom left plot of Figure~\ref{fig:nmrDemo};
this solution hardly resembles the exact one reported in Figure~\ref{fig:nmrDemo}a and,
more specifically, it has large oscillations and negative values in the part
that ideally should be zero.

To run \texttt{IRmrnsd} with the same test data and input options as the ones used for \texttt{IRcgls} we simply type
\begin{verbatim}
   [X_mrnsd, IterInfo_mrnsd] = IRmrnsd(A, bn, K, options);
\end{verbatim}
We recall that nonnegativity constraints are automatically imposed
within the \texttt{IRmrnsd} iterations.
{\tt IterInfo\_mrnsd} stores various pieces of information about the
behavior of this solver applied to this test problem.
In particular, we can access the relative error at each iteration
in {\tt IterInfo.Enrm}, which is displayed in Figure~\ref{fig:nmrDemo}b.
The relative residual at each iteration is stored in {\tt IterInfo.Rnrm};
this is displayed in Figure~\ref{fig:nmrDemo}c, together with a
horizontal line marking the relative noise level
(useful to visually inspect when the discrepancy principle is satisfied).
The ``best regularized solution" is saved in
{\tt IterInfo\_mrnsd.BestReg.X}, and the iteration where the error
is smallest can be found in {\tt IterInfo\_mrnsd.BestReg.It}.
This solution is shown in Figure~\ref{fig:nmrDemo}e.
The precise iteration {satisfying }the discrepancy principle, along with
its corresponding solution, can be obtained from
{\tt IterInfo\_mrnsd.StopReg.It} and {\tt IterInfo\_mrnsd.StopReg.X}, respectively.
This solution is shown in Figure~\ref{fig:nmrDemo}f.

The code used to generate the test problem and results described in this example is
provided in our package in the script {\tt EXnmr\_cgls\_mrnsd.m}.

\subsection{Computing Sparse Reconstructions}

This example illustrates how to use \texttt{IRirn} and \texttt{IRell1} to compute
\textit{approximately sparse} reconstructions -- in the sense that
the solution has many small values (that may consecutively be
truncated to zero).

The test problem is Gaussian image deblurring, and we choose
one of our synthetically generated images that is made up of
randomly placed small ``dots", with random intensities.
This test mage may be used, for example, to simulate stars
being imaged from ground based telescopes.
To generate the test problem, we use the options structure to specify
the  {\tt 'dotk'} synthetic image,
\begin{verbatim}
   PRoptions.trueImage = 'dotk';
   [A, b, x, ProbInfo] = PRblurgauss(PRoptions);
\end{verbatim}
and add $10\%$ white noise,
\begin{verbatim}
   NoiseLevel = 0.1;
   bn = PRnoise(b, NoiseLevel);
\end{verbatim}
We then compute the best solutions by means of \texttt{IRcgls}, which
cannot impose sparsity, as well as \texttt{IRell1} and \texttt{IRirn}, which are simplified drivers for
\texttt{IRhybrid\_fgmres} and \texttt{IRrestart}, respectively.
%which is a simplified driver for \texttt{IRhybrid\_fgmres}, and \texttt{IRirn}, which is a simplified driver for \texttt{IRrestart}.
Both \texttt{IRell1} and \texttt{IRirn} are designed to
make it easy to approximately enforce a 1-norm penalization on the solution,
leading to a reconstruction with many small elements.
%For all the solvers we use the ``no stop'' feature to continue the iterations
%beyond the iteration found by the default stopping rule, such that we
%can display the best solutions.

To illustrate the effect of the parameter-choice rule for the
projected problem in the hybrid method \texttt{IRhybrid\_fgmres}, we use
both GCV (which is the default) and the discrepancy principle. If GCV is used, then the iterations stop when the minimum of the iteration-dependent GCV function stabilizes or starts increasing within a given window. If the discrepancy principle is used, the iterations are stopped according to the strategy proposed in \cite{GaNo14}, and previously addressed in Section \ref{sec:StopR}. The regularization parameters for the inner iterations of \texttt{IRirn} are chosen by the discrepancy principle that
%, similarly to the \texttt{IRell1} case,
also acts as a stopping rule for the inner iterations; the default stopping criterion for the outer iterations is the stabilization of the norm of the solution at each restart. However, for all the solvers, we are interested in computing the best solutions, which may be found after the stopping rules are satisfied. For this reason we use the ``no stop'' feature in \texttt{IRcgls} and \texttt{IRell1}, and the ``no stop out'' feature in \texttt{IRirn}, so to ensure that the iterations are continued after the stopping criterion (for \texttt{IRcgls} and \texttt{IRell1}) and the outer stopping criterion (for \texttt{IRirn}) are satisfied.

Specifically, first
run {\tt IRcgls} for 80 iterations, using the true solution to compute error norms,
and turn {\tt NoStop} on:
\begin{verbatim}
   options = IRset('MaxIter', 80, 'x_true', x, 'NoStop', 'on');
   [Xcgls, info_cgls] = IRcgls(A, bn, options);
\end{verbatim}
Now compute a sparse solution using the default GCV rule for choosing the
regularization parameter of the projected problem,
\begin{verbatim}
   options.NoStop = 'on';
   [Xell1_GCV, info_ell1_GCV] = IRell1(A, bn, options);
\end{verbatim}
To change the default regularization parameter-choice rule
to the discrepancy principle, using the true {\tt NoiseLevel} with a safety
value for {\tt eta}, we use \texttt{IRset} as  follows:
\begin{verbatim}
   options = IRset(options, 'RegParam', 'discrep', ...
                            'NoiseLevel', NoiseLevel, 'eta', 1.1);
   [Xell1_DP, info_ell1_DP] = IRell1(A, bn, options);
\end{verbatim}
Finally, consider \texttt{IRirn} with the discrepancy principle used for the inner iterations, with \texttt{NoStopOut} turned on, and with 80 total iterations. This can be simply achieved as follows:
\begin{verbatim}
   options.NoStopOut = 'on';
   K = 80;
   [Xirn_DP,info_irn_DP] = IRirn(A, bn, K, options);
\end{verbatim}
Because of the interplay between the inner iterations (whose number depends on the discrepancy principle) and outer iterations, we need to explicitly specify \texttt{K} to ensure that the maximum number of total iterations is 80.

\begin{figure}[htbp]
\begin{center}
\begin{tabular}{ccc}
\includegraphics[width=3.5cm]{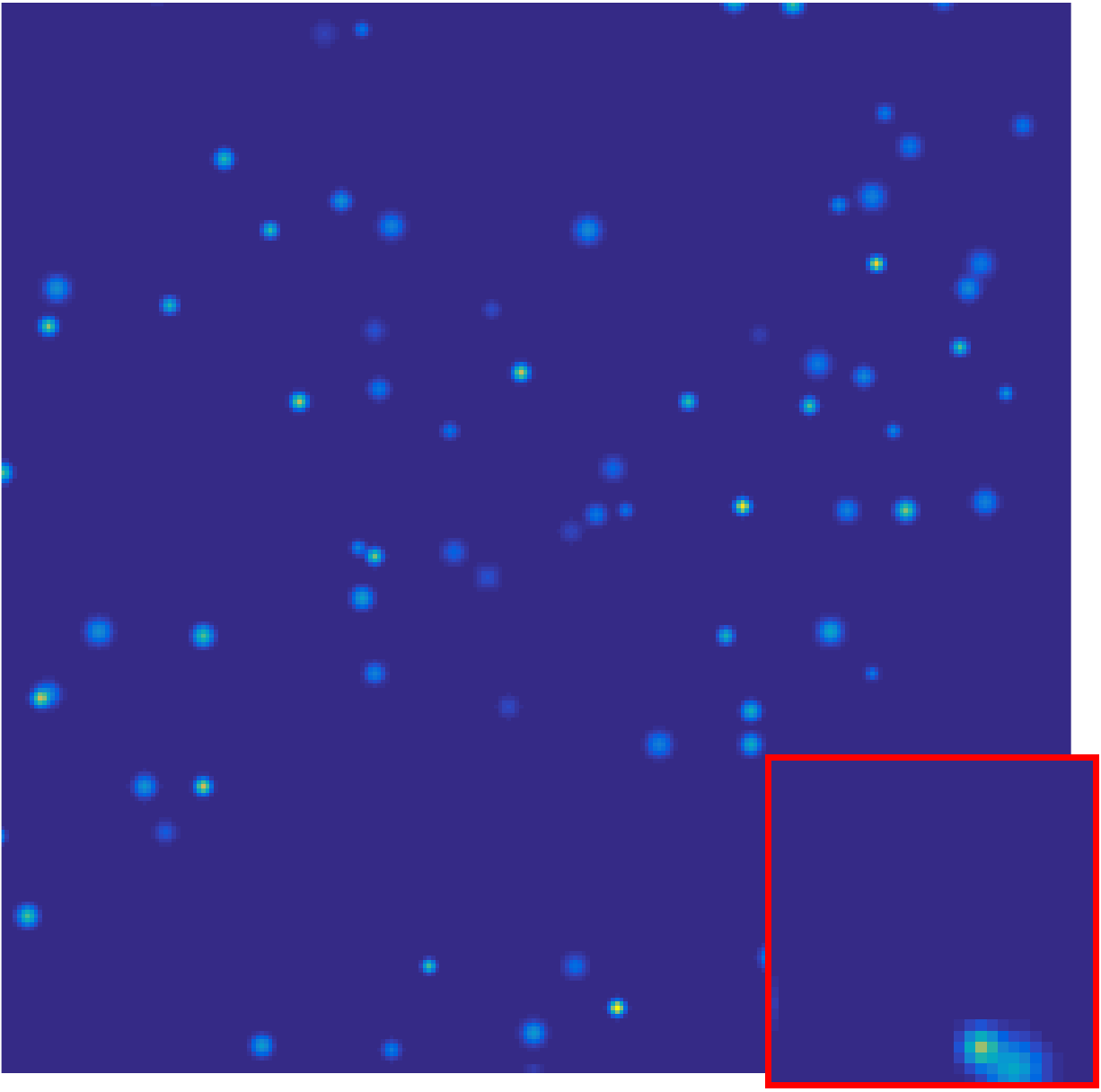} &
\includegraphics[width=3.5cm]{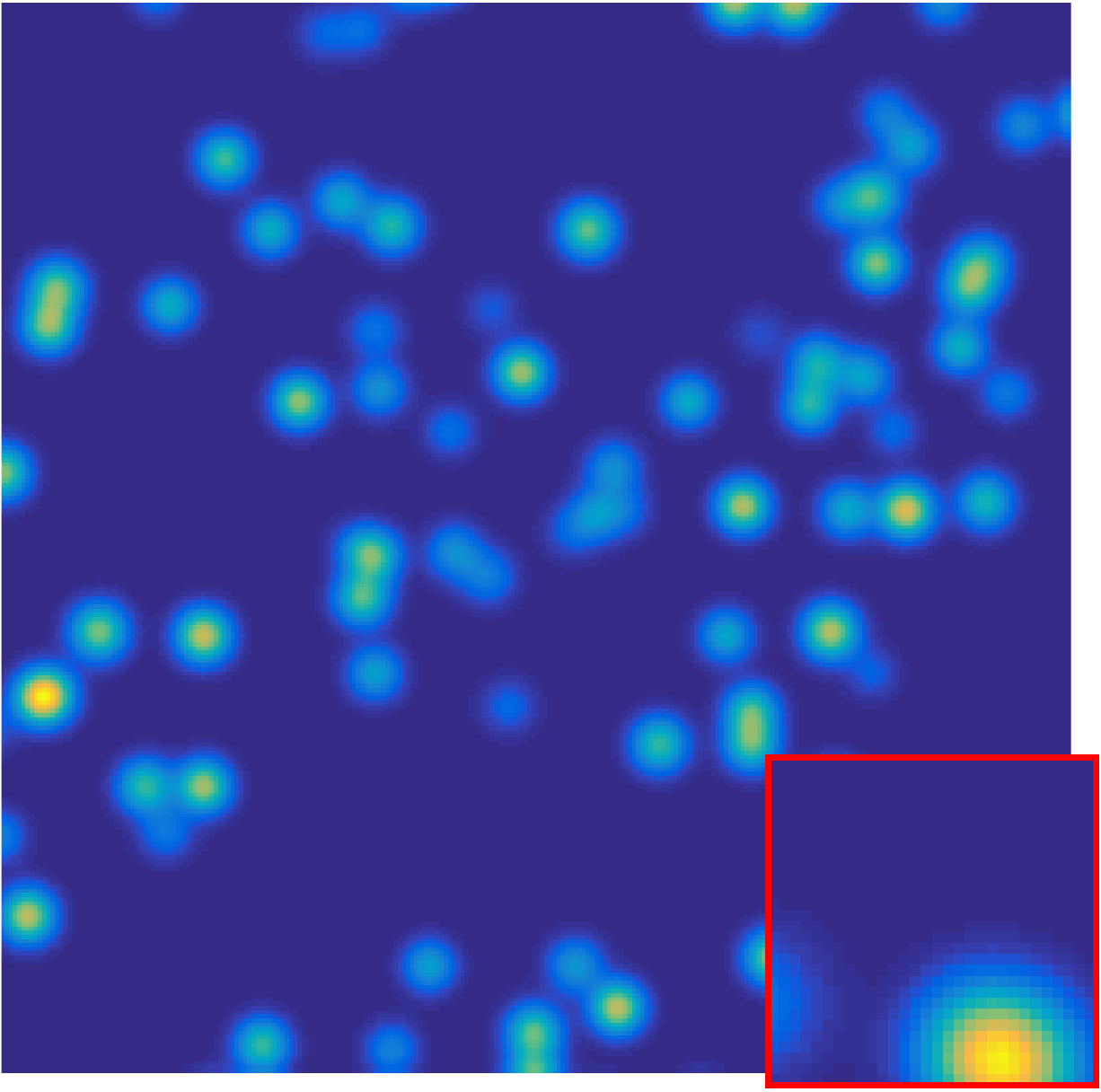} &
\includegraphics[width=3.5cm]{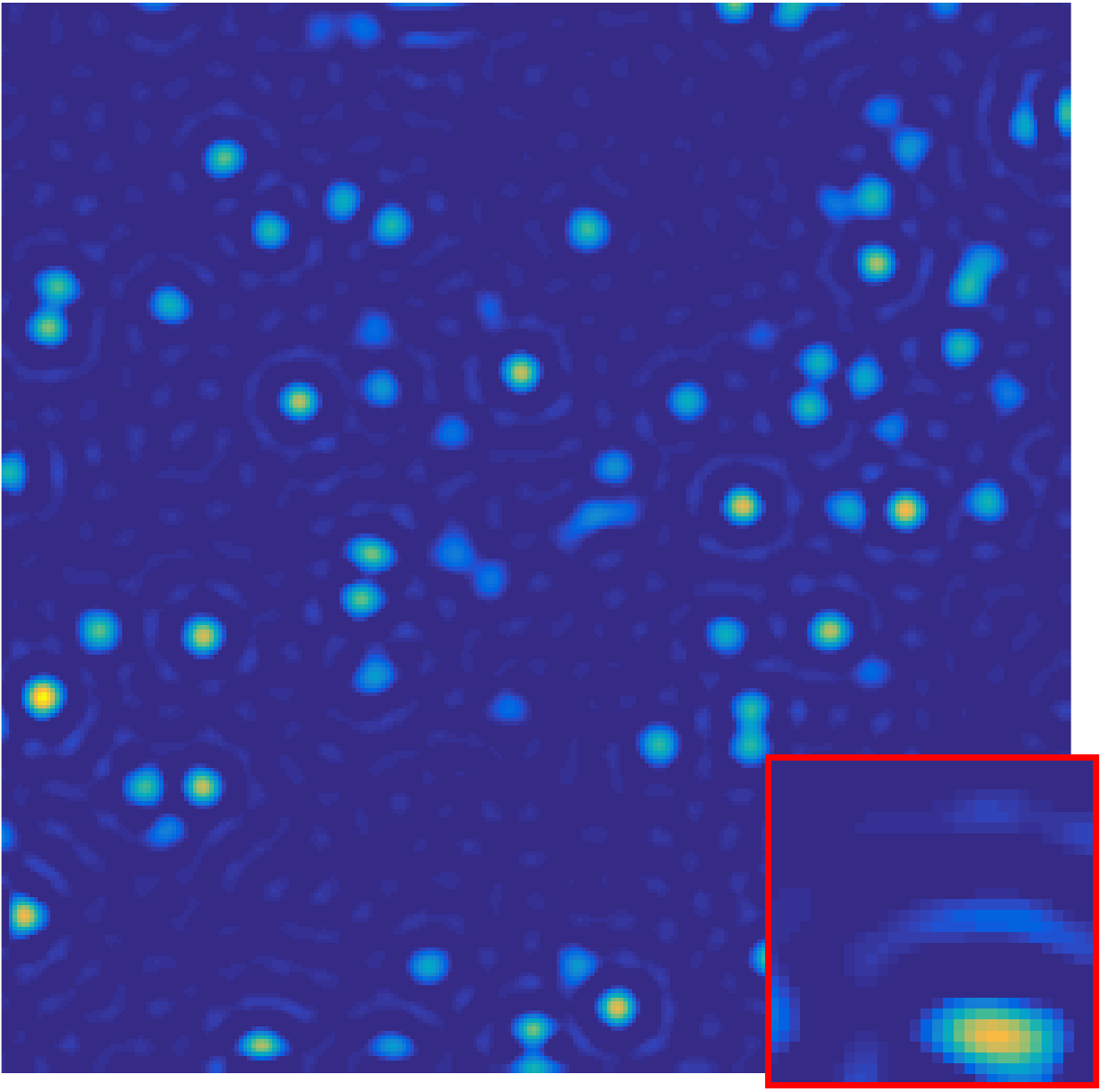}  \\
(a) & (b) & (c) \\
\includegraphics[width=3.5cm]{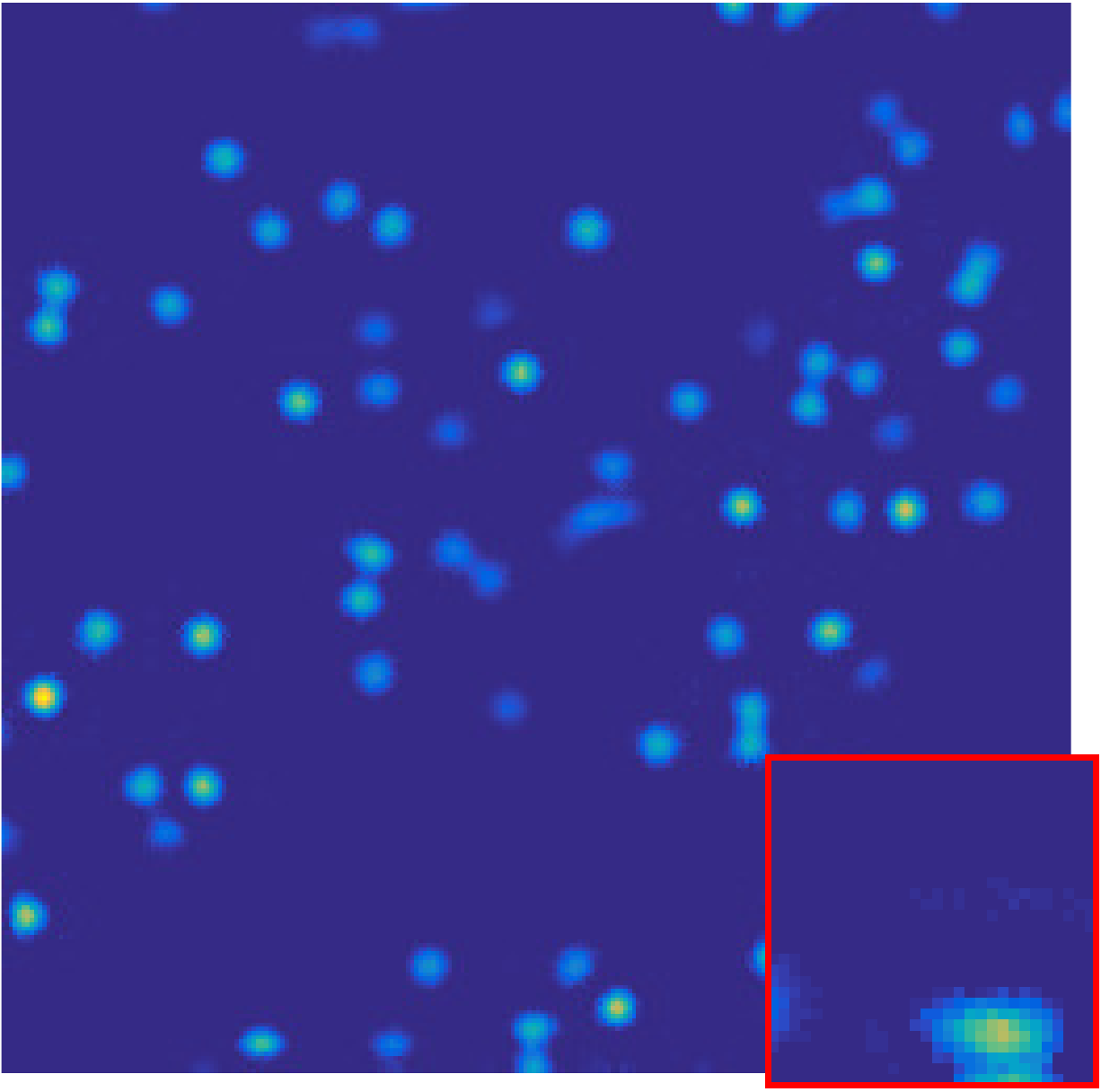} &
\includegraphics[width=3.5cm]{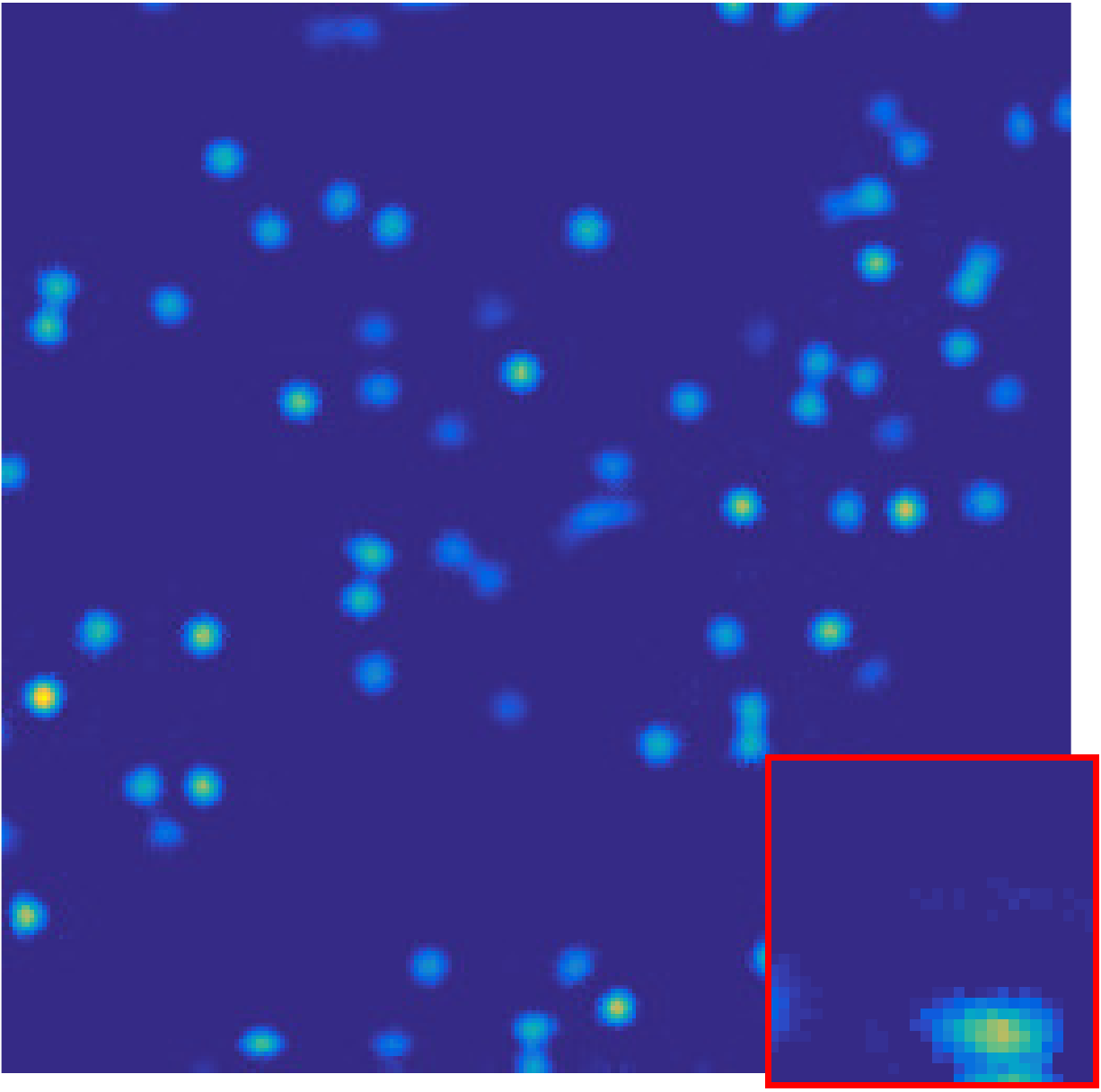} &
\includegraphics[width=3.5cm]{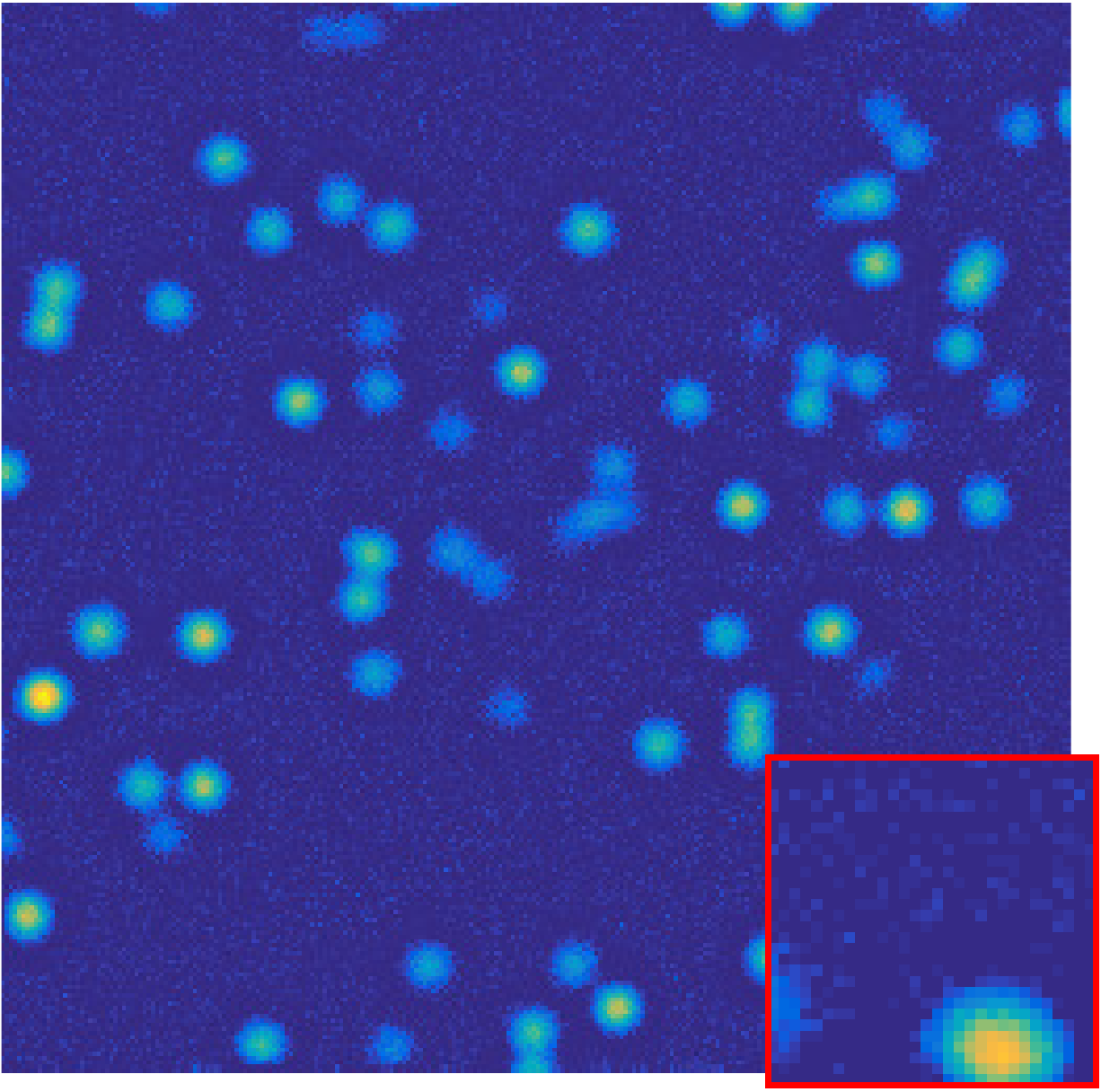}  \\
(d) & (e) & (f)
\end{tabular}
\caption{Illustration of sparse solutions to an image deblurring problem
generated with \texttt{PRblur} and \texttt{n} $=256$, and with a sparse test image
with synthetic stars.
(a) true solution \texttt{x}, (b) noisy blurred image \texttt{bn},
(c) best CGLS reconstruction ($k = 53$),
(d) best \texttt{IRell1} solution with the default GCV parameter-choice rule
for the projected problem ($k = 5$),
(e) best \texttt{IRell1} solution with the discrepancy principle
parameter-choice rule for the projected problem ($k = 51$),
(f) best \texttt{IRirn} solution with the discrepancy principle parameter-choice and stopping rule for the inner iterations ($k=4$).
All negative pixels are truncated to 0
and the inset figures zoom in on the bottom right $30\times 30$ corner of the image.}
\label{fig:EXsparity}
\end{center}
\end{figure}

Figures \ref{fig:EXsparity}a and \ref{fig:EXsparity}b show the true image and the
noisy blurred image, respectively.
The CGLS reconstruction shown in \ref{fig:EXsparity}c is clearly contaminated
by oscillations (``ringing effects'') around the reconstructed stars, and we
also see other artifacts that appear as ``freckles'' as discussed in \cite{HaJe08}.
For the reconstructions computed by \texttt{IRell1} shown in \ref{fig:EXsparity}d-e,
we see that the parameter-choice rule for the projected problem
indeed has an effect on the iterations -- in this example we obtain the best
reconstruction with the discrepancy principle. Also the reconstruction computed by \texttt{IRirn} shown in \ref{fig:EXsparity}f is successful in computing a sparse solution though, for this test problem, \texttt{IRell1} exhibits a better performance. This example demonstrates that the heuristic approach to computing a
sparse reconstruction in \texttt{IRell1} works well -- as long as one
can accept small elements rather than exact zeros in the reconstruction.

The code used to generate the test problem and results described in this
example is provided in our package in the script {\tt EXsparsity.m}.

\section{Conclusion}
\label{sec:Conclusions}

We gave an overview of a MATLAB software package \textsc{IR Tools}
that provides large-scale iterative regularization methods and
new large-scale test problems.
Our package allows the user to easily experiment with a variety of
well-documented iterative regularization methods in a flexible and
uniform framework, and at the same time our software can be used
efficiently for real-data problems.
We also provide a set of realistic large-scale 2D test problems
that replace the outdated ones from \textsc{Regularization Tools}
and that are valid alternatives to the ones available within other popular
MATLAB toolboxes and packages.

\section{Acknowledgements}

The authors are grateful to Julianne Chung for providing an implementation
of \textsc{HyBR}, which forms the basis of our \texttt{IRhybrid\_lsqr} function.
For further details, see \cite{ChungThesis,HyBR} and
{\color{blue}\url{http://www.math.vt.edu/people/jmchung/hybr.html}}.
We also thank Germana Landi for providing insight about the NMR relaxometry problem.

The satellite image in our package, shown in Fig.~\ref{fig:SpeckleBlurExample},
is a test problem that originated from
the US Air Force Phillips Laboratory, Lasers and Imaging Directorate,
Kirtland Air Force Base, New Mexico.  The image is from a computer simulation
of a field experiment showing a satellite  as taken from a ground
based telescope.  This data has been used widely in the literature
for testing algorithms for ill-posed image restoration problems; see,
for example \cite{RoWe96}.

Our package also includes a picture of NASA's Hubble Space Telescope as
shown in Fig,~\ref{fig:EXblur_cgls_images}.
The picture is in the public domain and can be obtained from
{\color{blue}\url{https://www.nasa.gov/mission_pages/hubble/story/index.html}}.

\end{document}